\documentclass{article}
\title{Dissipativity in infinite horizon optimal control and dynamic programming}
\author{David Angeli\footnote{David Angeli is with the Dept.\ of Electrical and Electronic Engineering of Imperial College London and with the Dept. of Information Engineering, University of Florence, Italy.}\; and Lars Gr\"une\footnote{Lars Gr\"une is with the Dept.\ of Mathematics of University of Bayreuth, Germany.}}
 
\usepackage{amssymb,amsmath,graphicx}
\usepackage{appendix}

\numberwithin{equation}{section}
\numberwithin{figure}{section}

\usepackage{color}
\newcommand{\LG}{\textcolor{blue}}

\newcommand{\comment}[1]{}

\newtheorem{thm}{\textbf{Theorem}}[section]
\newtheorem{assum}[thm]{\textbf{Assumption}}
\newtheorem{defn}[thm]{\textbf{Definition}}
\newtheorem{lem}[thm]{\textbf{Lemma}}
\newtheorem{prop}[thm]{\textbf{Proposition}}
\newtheorem{rem}[thm]{\textbf{Remark}}
\newtheorem{cor}[thm]{\textbf{Corollary}}

\def\qed{
\vbox{\hrule\hbox{\vrule\hbox to 5pt{\vbox to 8pt{\vfil}\hfil}\vrule}\hrule}}
\def\endproof{\unskip \nobreak \hskip0pt plus 1fill \qquad \qed\newline}
\def\wT{\widetilde T}

\begin{document}
\maketitle
\begin{abstract}
In this paper we extend dynamic programming techniques to the study of discrete-time infinite horizon optimal control problems on compact control invariant sets with state-independent best asymptotic average cost. To this end we analyse the interplay of dissipativity and optimal control, and propose novel recursive approaches for the solution of so called shifted Bellman Equations.  
\end{abstract}
\section{Introduction}
Dynamic programming (DP) is a cornerstone of control theory which allows to solve (in feedback form) optimal control problems formulated on horizons of increasing length through a suitable recursive formula for the computation of the so called value function, \cite{bellman}. 

Remarkably, dynamic programming allows to study problems formulated both on a finite horizon or an infinite one, the latter achieved under suitable technical assumptions, by studying the asymptotic properties of the recursion or by computing its fixed points. 
By now, the subject of dynamic programming and infinite horizon optimal control has been studied in depth by many authors and several monographs on the subject exist both in the control domain   \cite{bertsekas1,bertsekas3,infinite} and in economics, \cite{stokey,economic2}. 

While, in its naive form, DP is often associated to the curse of dimensionality, which may hinder its applicability to scenarios of practical relevance, the topic of its approximate and efficient numerical treatment has also gathered significant impetus, in particular in the context of machine learning, \cite{bertsekas2}. Indeed, the dynamic programming or Bellman Equation is at the core of any (deep) reinforcement learning algorithm \cite{BarS18,ADBB17}. 

The link of optimal control to dissipativity was already established by Willems in the seminal papers 
\cite{willems1,willems2} and in parallel in the 
study of nonlinear inverse optimal regulators for nonlinear systems, \cite{moylan}. 
However, it was only brought to the forefront of the discourse on optimisation-based control in recent years, \cite{dissipativity3,dissipativity4}, thanks to its surprising connections to closed-loop stability of
Economic Model Predictive Control \cite{dissipativity1,dissipativity2} and long-run average optimal control, \cite{finlay,gaitsgory0}. 
In particular, \cite{dissipativity1} proposes a notion of optimal operation at steady-state and provides a sufficient conditions for this property to hold based on dissipativity of the associated systems' dynamics with respect to a suitable supply function. The converse statement is investigated in \cite{dissipativity2} where an additional controllability assumption is needed in order to prove necessity of dissipativity. 
While generalizations of the above results, its relation to the so-called turnpike property, and extensions to periodic optimal solutions are provided in several subsequent works (i.e.  \cite{MULLERgrune} and \cite{muller}), the connection to Dynamic Programming and infinite horizon optimal control has remained elusive, due to restrictive technical assumptions needed to make sense of undiscounted cost functionals.

In this paper we further explore connections between dissipativity and infinite horizon optimal control problems, while proposing new formulations and iterative methods for their solutions that significantly expand the class of problems which can be meaningfully addressed by this approach.
Our main contributions are
\begin{itemize}
    \item introducing a terminal penalty in infinite horizon optimal control, in the form of suitable \emph{storage functions} with negative sign;
    \item proposing a shifted Bellman Equation to be used in optimal control problems with non-zero (yet state-independent) optimal long run average performance (this includes systems with periodic, almost periodic or even chaotic regimes of operation allowing general time-varying asymptotic cost along optimal solutions);
    \item proposing two novel recursions whose fixed points are solutions of shifted Bellman Equation (of any shift);
    \item analysing the convergence properties of such recursions under fairly general technical assumptions, allowing simultaneous computation of the best average performance and of the associated value function;
    \item tackling the non-existing trade-off between transient cost and asymptotic average performance.
\end{itemize}
The rest of the paper is organized as follows: Section \ref{preliminary} introduces the problem formulation, basic notations and some preliminary results, 
Section \ref{shiftedsection} introduces the shifted Bellman Equation and the novel recursion operators whose properties are investigated in Section \ref{propertiesoperators}. 
Section \ref{equicont} provides a general convergence result under suitable conditions on the controllability of the system's dynamics, while Section \ref{additionlars} relaxes some continuity assumptions needed for convergence analysis approaching the recursion from specific initialisations. Examples and counter-examples are shown in Section \ref{exandcex} and Section \ref{finalsection} draws some conclusions and points to further open research directions. Important intermediate technical results are collected in the appendix in Section \ref{preliminaryresults}.

\section{Problem formulation and preliminary results}
\label{preliminary}
Consider the discrete-time finite dimensional nonlinear control system described by the following difference equations:
\begin{equation}
\label{systemmap}
x (t+1) = f(x(t),u(t))
\end{equation}
where $x(t) \in \mathbb{X} \subset \mathbb{R}^n$ is the state-variable, taking values in some compact control invariant set $\mathbb{X}$, $u(t) \in \mathbb{R}^m$ is the control input and
$f: \mathbb{Z} \rightarrow \mathbb{X}$, is the continuous transition map. 
We denote by $\mathbb{U}(\cdot): \mathbb{X} \rightarrow 2^{\mathbb{R}^m} $ the upper semicontinuous set-valued mapping defined below:
\begin{equation}
\label{feasiblecontrol}
\mathbb{U} (x) := \{ u \in \mathbb{R}^m: (x,u) \in \mathbb{Z} \},
\end{equation}
 which corresponds to the set of feasible control inputs in state $x$, given the compact state/input constraint
 set $\mathbb{Z}$.
 Moreover, we assume, without loss of generality,
 \begin{equation}
     \label{controlinvariance}
     f(x, \mathbb{U}(x) ) \subset \mathbb{X},
     \end{equation}
for all $x \in \mathbb{X}$.
 For an input sequence $\textbf{u}= \{ u(t)\}_{t=0}^{\infty}$, we denote by $\phi(t,x,\textbf{u})$ the state at time $t$, from initial condition $x(0)=x$, as given by iteration (\ref{systemmap}).
 We also extend definition (\ref{feasiblecontrol}), to allow feasible control sequences of length $T$, as follows:
 \begin{equation}
\mathbb{U}_T (x) := \{ \textbf{u}=\{ u(t) \}_{t=0}^{T-1} \in \mathbb{R}^{m T}: 
  (\phi(t,x,\textbf{u}), u(t)) \in \mathbb{Z}, \forall \, t \in \{0, \ldots, T-1\} \}.
\end{equation}
Our contribution is twofold; namely, to define optimal control problems over an infinite horizon within a significantly larger set of systems dynamics and associated cost functional than is currently possible to address by existing formulations, and, at the same time, to propose a dynamic programming approach for their solution.
To this end we consider a continuous stage cost, $\ell (x,u): \mathbb{Z} \rightarrow \mathbb{R}$, and formulate the following cost functional:
\begin{equation}
J^{\psi}_T (x(\cdot),u(\cdot) ) : = \sum_{t =0}^{T-1} \ell (x(t),u(t)) + \psi (x(T))
\end{equation}
where $\psi: \mathbb{X} \rightarrow \mathbb{R}$ is a continuous function called the terminal cost. Terminal costs significantly affect the solution of an optimal control problem and a key insight of our paper will be providing guidelines for their selection in order to allow the formulation of infinite horizon optimal control problems.
A finite horizon optimal control problem is then defined as follows:
\begin{equation}
\label{finite}
\begin{array}{rl}
V^{\psi}_{T} (x) : = \min_{x(\cdot), u( \cdot)} &  J_T (x(\cdot),u(\cdot)) \\
\textrm{subject to} & \\
x(0) & = x \\
x(t+1) & = f(x(t),u(t)) \qquad t \in \{0,1, \ldots, T-1 \} \\
 (x(t),u(t))& \in \mathbb{Z} \qquad  t \in \{0,1, \ldots, T-1 \} \\
x(T) & \in \mathbb{X}
\end{array}
\end{equation}
For each value of the initial condition $x \in \mathbb{X}$, a solution of (\ref{finite}) is guaranteed to exist thanks to the compactness and non-emptiness properties of the feasible set and continuity of the cost function.

On the other hand,
when the control problem has no natural termination time, one might want to define an infinite horizon optimisation problem. This has often the additional appealing feature of being achieved through implementation of a time-invariant feedback policy.
However, making sense of an infinite horizon formulation of (\ref{finite}) typically entails strong assumptions on the kind of system's dynamics and cost functional that are allowed.  

One strategy for avoiding such kind of limitations is, at least in practice, to introduce a discounting factor $0< \gamma < 1$ in the cost function:
\begin{equation}
J_\gamma (x(\cdot),u(\cdot) ) : = \sum_{t =0}^{\infty} \gamma^t \ell (x(t),u(t)), 
\end{equation}
which for $\gamma \approx 1$ provides a good approximation to some form of infinite horizon (average) cost. While this approach has some appealing features, for instance making optimal solutions invariant with respect to translation of $\ell$ by any finite constant value, having to settle on a specific value of $\gamma$ less than unity is unsatisfactory as it always leaves open the question of how optimal control policies would be affected by variations in $\gamma$, i. e. if higher values were to be considered.
Moreover, as shown later in Section \ref{motivatingexample}, adoption of a discounting factor may introduce non-existent trade-offs between optimisation of steady-state and transient costs. 

An alternative approach is to resort to average, rather than summed costs:
\begin{equation}
J^{\textrm{avg}} (x(\cdot),u(\cdot) ) : = \limsup_{T \rightarrow + \infty} \frac{ \sum_{t =0}^{T-1} \ell (x(t),u(t)) }{T}.  
\end{equation}
Taking the average yields well-defined costs even when summed costs would be divergent to $\pm \infty$, or are non-convergent (for instance oscillating), which constitute the main obstructions in the definition of infinite horizon control problems for general dynamics and costs. On the other hand, time-shift invariance of average costs along any solution, implies that this approach disregards \emph{transient costs}, which therefore won't be minimised and might be arbitrarily large even for optimal feedback policies (see again example in Section
\ref{motivatingexample}. 

Our proposed solution and novel contribution is to provide fairly general conditions on the terminal cost $\psi$ to make sure that the functional:
\[      V^{\psi}_{\infty} (x) := \lim_{T \rightarrow + \infty} V^{\psi}_{T} (x)   \]
is well-defined. To this end the notion of dissipativity will play an interesting role. This notion was originally introduced by Willems in \cite{willems1,willems2} and has recently received a surge in interest for its crucial role in the analysis of closed-loop Economic Model Predictive Control schemes \cite{dissipativity1,dissipativity2,dissipativity3,dissipativity4}.
In a nutshell a system as (\ref{systemmap}), is said to be dissipative with respect to the supply function $\ell(x,u)$, if there exists a continuous storage function 
$\lambda: \mathbb{X} \rightarrow \mathbb{R}$ such that:
\begin{equation}
\lambda ( f(x,u) ) \leq \lambda (x) + \ell(x,u) \qquad \forall \, (x,u) \in \mathbb{Z}. \label{eq:diss}
\end{equation}
This inequality is normally interpreted in ``energetic'' terms as enforcing, for a dissipative system, that energy stored within, at the next state, cannot
exceed the energy at the current state plus the energy externally supplied through the supply function $\ell(x,u)$.
In the context of optimal control, where the objective is to minimize a cost functional, $\lambda(x)$ can be interpreted as the value of the state $x$ and the dissipation inequality guarantees that the gain in value for any feasible control action $u$ and state $x$ cannot exceed the corresponding stage cost.
Notice that, while optimal control sequences over any finite control horizon (or over infinite control horizon with forgetting factor $\gamma$) are invariant with respect to cost translation, viz. $\tilde{ \ell } (x,u ) := \ell (x,u) - c$ for any constant $c \in \mathbb{R}$, dissipativity is not a shift-invariant property.
In fact, it can always be guaranteed by a sufficiently negative value of $c$, given compactness of $\mathbb{Z}$. Trivially,
if $\tilde{\ell} (x,u) \geq 0$ for all $(x,u) \in \mathbb{X}$, dissipativity is ensured just by defining $\lambda(x)=0$ for all $x \in \mathbb{X}$.
Our first result is stated below.
\begin{prop}
\label{firstprop}
Assume that system (\ref{systemmap}) is dissipative with continuous storage function $\lambda(\cdot)$ with respect to the supply $\ell(x,u)$, and let
$\psi(x) = - \lambda(x)$, then the limit:
\begin{equation}
\label{limitcost}
\lim_{T \rightarrow + \infty} V^{\psi}_T (x)
\end{equation}  
exists for all $x \in \mathbb{X}$, possibly assuming the value $+\infty$.
\end{prop}
\emph{Proof.}
Consider any feasible solution $x^*_{T+1}$, $u^*_{T+1}$  (with $x^*(0) = x$) which achieves the optimal cost $V^{\psi}_{T+1} (x)$.
By definition,
\[ V^{\psi}_{T+1} (x) = \sum_{t = 0}^{T} \ell ( x^*(t), u^*(t) ) - \lambda(x^*(T+1) ) \]
\[ = \left ( \sum_{t = 0}^{T-1} \ell ( x^*(t), u^*(t) ) \right ) + \ell (x^*(T),u^*(T))  - \lambda(x^*(T+1) ) \]
\[ \geq  \left ( \sum_{t = 0}^{T-1} \ell ( x^*(t), u^*(t) ) \right ) - \lambda (x^*(T) ) \geq V^{\psi}_{T} (x), \]
where the first inequality holds by the dissipativity assumption, and the second because $x^*, u^*$ is a feasible solution also over the shorter horizon $[0,T]$.
Hence, $V^{\psi}_{T} (x)$ is monotone non-decreasing with respect to $T$ and the limit (\ref{limitcost}) exists. 

It is important to realise that Proposition \ref{firstprop} only guarantees existence of the limit, not actual boundedness of the cost $V^{\psi}_{\infty} (x)$. In fact,
typically the cost would be $+ \infty$ unless a suitably shifted version of $\ell(x,u)$ is considered.
In particular, there is only a single value of this shift that might result in a finite cost. This can be found, by alternative means, looking for the optimal average cost,
\begin{equation}
\begin{array}{rl}
 V^{\textrm{avg}} (x) =  \inf_{x(\cdot), u ( \cdot) } & J^{\textrm{avg}} (x(\cdot),u(\cdot) ) \\
\textrm{subject to} & \\
x(0) & = x \\
x(t+1) & = f(x(t),u(t)) \qquad t \in \mathbb{N} \\
 (x(t),u(t)) & \in \mathbb{Z} \qquad  t \in \mathbb{N}. 
\end{array}
\end{equation}
Under suitable technical conditions, for instance global controllability assumptions, the optimal cost is independent of $x$, and its value can be found \cite{gaitsgory1,gaitsgory0} by an infinite dimensional linear program, viz. by solving the following optimisation problem:
\begin{equation}
\label{infinitelp}
\begin{array}{l}
V^{\textrm{avg}} =  \sup_{\lambda(\cdot) \in \mathcal{C} (\mathbb{X}) } \; c \\
\qquad \qquad \textrm{subject to} \\
\qquad \qquad \lambda (f(x,u)) \leq \lambda (x) + \ell(x,u) - c \qquad \quad \forall \, (x,u) \in \mathbb{Z} \\
\end{array}
\end{equation}
where:
\[ \mathcal{C} ( \mathbb{X} ) := \{ \lambda: \mathbb{X} \rightarrow \mathbb{R}: \lambda \textrm{ is continuous } \}. \]

We note that this approach has similarities to the effective Hamiltonian approach in continuous-time ergodic optimal control, see \cite{AlBM07}.
Dynamic programming allows to solve optimal control problems through iteration of a suitably defined operator, which computes the optimal cost for increasing values of the control horizon.
To this end, for summed costs without exponential rescaling, the following Bellman operator is normally defined: $T: \mathcal{C} ( \mathbb{X} ) \rightarrow \mathcal{C} (\mathbb{X} )$.
\begin{equation}
\label{bellman}
T \psi (x):= \min_{u \in \mathbb{U} (x)} \ell (x,u) + \psi ( f(x,u) ). 
\end{equation}

The following result characterizes $V_{\infty}^{\psi}(x)$ as a fixed-point of the Bellman operator.
\begin{prop}
\label{limitofvt}
Assume that $\psi= - \lambda$ for some storage function  $\lambda \in \mathcal{C} ( \mathbb{X} )$ and that the following limit exists and is finite:
\begin{equation}
V_{\infty}^{\psi} (x ) = \lim_{T \rightarrow + \infty} V_T^{\psi} (x). 
\end{equation}
Then, $V_{\infty}^{\psi}$ is a lower semi-continuous solution of the Bellman Equation, viz. $T V_{\infty}^{\psi} (x) = V_{\infty}^{\psi} (x)$.
\end{prop}
\emph{Proof.} To see this, recall that $V_T^{\psi}(x)$ is non-decreasing with respect to $T$.
Hence:
\[ \liminf_{x \rightarrow x_0} V_{\infty}^{\psi} (x) = \liminf_{x \rightarrow x_0} \lim_{T \rightarrow + \infty} V_T^{\psi} (x)  \geq \liminf_{x \rightarrow x_0} V_T^{\psi} (x) = V_T^{\psi} (x_0) \qquad    \forall \, T \in \mathbb{N} \]
Since $T$ is arbitrary, we see that:
\[  \liminf_{x \rightarrow x_0} V_{\infty}^{\psi} (x) \geq \lim_{T \rightarrow + \infty} V_T^{\psi} (x_0) = V_{\infty}^{\psi} (x_0). \]
This proves that $V^{\psi}_{\infty}$ is lower semicontinuous. Hence the minimum of
\[ \min_{u \in \mathbb{U}(x)} \ell (x,u) +  V_{\infty}^{\psi} ( f(x,u) ),  \]
is achieved, for some optimal feedback policy $u^*(x)$.
Moreover it fulfills:
\begin{eqnarray*} T V_{\infty}^{\psi} (x) & = &  \ell (x,u^*(x)) +  V_{\infty}^{\psi} ( f(x,u^*(x)) ) \;   = \; \lim_{T \rightarrow + \infty} \ell (x,u^*(x)) + V_T^{\psi} ( f(x,u^*(x))) \\
& \geq & \lim_{T \rightarrow + \infty} V_{T+1}^{\psi} (x) = V_{\infty}^{\psi} (x). \end{eqnarray*}

On the other hand:
\[ V_{\infty}^{\psi} (x) =  \lim_{\tau \rightarrow + \infty} V_{\tau+1}^{\psi} (x) = \lim_{\tau \rightarrow + \infty} T V_{\tau}^{\psi} (x) = \lim_{\tau \rightarrow + \infty} \min_{u \in \mathbb{U}(x) } \ell(x,u) + V_{\tau}^{\psi} ( f(x,u) ). \]
Let $x \in \mathbb{X}$ be fixed and arbitrary. Since $V_\tau$ is continuous in $x$, for each $\tau>0$ and the current fixed value of $x$ there exists a minimizer $u_\tau (x) \in \mathbb{U}(x)$ for this last expression. Since $\mathbb{U}(x)$ is compact, we find a sequence $\tau_n\to\infty$ (possibly $x$-dependent) such that $u_{\tau_n}$ converges to a control value $u_\infty (x) \in \mathbb{U}(x)$. For each $T>0$ this implies 
\begin{eqnarray}
\label{inequalitypreliminary}
V_{\infty}^{\psi} (x) & = & \lim_{n \rightarrow + \infty} \ell(x,u_{\tau_n}(x)) + V_{\tau_n}^{\psi} ( f(x,u_{\tau_n}(x)) )\\
& \ge & \lim_{n \rightarrow + \infty} \ell(x,u_{\tau_n}(x)) + V_{T}^{\psi} ( f(x,u_{\tau_n}(x)) ) \; = \; \ell(x,u_{\infty}(x)) + V_{T}^{\psi} ( f(x,u_{\infty}(x)) ). \nonumber
\end{eqnarray}
Since $V_T^{\psi} (x) \rightarrow V_{\infty}^{\psi} (x)$, for each $\varepsilon>0$ there exists $T_\varepsilon (x)>0$ such that $V_{T_\varepsilon(x)}^{\psi} (x) \geq V^{\psi}_{\infty} (x) - \varepsilon$. 
Hence we see, starting from (\ref{inequalitypreliminary}):
\begin{eqnarray*} V_{\infty}^{\psi} (x) & \ge &  \ell(x,u_{\infty}(x)) + V_{T_\varepsilon (f(x,u_{\infty}(x)))}^{\psi} ( f(x,u_{\infty}(x)) ) \; \ge \; \ell(x,u_{\infty}(x)) + V_{\infty}^{\psi} ( f(x,u_{\infty}(x)) ) - \varepsilon \\
& \ge & \inf_{u \in \mathbb{U}(x) } \ell(x,u) + V_{\infty}^{\psi} ( f(x,u) ) - \varepsilon \; = \; T V_{\infty}^{\psi} (x) - \varepsilon. 
\end{eqnarray*}
Since $x\in\mathbb{X}$ and $\varepsilon>0$ were arbitrary, the assertion $V_{\infty}^{\psi} (x) \geq T V_{\infty}^{\psi} (x)$ follows for all $x \in \mathbb{X}$.

\section{Shifted Bellman Equation and operators}
\label{shiftedsection}
In the literature, different constructive approaches for computing storage functions are described, above all the classical constructions of the available storage and the required supply, which go back to \cite{willems1} and are easily adapted to the discrete-time case (see, e.g., \cite{dissipativity2,MULLERgrune} for the available storage). For this reason, a possible, but ultimately unsatisfactory, way to approach an infinite horizon optimal control problem would be according to the following steps:
\begin{enumerate}
    \item Computing the minimal average cost, $V^{\textrm{avg}}$; 
\item Defining a shifted stage cost, $\tilde{\ell} (x,u) = \ell (x,u) - V^{\textrm{avg}}$, so as to yield $0$ optimal average; 
\item Compute a storage function $\lambda$ for the supply function $\tilde{\ell}(x,u)$;
\item Defining $\psi := - \lambda$ as a terminal penalty for the infinite horizon optimal control problem, with shifted stage costs $\tilde{\ell}$;
\item Use the standard Bellman iteration to asymptotically compute the value function over an infinite horizon or directly looking for a fixed point of the associated Bellman Equation.
\end{enumerate}
This procedure is non ideal for several reasons: first of all, computation of the optimal average cost involves a limiting operation, and therefore typically only approximate values for $V^{\textrm{avg}}$ can ever be achieved.
However, using approximate values in the iteration of the Bellman operator, yields diverging optimal costs over an infinite horizon, either to $\pm \infty$, depending on whether the optimal average cost has been over or underestimated. 
In addition, Step 3 is bound to fail whenever the average optimal cost $V^{\textrm{avg}}$ has been overestimated (in other words a storage function might exist only for $\tilde{\ell}(x,u)=\ell(x,u) - c$ where $c \leq V^{\textrm{avg}}$). 

The goal of this section is to propose operators, the \emph{$\min$-shifted and $\max$-shifted Bellman operator}, whose iteration would converge to the optimal infinite horizon cost, and, at the same time, yield as a by-product the optimal average cost. 

To this end, we need additional notation.
Given $\psi_1: \mathbb{X} \rightarrow \mathbb{R}$ and $\psi_2: \mathbb{X} \rightarrow \mathbb{R}$, continuous, we define the following:
\begin{equation}
c ( \psi_1, \psi_2 ) := \frac{1}{2} \max_{x \in \mathbb{X}} [\psi_1 (x) - \psi_2 (x)]   + \frac{1}{2} \min_{x \in \mathbb{X} } [\psi_1 (x) - \psi_2(x)].
\end{equation}
The following distance notion is also defined:
\[ d( \psi_1, \psi_2 ) :=   \min_{c \in \mathbb{R} } \|  \psi_1 - \psi_2 + c \|_{\infty}. \]
Notice that $d( \psi_1 + c_1 , \psi_2 + c_2 ) = d ( \psi_1, \psi_2 )$ for all $c_1, c_2 \in \mathbb{R}$.
Moreover:
\[ d( \psi_1, \psi_2 ) = \| \psi_1 -  \psi_2 - c ( \psi_1,\psi_2 ) \|_{\infty}. \] 
In fact, an equivalent alternative definition for $d(\psi_1,\psi_2)$ is as follows:
\[ d( \psi_1,\psi_2 ) =
\frac{1}{2} \max_{x \in \mathbb{X}} [ \psi_1 (x) - \psi_2 (x) ]
- \frac 12 \min_{x \in \mathbb{X}} [ \psi_1 (x) - \psi_2 (x) ]. \]

Recall the Bellman operator $T: \mathcal{C} ( \mathbb{X} ) \rightarrow \mathcal{C} (\mathbb{X} )$ previously introduced:
\[ T \psi := \min_{u \in \mathbb{U} (x)} \ell (x,u) + \psi ( f(x,u) ). \]
\begin{defn} Define the $\min$-shifted Bellman operator $\hat{T}: \mathcal{C} ( \mathbb{X}) \rightarrow \mathcal{C} ( \mathbb{X} )$ as:
\begin{equation}
\label{shiftedoperator}
\hat{T} \psi := \min \{ \psi, T \psi + c ( \psi, T \psi ) \}. 
\end{equation}
\end{defn}
Similarly, we may consider the following operator.
\begin{defn} Define the $\max$-shifted Bellman operator 
 $\check{T}: \mathcal{C} ( \mathbb{X}) \rightarrow \mathcal{C} ( \mathbb{X} )$ as:
\begin{equation}
\label{shiftedoperator2}
\check{T} \psi := \max \{ \psi, T \psi + c ( \psi, T \psi ) \}. 
\end{equation}
\end{defn}
It is straightforward to see that:
\[ \psi(x) \geq \hat{T} ( \psi ) (x) \geq \hat{T}^2 ( \psi ) (x) \geq \ldots \geq \hat{T}^k ( \psi ) (x) \geq \ldots \]
for all $k \in \mathbb{N}$. 
Opposite inequalities hold in the case of the $\check{T}$ operator:
\[ \psi(x) \leq \check{T} ( \psi ) (x) \leq \check{T}^2 ( \psi ) (x) \leq \ldots \leq \check{T}^k ( \psi ) (x) \leq \ldots \]
\begin{rem}
By induction, and exploiting the $\min$ commutativity property, the following formula can be proved (see Appendix \ref{formulaproof}):
\begin{equation}
    \label{minformula}
\hat{T}^k \psi (x) = \min_{\tau \in \{0, \ldots, k \}}  \left \{ T^{\tau} \psi(x) + \min_{S \subseteq \{0,\ldots,k-1 \}, |S|= \tau  }  \sum_{s \in S} c ( \hat{T}^s \psi, T \hat{T}^s \psi )  \right \}. 
\end{equation}
Along similar lines the following inequality can be shown by induction for the $\check{T}$ operator:
\begin{equation}
\label{maxformula}
\check{T}^k \psi (x) \geq \max_{\tau \in \{0, \ldots, k \}}  \left \{ T^{\tau} \psi(x) + \max_{ S \subseteq \{0,\ldots,k-1 \}, |S|= \tau }  \sum_{s \in S} c ( \check{T}^s \psi, T \check{T}^s \psi )  \right \}. 
\end{equation}
\end{rem}

The following result holds:
\begin{prop}
\label{fixedpointequalshifted}
\LG{A} function $\bar{\psi} (x) \in \mathcal{C}( \mathbb{X})$ is a fixed point of $\hat{T}$ \LG{ or $\check{T}$} if and only if there exists $c \in \mathbb{R}$ such that $\bar{\psi}$ is a fixed point of the following shifted Bellman Equation:
\begin{equation}
\label{shiftedBE}
T \bar{\psi} = \bar{ \psi } + c.
\end{equation}
 \end{prop}
\emph{Proof.} Assume that $\bar{ \psi }$ fulfills the shifted Bellman Equation (\ref{shiftedBE}).
Then, direct computation shows:
\[ \hat{T} \bar{ \psi} = \min \{ \bar{\psi}, T \bar{\psi} + c ( \bar{ \psi}, T \bar{ \psi} ) \}  =  \min \{ \bar{\psi},  \bar{\psi} + c + c ( \bar{ \psi}, \bar{ \psi} + c ) \}  = \min \{ \bar{ \psi} , \bar{ \psi} \} = \bar{ \psi}, \]
where the equality follows since by definition $c ( \bar{ \psi} , \bar{ \psi } + c ) = -c$.
Conversely, assume $\hat{T} \bar{ \psi} = \bar{ \psi }$:
\[ \bar{ \psi } = \min \{ \bar{ \psi } , T \bar{ \psi} + c ( \bar{\psi}, T \bar{ \psi } ) \}. \]
Hence, the following inequality holds:
\begin{equation}
\label{consequence1}
 \bar{ \psi } (x) \leq T \bar{\psi} (x) + c ( \bar{ \psi} , T \bar{ \psi} ) \; \quad \forall \, x \in \mathbb{X}. 
\end{equation}
We claim that more is true, namely:
\begin{equation}
\label{isshifted}
 \bar{ \psi } (x) - T \bar{ \psi} (x) = c ( \bar{ \psi} , T \bar{ \psi} ) \qquad \forall \, x \in \mathbb{X}.
\end{equation}
Assume by contradiction:
\[ \min_{x \in \mathbb{X} }   \bar{ \psi } (x) - T \bar{ \psi} (x) < c ( \bar{ \psi}, T \bar{ \psi} ), \]
where the $\min$ exists by continuity of $\bar{\psi}$, $T \bar{ \psi}$ and compactness of $\mathbb{X}$.
By inequality (\ref{consequence1}) we also know that:
\[  \max_{x \in \mathbb{X} }   \bar{ \psi } (x) - T \bar{ \psi} (x) \leq c ( \bar{ \psi}, T \bar{ \psi} ). \]
Taking a convex combination of the two previous inequalities yields:
\[ c( \bar{ \psi}, T \bar{\psi} ) = \frac{1}{2}  \min_{x \in \mathbb{X} }   \bar{ \psi } (x) - T \bar{ \psi} (x) + \frac{1}{2}  \max_{x \in \mathbb{X} }   \bar{ \psi } (x) - T \bar{ \psi} (x) < c ( \bar{ \psi}, T \bar{ \psi} ), \]
which is a contradiction.
Hence, (\ref{isshifted}) holds, and $\bar{ \psi}$ is solution of a shifted Bellman Equation.  A similar proof applies to the $\check{T}$ operator.\\
\section{Properties of $T$, $\hat{T}$ and $\check{T}$ operators}
\label{propertiesoperators}
Throughout this section we recall some useful properties of the $T$ operator and additionally provide original derivations for the properties of the $\hat{T}$ and $\check{T}$ operators.
Some of the properties listed below are well known and can be found in \cite{bertsekas3}:
\begin{itemize}
\item Monotonicity: 
\[ [ \psi_1 (x) \leq \psi_2 (x), \; \forall \, x \in \mathbb{X}] \Rightarrow [ T \psi_1 (x) \leq T \psi_2 (x) , \; \forall \, x \in \mathbb{X} ] \]
\item Translation invariance:
\[ T ( \psi + c ) = T \psi + c, \]
for any constant $c \in \mathbb{R}$;
\item Minimum commutativity, for finite index set $K$:
\[   T \left ( \min_{k \in K} \{ \psi_k  \} \right  ) = \min_{k \in K}  \{ T \psi_k \} \]
To see the last one, notice:
\begin{eqnarray*} T \left ( \min_{k \in K} \{ \psi_k \} \right ) & = & \min_{u \in \mathbb{U}(x) } \ell (x,u) + \min_{k \in K}  \{ \psi_k (f(x,u)) \} \; = \; \min_{u \in \mathbb{U}(x) } \min_{k \in K} \{ \ell(x,u) + \psi_k (f(x,u) \}\\
& = & \min_{k \in K} \min_{u \in \mathbb{U}(x)} \ell(x,u) + \psi_k (f(x,u)) \; = \; \min_{k \in K} \{ T \psi_k \}. \end{eqnarray*}
\item Concavity:\\
 For all $\alpha \in [0,1]$ and any $\psi_1, \psi_2$ it holds:
\[ T ( \alpha \psi_1 + (1-\alpha) \psi_2 ) \geq \alpha T \psi_1 + (1- \alpha) T \psi_2. \]
To see this, notice:
\begin{eqnarray*}  T ( \alpha \psi_1 + (1-\alpha) \psi_2 ) & = & \min_{u \in \mathbb{U}(x) }\ell(x,u) + \alpha \psi_1 (f(x,u)) + (1- \alpha) \psi_2 (f(x,u)) \\
& = & \min_{u \in \mathbb{U} (x) } \alpha [ \ell(x,u) + \psi_1 ( f(x,u)) ] + (1- \alpha) [ \ell(x,u) + \psi_2 (f(x,u) ) ] \\
& \geq &  \min_{u \in \mathbb{U} (x) } \alpha [ \ell(x,u) + \psi_1 ( f(x,u)) ] + \min_{u \in \mathbb{U} (x) }(1- \alpha) [ \ell(x,u) + \psi_2 (f(x,u) ) ]\\
& = & \alpha T \psi_1 + (1- \alpha) T \psi_2. \end{eqnarray*}
\item Max-super-commutativity: the following inequality holds:
\[  T \max \{ \psi_1, \psi_2 \} \geq
\max \{ T \psi_1 , T \psi_2 \}, \]
and by induction, for any finite set $K$:
\[ T \left ( \max_{k \in K} \{ \psi_k(x) \} \right ) \geq \max_{k \in K} \{ T \psi_k(x) \}. \]
\item Non-expansiveness: monotonicity and shift-invariance can be exploited to show the following inequality, expressing (incremental) non-expansiveness of the $T$ operator:
\[ d (T \psi_1, T \psi_2 ) \leq d ( \psi_1, \psi_2 ), \qquad \forall \, \psi_1, \psi_2 \in \mathcal{C}( \mathbb{X} ). \]
\end{itemize}
Next we derive some useful properties of the $\hat{T}$ and $\check{T}$ operators.
Notice that for all $c_1,c_2 \in \mathbb{R}$ the following holds:
\[   c ( \psi_1+c_1 , \psi_2+c_2 ) = c( \psi_1, \psi_2) + c_1 - c_2. \]
Hence the following translation invariance can be seen:
\[  \hat{T} ( \psi + c ) = \hat{T} \psi  + c, \]
for all $c \in \mathbb{R}$. In fact,
\begin{eqnarray*}
\hat{T} ( \psi + c ) & = & \min \{ \psi + c, T ( \psi + c) + c ( \psi + c, T ( \psi + c ) ) \} \; = \;
\min \{ \psi + c, T \psi + c + c ( \psi + c , T \psi + c ) \} \\
&  = & \min \{ \psi + c, T \psi + c + c ( \psi  , T \psi  ) \} \; = \; \min   \{ \psi , T \psi  + c ( \psi  , T \psi  ) \} + c \; = \; \hat{T} \psi  + c. 
\end{eqnarray*}
The same property holds for $\check{T}$.
The next proposition states that all solutions of a shifted Bellman Equation share the same shift value.
\begin{prop}
\label{sameaverage}
Let $\psi_1$ and $\psi_2$ be continuous solutions of the shifted Bellman Equation (\ref{shiftedBE}), viz.
$T \psi_1 + c_1 = \psi_1$ and $T \psi_2 + c_2 = \psi_2$ for suitable constants $c_1$ and $c_2$.
Then, $c_1=c_2$.
\end{prop}
\emph{Proof.} See Appendix \ref{sameaverageproof}. \\
 
We show later, by means of an example, that while the shift is uniquely defined for all solutions of the shifted Bellman Equation, it is not true in general that $d ( \psi_1, \psi_2) =0$, i.e. there may be multiple solutions of the shifted Bellman Equation, even after taking into account translation invariance. 
In the remainder of this section, we describe a situation in which the solution of the shifted Bellman Equation is unique, up to the addition of a constant. Again, a dissipativity inequality plays a role, but now a stronger one than \eqref{eq:diss}. For an equilibrium $(x^e,u^e)$ we call the system {\em strictly dissipative}, if there exists a storage function $\lambda: \mathbb{X} \rightarrow \mathbb{R}$, bounded from below, and\footnote{As usual we define $\mathcal K$ as the set of continuous functions $\alpha:[0,\infty)\to[0,\infty)$ that are strictly increasing with $\alpha(0)=0$.} $\alpha\in\mathcal{K}$ such that
\begin{equation}
\lambda ( f(x,u) ) \leq \lambda (x) + \ell(x,u) - \ell(x^e,u^e)-\alpha(\|x-x^e\|) \qquad \forall \, (x,u) \in \mathbb{Z}. \label{eq:sdiss}
\end{equation}
We note that a positive definite stage cost, i.e., an $\ell$ satisfying $\ell(x,u) \ge \alpha(\|x-x^e\|)$ for all $(x,u)\in\mathbb{Z}$ and $\ell(x^e,u^e)=0$, satisfies the inequality \eqref{eq:sdiss} for $\lambda\equiv 0$. For this kind of stage costs, the following proposition holds.

\begin{prop} Suppose the stage cost $\ell$ satisfies $\ell(x,u) \ge \alpha(\|x-x^e\|)$ for all $(x,u)\in\mathbb{Z}$ and some $\alpha\in\mathcal{K}$, and $\ell(x^e,u^e)=0$. 
Then, up to the addition of a constant, there exists at most one continuous solution of the shifted Bellman Equation.
\label{prop:uniquepd}\end{prop}
\emph{Proof.} Let $\psi_1$ and $\psi_2$ be two continuous solutions of the shifted Bellman Equation \eqref{shiftedBE} that are bounded from below. By adding suitable constants, we can assume that $\psi_1(x^e)=\psi_2(x^e)=0$. From \eqref{bellman} we obtain that
\[ \psi_i(x^e) + c =  T\psi_i(x^e) = \min_{u \in \mathbb{U} (x^e)} \ell (x^e,u) + \psi_i ( f(x^e,u) ) \le \ell (x^e,u^e) + \psi_i ( f(x^e,u^e) ) = \psi_i(x^e), \]
implying $c\le 0$.

For each $x\in\mathbb{X}$, let $u_i^*(x)\in\mathbb{U}(x)$ be a control that realizes the minimum in the Bellman operator \eqref{bellman} for $\psi=\psi_i$, $i=1,2$. Such a $u_i^*(x)$ exists because $\ell$, $f$, and $\psi_i$ are continuous and $\mathbb{U}(x)$ is compact.
Then from the shifted Bellman Equation we obtain that 
\[ \psi_i(x) + c = \ell(x,u_i^*(x)) + \psi_i(f(x,u_i^*(x))),\]
implying
\begin{equation} \psi_i(f(x,u_i^*(x))) = \psi_i(x) + c - \ell(x,u_i^*(x)) \le \psi_i(x) - \alpha(\|x-x^e\|).\label{eq:decay}\end{equation}
Now, given $x_i^*(0)\in\mathbb{X}$, by $x_i^*(k)$ we denote the sequence generated by $x_i^*(k+1) = f(x_i^*(k),u_i^*(x_i^*(k))$. Then \eqref{eq:decay} implies 
\[ \psi_i(x_i^*(k)) \le \psi_i(x_i^*(0)) - \sum_{k'=0}^{k-1} \alpha(\|x_i^*(k')-x^e\|).\]
Since $\psi_i$ is bounded from below in $\mathbb{X}$, this sum must converge, implying that $\alpha(\|x_i^*(k)-x^e\|)\to 0$ and thus $x_i^*(k)\to x^e$ as $k\to\infty$. Since $\psi_i(x^e)=0$ and $\psi_i$ is continuous, we also obtain $\psi_j(x_i^*(k))\to 0$ as $k\to\infty$ for $i=1,2$ and $j=1,2$.

Now pick an arbitrary $x\in\mathbb{X}$. We show that for each $\varepsilon>0$ and for both choices $i=1$, $j=2$ and $i=2$, $j=1$ we have 
\begin{equation} \psi_j(x) - \psi_i(x) < \varepsilon
\label{eq:epssmall}    
\end{equation}
holds, which shows $\psi_1(x)=\psi_2(x)$ and thus the assertion.

To this end, consider the sequence $x_i^*(k)$ with $x_i^*(0)=x$. For each $k\ge 0$ we obtain, using that $c$ must be the same in the shifted Bellman Equation for $\psi_j$ and $\psi_i$ due to Proposition \ref{sameaverage},
\begin{eqnarray*} && \psi_j(x_i^*(k)) - \psi_i(x_i^*(k)) = T\psi_j(x_i^*(k)) + c - (T\psi_i(x_i^*(k)) + c)\\[2ex]
&& = \; \underbrace{\min_{u \in \mathbb{U} (x_i^*(k))} \ell (x_i^*(k),u) + \psi_j ( f(x_i^*(k),u) )}_{\le \, \ell (x_i^*(k),u_i^*(x_i^*(k))) + \psi_j ( f(x_i^*(k),u_i^*(x_i^*(k)))} \; - \underbrace{\min_{u \in \mathbb{U} (x_i^*(k))} \ell (x_i^*(k),u) + \psi_i ( f(x_i^*(k)e,u) )}_{= \, \ell (x_i^*(k),u_i^*(x_i^*(k))) + \psi_i ( f(x_i^*(k),u_i^*(x_i^*(k)))}\\[2ex]
&& \le \; \psi_j ( f(x_i^*(k),u_i^*(x_i^*(k)))  - \psi_j ( f(x_i^*(k),u_i^*(x_i^*(k))) \; = \; \psi_j(x_i^*(k+1)) - \psi_i(x_i^*(k+1)).
\end{eqnarray*}
Iterating this inequality we thus obtain
\[ \psi_j(x) - \psi_i(x) \le \psi_j(x_i^*(k)) - \psi_i(x_i^*(k))\]
for all $k\ge 0$. Since we know that $\psi_j(x_i^*(k))\to 0$ and $\psi_i(x_i^*(k))\to 0$ as $k\to\infty$, there is $k\in\mathbb{N}$ such that both $|\psi_j(x_i^*(k))|<\varepsilon/2$ and $|\psi_j(x_i^*(k))|<\varepsilon/2$ hold, implying $\psi_j(x_i^*(k)) - \psi_i(x_i^*(k))<\varepsilon$ and thus \eqref{eq:epssmall}.
\endproof

Now for a strictly dissipative system satisfying \eqref{eq:sdiss} we consider the ``rotated'' stage cost 
\begin{equation}\label{eq:rotell}
    \tilde \ell(x,u) = \ell(x,u) - \ell(x^e,u^e) + \lambda(x) - \lambda(f(x,u))
\end{equation}
and observe that it satisfies the conditions on $\ell$ from Proposition \ref{prop:uniquepd}. The corresponding Bellman operator defined by
\[\wT \psi(x) :=  \min_{u \in \mathbb{U} (x)} \tilde \ell (x,u) + \psi ( f(x,u) )\]
satisfies the following property. 
\begin{lem}\label{lem:tildeT1} For any continuous function $\lambda:\mathbb{X}\to \mathbb{R}$ the identity
\[ \wT\psi = T(\psi-\lambda)+\lambda - \ell(x^e,u^e)\]
holds. Particularly, if $\psi$ is a solution of the shifted Bellman Equation for $T$ and some $c$, then $\tilde\psi = \psi+\lambda$ is a solution of the shifted Bellman Equation for $\wT$ and $\tilde c = c - \ell(x^e,u^e)$.
\end{lem}
\emph{Proof.} For all $x\in\mathbb{X}$ we have that
\begin{eqnarray*} 
    \wT\psi(x) & = & \min_{u \in \mathbb{U} (x)} \{ \tilde\ell (x,u) + \psi ( f(x,u) ) \}\\
    & = & \min_{u \in \mathbb{U} (x)} \{\ell(x,u) - \ell(x^e,u^e) + \lambda(x) - \lambda(f(x,u)) + \psi ( f(x,u) ) \} \\
    & = & \min_{u \in \mathbb{U} (x)} \{ \ell(x,u)  + \psi ( f(x,u) ) - \lambda(f(x,u))\} + \lambda(x) - \ell(x^e,u^e)\\
    & = & T(\psi-\lambda)(x) +\lambda(x) - \ell(x^e,u^e).
\end{eqnarray*}
This proves the first statement. Now, if $\psi$ is a solution of the shifted Bellman Equation for $T$, then
\[\wT\tilde\psi = T(\tilde \psi - \lambda) + \lambda - \ell(x^e,u^e) = T\psi + \lambda - \ell(x^e,u^e) = \psi + c + \lambda - \ell(x^e,u^e) = \tilde \psi + c - \ell(x^e,u^e),\]
i.e.\ $\tilde \psi$ is a solution of the shifted Bellman Equation for $\wT$.
\endproof

\begin{thm}
   Consider an optimal control problem for which strict dissipativity \eqref{eq:sdiss} holds with a continuous storage function $\lambda$. 
   Then, up to the addition of a constant, there exists at most one continuous solution of the shifted Bellman Equation. 
\label{thm:unique}\end{thm}
{\em Proof.}
Let $\psi_1$ and $\psi_2$ be two solutions of the shifted Bellman Equation satisfying the assumption. Then $\tilde\psi_i = \psi_i+\lambda$, $i=1,2$ satisfy the assumption of Proposition \ref{prop:uniquepd} since $\lambda$ is continuous and bounded from below. Hence, applying Proposition \ref{prop:uniquepd} to $\wT$ yields that $\psi_1 + c - \ell(x^e,u^e)$ and $\psi_2 + c - \ell(x^e,u^e)$ coincide up to the addition of a constant, implying the same for $\psi_1$ and $\psi_2$.
\endproof 

We note that non-strict dissipativity is not enough to obtain this uniqueness result up to additions of constants, as the example in Subsection \ref{nondissipativeex} shows.

\section{Convergence analysis under equicontinuity}
\label{equicont}
In order to prove convergence of the $\hat{T}$ and $\check{T}$ iterations to a fixed point of the shifted Bellman Equation we restrict the dynamics to fulfill suitable equicontinuity assumptions. Moreover, we provide sufficient conditions, in the form of controllability assumptions, which lead to the needed equicontinuity properties both for the iteration $T^k \psi$ and $\hat{T}^k \psi$.

In order to have convergence guarantees for a sequence of functions, the following notion of equicontinuity is adopted.
\begin{defn} A sequence of functions $\{ \psi_k \}_{k=0}^{+ \infty}$, $\psi_k : \mathbb{X} \rightarrow \mathbb{R}$ is said to be
equicontinuous, if there exists a function $\gamma \in \mathcal{K}_{\infty}$ such that:
\[ \forall \, k \in \mathbb{N}, \; \forall \, x_1, x_2 \in \mathbb{X}: \qquad | \psi_k (x_1) - \psi_k (x_2 ) |\leq \gamma (| x_1 - x_2|). \]
\end{defn}
To carry out our analysis, we will need the following assumption.
\begin{assum}
\label{tkequi}
The sequence $\{ T^k \psi \}_{k = 0}^{+\infty}$ is equicontinuous.
\end{assum}

The following lemma shows that this assumption immediately carries over to $\hat{T}^k\psi$.
\begin{lem}
\label{thatequi}
The sequence $\{ \hat{T}^k \psi \}_{k=0}^{+ \infty}$ is equicontinuous provided $\{ T^k \psi\}_{k=0}^{+ \infty}$ is such.
\end{lem}
\emph{Proof.}
The lemma is a simple consequence of formula (\ref{minformula}). In particular, equicontinuity holds with the same function $\gamma$, i.e. $| \hat{T}^k \psi (x_1) - \hat{T}^k \psi (x_2)| \leq \gamma(|x_1-x_2|)$ provided
$|T^k \psi(x_1) - T^k \psi (x_2)| \leq \gamma (|x_1-x_2|)$. \\

Our main convergence results under equicontinuity are now stated in the following two theorems.

\begin{thm}
\label{mainconvergence}
Let $\psi \in \mathcal{C} ( \mathbb{X} )$ be such that $T^k \psi (x)$ fulfill Assumption \ref{tkequi}. Then, if a continuous fixed point of the shifted Bellman Equation exists, the sequence $\hat{T}^k \psi (x)$ converges uniformly to one such fixed point.
\end{thm}

\emph{Proof.} Consider the sequence $[\hat{T}^k \psi]_n$. By Lemma \ref{boundednessT} this sequence is bounded since:
\begin{eqnarray*} 0 & \leq & [\hat{T}^k \psi(x)]_n \; \leq \; \max_{x \in \mathbb{X}} \hat{T}^k \psi(x) - \min_{x \in \mathbb{X}} \hat{T}^k \psi (x) \\
&  \leq &  \left [ \max_x \bar{ \psi} (x) - \min_x \bar{ \psi} (x) \right ] + \left [    \max_x [ \psi(x) - \bar{ \psi } (x) ]- \min_x [ \psi(x) - \bar{ \psi } (x) ]        \right ].\end{eqnarray*}%
Moreover, by Lemma \ref{thatequi} it is equicontinuous. Hence, by the Arzela-Ascoli Theorem, it admits a non empty set of accumulation points (with respect to the uniform topology),
\[ \omega(\psi) := \{ \bar{\psi} \in \mathcal{C}(\mathbb{X}): \exists \{ k_n \}_{n=1}^{+\infty}, k_n \rightarrow + \infty:
\bar{\psi} = \lim_{n \rightarrow + \infty} \hat{T}^{k_n} \psi \}. \]
Moreover, each accumulation point in $\omega(\psi)$ is continuous and fulfills the same continuity inequality, 
\begin{equation}
\label{inequalityinomega} 
 |\psi(x_1)- \psi(x_2)| \leq \gamma (|x_1-x_2|)  
\end{equation}
By Lemma \ref{newlyap}, the function $W( [\psi]_n)= W (\psi) := d(\psi,T \psi)$ is non-increasing along the iteration of $\hat{T}$, viz.
$W(\hat{T}^k \psi)$ is a non-increasing sequence, bounded from below by $0$.
In addition $W$ is continuous in the topology of uniform convergence.
Hence, the limit $\lim_{k \rightarrow + \infty} W (\hat{T}^k \psi)$ exists, and we denote it by $\bar{W}$.
Because of continuity of $W$ and uniform convergence to the limit points we also have
$W ( \bar{\psi} ) = \bar{W}$ for all $\bar{\psi} \in \omega( \psi)$.
Notice that $\omega(\psi)$ is invariant with respect to $\hat{T}$. Hence, for any 
$\bar{\psi} \in \omega(\psi)$ and any $k \in \mathbb{N}$ we have
$W ( \hat{T}^k \bar{\psi} ) = \bar{W}$. 
By combined inequalities (\ref{goingup}) and (\ref{goingdown}) we see that $W ( \hat{T}^k \bar{\psi} )$ can be constant only provided 
 $\min_{x \in \mathbb{X}} \hat{T}^k \bar{\psi} (x) - T \hat{T}^k \bar{\psi}(x)$ and  $\max_{x \in \mathbb{X}} \hat{T}^k \bar{\psi} (x) - T \hat{T}^k \bar{\psi}(x)$ are constant with respect to $k$.
 By Corollary \ref{convergeupper}, the sequence $\hat{T}^k \bar{\psi}$ is bounded and converges monotonically to an upper semi-continuous limit.
 Notice that, by invariance of $\omega(\psi)$ and the fact that all elements of $\omega(\psi)$ fulfill inequality (\ref{inequalityinomega}),  equicontinuity of $\hat{T}^k \bar{\psi}$ follows. Hence
 the limit $\psi_{\infty} (x) := \lim_{k \rightarrow + \infty} \hat{T}^k \bar{\psi} (x)$ not only exists (as previously established), but is also continuous and, by Dini's Theorem, convergence is uniform in $\mathbb{X}$. By continuity of the $\hat{T}$ operator with respect to uniform convergence,
  $\psi_{\infty}(x)$ is a fixed point of the shifted Bellman Equation (cf.\ Lemma \ref{fixcontinuous}) and
 $0=d( \psi_{\infty}, T \psi_{\infty})=d( \bar{\psi}, T \bar{\psi} )$. This shows that any element of $\omega(\psi)$ is an equilibrium of the shifted Bellman Equation.
 We only need to show that $\omega(\psi)$ is a singleton. This follows because of Lemma \ref{Lyapunovoperator2}. Indeed, the distance to any element $\bar{\psi}$ of $\omega(\psi)$ is non increasing along the iteration $\hat{T}^k \psi$. Since such distance is converging to $0$ along some subsequence $\hat{T}^{k_n} \psi$, then it is converging to $0$ along the sequence $\hat{T}^k \psi$ itself.  \endproof

Due to the lack of an analogue to formula (\ref{minformula}) for the $\check{T}$ operator, there is no simple way of proving a version of Lemma \ref{thatequi} for $\check{T}^k \psi$. As a consequence, the analogue of Theorem \ref{mainconvergence} for $\check{T}$ is stated by directly assuming equicontinuity of $\check{T}^k \psi$.
 
\begin{thm}
\label{mainconvergence2}
Let $\psi \in \mathcal{C} ( \mathbb{X} )$ be such that $\check{T}^k \psi (x)$ fulfills Assumption \ref{tkequi}. Then, if a continuous fixed point of the shifted Bellman Equation exists, the sequence $\check{T}^k \psi (x)$ converges uniformly to one such fixed point.
\end{thm} 
 
\emph{Proof.} Consider the sequence $[\check{T}^k \psi]_n$. This sequence is bounded since:
\begin{eqnarray*}
    0 & \leq & [\check{T}^k \psi(x)]_n \; \leq \; \max_{x \in \mathbb{X}} \check{T}^k \psi(x) - \min_{x \in \mathbb{X}} \check{T}^k \psi (x) \\
    & \leq & \left [ \max_x \bar{ \psi} (x) - \min_x \bar{ \psi} (x) \right ] + \left [    \max_x [ \check{T}^k \psi(x) - \bar{ \psi } (x) ]- \min_x [ \check{T}^k \psi(x) - \bar{ \psi } (x) ]        \right ]\\
& \leq & \left [ \max_x \bar{ \psi} (x) - \min_x \bar{ \psi} (x) \right ] + \left [    \max_x [ \psi(x) - \bar{ \psi } (x) ]- \min_x [ \psi(x) - \bar{ \psi } (x) ]        \right ],
\end{eqnarray*}
where the last inequality follows by Lemma \ref{Lyapunovoperator}.
Moreover, by assumption, it is equicontinuous. Hence, by Arzela-Ascoli Theorem, it admits a non empty set of limit points (with respect to the uniform topology),
\[ \omega(\psi) := \{ \bar{\psi} \in \mathcal{C}(\mathbb{X}): \exists \{ k_n \}_{n=1}^{+\infty}, k_n \rightarrow + \infty:
\bar{\psi} = \lim_{n \rightarrow + \infty} \check{T}^{k_n} \psi (x) \}. \]
Note that each limit point in $\omega(\psi)$ is continuous and fulfills the same continuity inequality, 
\begin{equation}
\label{inequalityinomega2} 
 |\psi(x_1)- \psi(x_2)| \leq \gamma (|x_1-x_2|)  
\end{equation}
By Lemma \ref{newlyap2}, the function $W( [\psi]_n)= W (\psi) := d(\psi,T \psi)$ is non-increasing along the iteration of $\check{T}$, viz.
$W(\check{T}^k \psi)$ is a non-increasing sequence, bounded from below by $0$.
In addition $W$ is continuous in the topology of uniform convergence.
Hence, the limit $\lim_{k \rightarrow + \infty} W (\check{T}^k \psi)$ exists, and we denote it by $\bar{W}$.
Because of continuity of $W$ and uniform convergence to the limit points we also have
$W ( \bar{\psi} ) = \bar{W}$ for all $\bar{\psi} \in \omega( \psi)$.
Notice that $\omega(\psi)$ is invariant with respect to $\check{T}$. Hence, for any 
$\bar{\psi} \in \omega(\psi)$ and any $k \in \mathbb{N}$ we have
$W ( \check{T}^k \bar{\psi} ) = \bar{W}$. 
By combined inequalities (\ref{goingdown2}) and (\ref{goingup2}) we see that $W ( \check{T}^k \bar{\psi} )$ can be constant only provided 
 $\min_{x \in \mathbb{X}} \check{T}^k \bar{\psi} (x) - T \check{T}^k \bar{\psi}(x)$ and  $\max_{x \in \mathbb{X}} \check{T}^k \bar{\psi} (x) - T \check{T}^k \bar{\psi}(x)$ are constant with respect to $k$.
 By Corollary \ref{convergelower}, the sequence $\check{T}^k \bar{\psi}$ is bounded and converges monotonically to a lower semi-continuous limit.
 Notice that, by invariance of $\omega(\psi)$ and the fact that all elements of $\omega(\psi)$ fulfill inequality (\ref{inequalityinomega2}) follows equicontinuity of $\check{T}^k \bar{\psi}$, hence
 the limit $\psi_{\infty} (x) := \lim_{k \rightarrow + \infty} \check{T}^k \bar{\psi} (x)$ not only exists (as previously established), but is also continuous and, by Dini's Theorem, convergence is uniform in $\mathbb{X}$. By continuity of the $\check{T}$ operator with respect to uniform convergence,
  $\psi_{\infty}(x)$ is a fixed point of the shifted Bellman Equation and
 $0=d( \psi_{\infty}, T \psi_{\infty})=d( \bar{\psi}, T \bar{\psi} )$. This shows that any element of $\omega(\psi)$ is an equilibrium of the shifted Bellman Equation.
 We only need to show that $\omega(\psi)$ is a singleton. This follows because of Lemma \ref{Lyapunovoperator}. Indeed, the distance to any element $\bar{\psi}$ of $\omega(\psi)$ is non increasing along the iteration $\check{T}^k \psi$. Since such distance is converging to $0$ along some subsequence $\check{T}^{k_n} \psi$, then it is converging to $0$ along the sequence $\check{T}^k \psi$ itself.\endproof

In the remainder of this section we derive a sufficient condition for Assumption \ref{tkequi}, which is based on a controllability condition.

\begin{defn}
\label{incrementalcontrollability}
Given a system as in (\ref{systemmap}) and the associated state and input constraint sets $\mathbb{X}$ and $\mathbb{U}(x)$, we say that the system fulfills Uniform Incremental Continuous Controllability, if there exists $N \in \mathbb{N}$, and a class $\mathcal{K}_{\infty}$ function $\delta$, such that, for all $x_1, x_2 \in \mathbb{X}$, and for all $\textbf{u}_1 \in \mathbb{U}_N (x_1)$, there exists $\textbf{u}_2 \in \mathbb{U}_N (x_2)$ such that $\phi(N,x_1,\textbf{u}_1) = \phi(N,x_2,\textbf{u}_2)$, and in addition:
$\| \textbf{u}_1 - \textbf{u}_2 \| \leq \delta (|x_1-x_2|)$.
\end{defn}
A milder controllability assumption can be formulated by considering continuity with respect to the cost alone, rather than the control input.
To this end, let $J_N(x, \textbf{u})$, for $x \in \mathbb{X}$ and $\textbf{u} \in \mathbb{U}_N (x)$ denote the following:
\[ J_N (x, \textbf{u} ) = \sum_{t=0}^{N-1} \ell (\phi(t,x,\textbf{u}),u(t) ). \]
\begin{defn}
\label{costcontrollability}
Given a system as in (\ref{systemmap}) and the associated state and input constraint sets $\mathbb{X}$ and $\mathbb{U}(x)$, we say that the system fulfills Uniform Incremental Controllability Continuous in Cost, if there exists $N \in \mathbb{N}$, and a class $\mathcal{K}_{\infty}$ function $\delta$, such that, for all $x_1, x_2 \in \mathbb{X}$, and for all $\textbf{u}_1 \in \mathbb{U}_N (x_1)$, there exists $\textbf{u}_2 \in \mathbb{U}_N (x_2)$ such that $\phi(N,x_1,\textbf{u}_1) = \phi(N,x_2,\textbf{u}_2)$, and in addition:
$| J_N (x_1,\textbf{u}_1) - J_N(x_2,\textbf{u}_2) | \leq \delta (|x_1-x_2|)$.
\end{defn}
\begin{rem}
Notice that Uniform Incremental Continuous Controllability implies Uniform Incremental Controllability Continuous in Cost. This is because the considered stage-cost function and the dynamics are both continuous, moreover cost is considered only over a finite interval of length $N$. The converse implication is not true in general. 
\end{rem}

The following proposition now shows that Uniform Incremental Controllability Continuous in Cost implies the equicontinuity in Assumption \ref{tkequi} required in Theorem \ref{mainconvergence}.

\begin{prop} Assume that system (\ref{systemmap}) fulfills the controllability assumption in Definition \ref{costcontrollability}.
Then, for any continuous function $\psi: \mathbb{X} \rightarrow \mathbb{R}$, the sequence $\{ T^k \psi \}_{k=0}^{+ \infty}$ is equicontinuous, i.e., Assumption \ref{tkequi} is fulfilled.
\end{prop}
\emph{Proof.}
Consider any $k\in \mathbb{N}$, and arbitrary $x_1,x_2 \in \mathbb{X}$. Let
$\textbf{u}^*_1 \in \mathbb{U}_{k+N}(x_1)$ be any optimal control sequence corresponding to the optimal control problem
with terminal penalty function $\psi$ and horizon $k+N$, with initial condition $x_1$.
Then, from the optimality principle:
\begin{equation}
\label{cost1}
      T^{k+N} \psi (x_1) = J_N (x_1, \textbf{u}^*_1) + T^k \psi (\phi(N,x_1,\textbf{u}^*_1) ). 
      \end{equation}
Let now, $\textbf{u}_2$ be as in Definition \ref{costcontrollability}. Clearly, applying $\textbf{u}_2$ is, in general, suboptimal from initial condition $x_2$. Hence, the inequality below holds:
\begin{equation}
    \label{cost2}
    T^{k+N} \psi (x_2) \leq J_N (x_2, \textbf{u}_2 ) + T^k \psi (\phi(N,x_2,\textbf{u}_2) ). 
\end{equation}    
Combining equations (\ref{cost1}) and (\ref{cost2}) yields:
\begin{eqnarray*} T^{k+N} \psi (x_2 ) - T^{k+N} \psi (x_1)  &\leq& J_N (x_2, \textbf{u}_2 )  + T^k \psi (\phi(N,x_2,\textbf{u}_2) ) - 
J_N (x_1, \textbf{u}^*_1) - T^k \psi (\phi(N,x_1,\textbf{u}^*_1) ) \\
&=&  J_N (x_2, \textbf{u}_2 ) - J_N (x_1, \textbf{u}^*_1) \; \leq \; \delta (|x_1-x_2|), \end{eqnarray*} 
where the first equality follows because $\phi(N,x_1,\textbf{u}^*_1) = \phi(N,x_2,\textbf{u}_2)$, and the last inequality
from Definition \ref{costcontrollability}. Symmetric inequalities can be obtained swapping $x_1$ and $x_2$, yielding 
$|T^{k+N} \psi(x_1) - T^{k+N} \psi (x_2)| \leq \delta(|x_1-x_2|)$.
This shows that equicontinuity holds on the tail of the sequence $T^k \psi$. However, $\{ T^k \psi \}_{k=0}^{N-1}$ is a finite family of continuous functions defined over a compact set (thus also fulfilling an equicontinuity property), and therefore equicontinuity of the whole sequence follows.\endproof
 
Unfortunately, due to the lack of a counterpart of Lemma \ref{thatequi}, we currently do not have a controllability condition for ensuring the equicontinuity needed in Theorem \ref{mainconvergence2} for the $\check{T}$ operator.

\section{Convergence analysis without continuity}
\label{additionlars}
\def\N{\mathbb{N}}
\def\R{\mathbb{R}}
\def\U{\mathbb{U}}
\def\X{\mathbb{X}}
\def\hwT{\hspace*{1mm}\hat{\hspace*{-1mm}\widetilde{T}}}
\def\hwTk{\hspace*{-1mm}\hat{\hspace*{1mm}\widetilde T^k}}

In this section we provide a convergence result for the iteration using the $\hat T$ operator without assuming any continuity. This is possible if we assume a dissipativity condition and start the iteration from the negative storage function. The result can thus be seen as an extension of Proposition \ref{limitofvt} to the shifted Bellman Equation with nontrivial shift $c\ne 0$.

We first state a little auxiliary lemma, in which for any function $\psi:\X\to\R$ we define
\[ \psi^n(x) := \psi(x) - \min_{x\in\mathbb{X}} \psi(x).\]
We note that $\psi^n\ge 0$ and $\min_{x\in\X}\psi^n(x)=0$ as well as $(\psi+c)^n = \psi^n$ for all $c\in\R$.

\begin{lem}\label{lem:shift} For any $c\in\R$ it holds that 
\[(T(\psi  + c))^n = (T\psi)^n \quad \mbox{ and } \quad (\hat T(\psi  + c))^n = (\hat T\psi)^n.\]
\end{lem}
\emph{Proof.} We have that
\begin{eqnarray*} \hat T (\psi+c) & = & \min\{ \psi+c, \underbrace{T(\psi+c)}_{=T\psi + c} + \underbrace{c(\psi+c, T(\psi+c)}_{=c(\psi,T\psi)}\} \\
& = &  \min\{\psi, T\psi + c(\psi, T\psi)\}+c \; = \; \hat T\psi + c.
\end{eqnarray*}
This implies the assertion since $(\hat T \psi+c)^n = (\hat T\psi)^n$ for all $c\in\R$. A similar computation works for $T$ in place of $\hat T$.
\endproof

We now first consider the case where $\ell\ge 0$. To this end, we make the following assumption.

\begin{assum}\label{ass:N}
There exists a nonempty set $N\subset\X$ such that for any $\psi:\X\to\R$ with $\psi\ge 0$ and $\psi|_N\equiv 0$ we have that $T\psi|_N\equiv 0$.
\end{assum}

We note that this assumption is satisfied for instance if $\ell\ge 0$ and there is an equilibrium $(x^e,u^e)$ (i.e., $f(x^e,u^e)=x^e$) with $\ell(x^e,u^e)=0$. Then one can choose $N=\{x^e\}$.

\begin{lem} \label{lem:l0} Assume $\ell\ge 0$ and let Assumption \ref{ass:N} hold. Then for $\psi^0\equiv 0$ the sequence of functions $\psi^k := (\hat T^k \psi^0)^n$, $k\in\N$, satisfies the following properties for all $k\in\N$:
\[ \begin{array}{llll}
\mbox{\rm (a) } &  T^k \psi^0 \ge \psi^k, \qquad &  \mbox{\rm (b) } & T\psi^k \ge \psi^k, \\
\mbox{\rm (c) } & \psi^k|_N=0, &  \mbox{\rm (d) } & \psi^{k+1} \ge \psi^k.
\end{array}
\]
\end{lem}
\emph{Proof.} By applying Lemma \ref{lem:shift} inductively we see that $\psi^{k+1} = (\hat T \psi^k)^n$. Moreover, we observe for all $\psi:\X\to\R$ the equality
\begin{eqnarray*} (\hat T \psi)^n & = & \min\{\psi,T\psi+c(\psi,T\psi)\}^n \; = \; (\min\{\psi-c(\psi,T\psi),T\psi\}+c(\psi,T\psi) )^n \\
& = & \min\{\psi-c(\psi,T\psi),T\psi\}^n. \end{eqnarray*}
Now we prove (a)--(d) by induction over $k$. 

For $k=0$, (a) and (c) hold trivially, while (b) and (d) hold because $\psi^0\equiv 0$ and $T\psi^0\ge 0$ (since $\ell\ge 0$) and $\psi^1\ge 0$ (by definition of the $(\cdot)^n$ operator).

For $k\to k+1$, assume that (a), (b), and (c) hold for $\psi^k$. We now prove these three properties for $\psi^{k+1}$ and start with (c). By the above computation it holds that
\[ \psi^{k+1} =  ( \hat T \psi^k )^n= \min\{\psi^k-c(\psi^k,T\psi^k),T\psi^k\}^n.\]
By induction assumption (b) we have that $T\psi^k\ge \psi^k$ implying that $c(\psi^k,T\psi^k)\le 0$ and thus $\psi^k-c(\psi^k,T\psi^k)\ge 0$. Since $\ell\ge0$ and $\psi^k\ge 0$ we moreover have $T\psi^k\ge 0$. By induction assumption (c) we know that $\psi^k|_N\equiv 0$. Thus, Assumption \ref{ass:N} yields $T\psi^k|_N\equiv 0$. Together this implies that $\min\{\psi^k-c(\psi^k,T\psi^k),T\psi^k\}\ge 0$ and is equal to $0$ on $N$. This implies that
\begin{equation}\label{eq:non}
    \psi^{k+1} = \min\{\psi^k-c(\psi^k,T\psi^k),T\psi^k\}^n = \min\{\psi^k-c(\psi^k,T\psi^k),T\psi^k\}
\end{equation} 
and thus $\psi^{k+1}|_N\equiv 0$, i.e., (c) for $k+1$. 

Next we prove (b) for $k+1$. Using \eqref{eq:non} as well as the min commutativity and the translation invariance of $T$ we obtain
\begin{eqnarray*}
T\psi^{k+1} & = & T \min\{\psi^k-c(\psi^k,T\psi^k),T\psi^k\} \\
& = & \min\{T\psi^k-c(\psi^k,T\psi^k),TT\psi^k\}.
\end{eqnarray*}
Now using the induction assumption for (b) and the monotonicity of $T$ we obtain $T\psi^k\ge \psi^k$ and $TT\psi^k \ge T\psi^k$, implying, using \eqref{eq:non} once more
\[ \min\{T\psi^k-c(\psi^k,T\psi^k),TT\psi^k\} \ge \min\{\psi^k-c(\psi^k,T\psi^k),T\psi^k\} = \psi^{k+1}. \]
This shows (b) for $k+1$. From the induction assumption (a) and (b) and monotonicity of $T$ we obtain
\[ T^{k+1}\psi^0 = T T^k\psi^0 \ge T\psi^k \ge \psi^k, \]
which shows (a) for $k+1$. 

Finally, for showing (d), we use that the induction assumption for (b) yields $c(\psi^k,T\psi^k)\le 0$ and $T\psi^k\ge\psi^k$. Together with \eqref{eq:non} we obtain 
\[ \psi^{k+1} = \min\{\psi^k-c(\psi^k,T\psi^k),T\psi^k\} \ge \min\{\psi^k,\psi^k\} = \psi^k.\]
\endproof

\begin{prop}\label{prop:hatconv}
Assume $\ell\ge 0$, let Assumption \ref{ass:N} hold and assume that $V_\infty^{\psi^0}$ is finite for $\psi^0\equiv 0$. Then the sequence of functions $\psi^k=(\hat T^k \psi^0)^n$, $k\in\N$, converges to $V_\infty^{\psi^0}$, i.e., in particular to a solution of the Bellman Equation.
\end{prop}
\emph{Proof.} From Lemma \ref{lem:l0} it follows that $\psi^k$ is increasing and bounded from above by $V_\infty^{\psi^0}$. Hence, it converges to some limit function $\psi^\infty\le V_\infty^{\psi^0}$. 
Now from $T\psi^k\ge \psi^k$ we obtain that 
\[ c(\psi^k,T\psi^k) \leq - \frac{1}{2} \max_{\tilde x\in\X} [T\psi^k(\tilde x)-\psi^k(\tilde x)],\]
implying that 
\[ \psi^k(x)-c(\psi^k,T\psi^k)  \ge  \psi^k(x) + \frac{1}{2} \max_{\tilde x\in\X} [T\psi^k(\tilde x)-\psi^k(\tilde x) ] \ge \frac12 (\psi^k(x) + T\psi^k(x)).
\]
Since $T\psi^k\ge \psi^k$ we moreover obtain that $T\psi^k\ge \frac12 (\psi^k(x) + T\psi^k(x))$. Inserting these inequalities into \eqref{eq:non} then yields
\[ \psi^{k+1} = \min\{\psi^k-c(\psi^k,T\psi^k),T\psi^k\} \ge \frac12 (\psi^k + T\psi^k)\]
and using this inequality and $T(\psi_1/2 + \psi_2/2) \ge (T\psi_1)/2 + (T\psi_2)/2$ yields
\begin{eqnarray*}
    \psi^1 & \ge & \frac12 \psi^0 + \frac12 T\psi^0\\
    \psi^2 & \ge & \frac12\psi^1 + \frac12 T\psi^1 \; \ge \; \frac12 \left(\frac12 \psi^0 + \frac12 T\psi^0\right) + \frac12 \left(\frac12 T\psi^0 + \frac12 T^2\psi^0\right) \; = \; \frac14(\psi^0 + 2T\psi^0 + T^2\psi^0)\\
    \psi^3 & \ge & \frac12\psi^2 + \frac12 T\psi^2 \; \ge \; \frac12 \left(\frac12 \psi^1 + \frac12 T\psi^1\right) + \frac12 \left(\frac12 T\psi^1 + \frac12 T^2\psi^1\right) \; \ge \; \frac18(\psi^0+3T\psi^0+3T^2\psi^0+T^3\psi^0)\\
    & \vdots &
\end{eqnarray*}
which by induction yields the general formula
\[ \psi^k \ge \frac1{2^k}\sum_{l=0}^{k} \binom{k}{l}T^l\psi^0.\]
Since $\sum_{l=0}^{k} \binom{k}{l}=2^k$ grows exponentially in $k$ while for each fixed $p\in\N$ 
the sum $\sum_{l=0}^{p-1} \binom{k}{l}$ grows only polynomially in $k$, we have that
\[ \frac{\sum_{l=p}^k \binom{k}{l}}{2^k} = 1 - \underbrace{\frac{\sum_{l=0}^{p-1} \binom{k}{l}}{2^k}}_{\to 0} \to 1\]
as $k\to\infty$. Combining this with $T^q\psi^0\ge T^p\psi^0\ge 0$ for $q\ge p\ge 0$, we obtain that for each $C\in(0,1)$ and $p\in\N$ there is 
$k_{C,p}\in\N$ with 
\[  \psi^k \ge \frac1{2^k}\sum_{l=0}^{k} \binom{k}{l}T^l\psi^0 \ge \frac1{2^k}\sum_{l=p}^{k} \binom{k}{l}T^l\psi^0 \ge \frac1{2^k}\sum_{l=p}^{k} \binom{k}{l}T^p\psi^0 \ge C T^p\psi^0\]
for all $k\ge k_{C,p}$.
This implies that 
\[ \psi^\infty = \lim_{k\to\infty}\psi^k \ge C\lim_{p\to\infty}T^p \psi^0 = CV_\infty^{\psi^0}\]
for any $C\in(0,1)$. Since $C$ can be chosen arbitrarily close to $1$, this implies $\psi^\infty \ge V_\infty^{\psi^0}$,
which finishes the proof. \endproof

Now we extend our results to dissipative stage costs. The dissipativity inequality here is similar to \eqref{eq:diss}, where we explicitly include a shift of the cost function by $c$ in the inequality.

\begin{assum}\label{ass:diss}
There exists a continuous storage function $\lambda:\X\to\R$ and a value $c\in\R$ such that
\begin{equation}
\lambda ( f(x,u) ) \leq \lambda (x) + \ell(x,u) - c \qquad \forall \, (x,u) \in \mathbb{Z}
\end{equation}
\end{assum}

For such a function $\lambda$, similar to \eqref{eq:rotell} we define the rotated cost
\begin{equation}\label{eq:tell}
    \tilde \ell(x,u) = \ell(x,u) - c + \lambda(x) - \lambda(f(x,u))
\end{equation}
and the corresponding operators $\wT$ and $\hwT$. The next lemma extends Lemma \ref{lem:tildeT1}. 
\begin{lem}\label{lem:tildeT} For any continuous function $\lambda:\mathbb{X}\to \R$ and for all $k\in\N$ the identities
\[ \wT^k\psi = T^k(\psi-\lambda)+\lambda - kc \quad \mbox{ and } \quad \hwTk\psi = \hat T^k(\psi-\lambda)+\lambda \]
hold.
\end{lem}
\emph{Proof.} The first identity follows with an analogous proof as for Lemma \ref{lem:tildeT1} followed by induction over $k$. For the second identity we compute
\begin{eqnarray*} 
    \hwT\psi & = & \min\{ \psi,\wT\psi + c(\psi,\wT\psi)\} \\
    & = & \min\{ \psi, T(\psi-\lambda) + \lambda - c + \underbrace{c(\psi,T(\psi-\lambda) + \lambda - c)}_{=c(\psi-\lambda,T(\psi-\lambda)) + c}\}\\
    & = & \min\{ \psi - \lambda, T(\psi-\lambda) + c(\psi-\lambda,T(\psi-\lambda))\} + \lambda\\
    & = & \hat T(\psi-\lambda) + \lambda.
\end{eqnarray*}
From this, the statement for $\hwTk\psi$ follows by induction over $k$.
\endproof

\begin{assum}\label{ass:Ndiss}
There exists a nonempty set $N\subset\X$ such that for any $\psi:\X\to\R$ with $\psi\ge - \lambda$ and $\psi(x) = - \lambda(x)$ for all $x\in N$ we have that $T\psi(x) = c - \lambda(x)$ for all $x\in N$.
\end{assum}

Somewhat similar to Assumption \ref{ass:N}, for dissipative optimal control problems Assumption \ref{ass:Ndiss} holds with $N=\{x^e\}$ for an equilibrium $(x^e,u^e)$ with $\ell(x^e,u^e)=c$. This is because dissipativity implies $\tilde \ell \ge 0$ and Assumption \ref{ass:Ndiss} implies $\ell(x^e,u^e)=c$ implies $\tilde\ell(x^e,u^e)=0$. Together this yields for all $u\in\U(x^e)$ that
\[ \ell(x^e,u) + \psi(f(x^e,u)) \ge \ell(x^e,u) - \lambda(f(x^e,u)) = \tilde\ell(x^e,u) + c - \lambda(x^e) \ge c - \lambda(x^e),\]
while for $u=u^e$ we get
\[ \ell(x^e,u^e) + \psi(f(x^e,u^e)) = c -  \lambda(x^e), \]
implying that this is the minimum and hence $T\psi(x^e) = c -  \lambda(x^e)$. The situation just described in particular occurs for {\em strictly} dissipative problems, cf.\ eq.\ \eqref{eq:sdiss}.
\newline
\newline
\begin{thm} Assume that the optimal control problem is dissipative in the sense of Assumption \eqref{ass:diss}, that Assumption \ref{ass:Ndiss} holds and that there is $M>0$ with $T^k(\psi^0)\le M + ck$ for all $k\in\N$ and $\psi^0 = -\lambda$. Then the sequence of functions $\psi^k=(\hat T^k \psi^0)^n$, $k\in\N$, converges to a solution of the shifted Bellman Equation.
\end{thm}
\emph{Proof.} The assumptions together with Lemma \ref{lem:tildeT} imply that the operator $\hwT$ corresponding to the cost $\tilde\ell$ from \eqref{eq:tell} satisfies all assumptions of Proposition \ref{prop:hatconv}. Hence, for $\tilde\psi^0\equiv 0$ the sequence $\tilde \psi^k=(\hwTk \tilde \psi^0)^n$ converges to a solution $\tilde\psi^\infty$ of the Bellman Equation for $\tilde \ell$, i.e., $\wT\tilde\psi^\infty=\tilde\psi^\infty$. Because of Lemma \ref{lem:tildeT} and using that $(\psi+\phi)^n = (\psi^n+\phi)^n$ we obtain that 
\[ \psi^k = (\hat T^k (\tilde\psi^0 - \lambda))^n = 
(\hwTk (\tilde\psi^0) - \lambda)^n = ((\hwTk (\tilde\psi^0))^n - \lambda)^n = 
(\tilde \psi^k - \lambda)^n \] 
implying that 
\[ \psi^\infty = (\tilde \psi^\infty - \lambda)^n. \]
From this we get, again using Lemma \ref{lem:tildeT} and $w:=(\tilde \psi^\infty - \lambda)^n-\tilde \psi^\infty - \lambda$,
\begin{eqnarray*} T\psi^\infty & = & T (\tilde \psi^\infty - \lambda)^n \;
 = \; T (\tilde \psi^\infty - \lambda + w)\\
& = & T (\tilde \psi^\infty - \lambda) + w 
\; = \; \wT \tilde\psi^\infty - \lambda + w + c\\
& = & \tilde\psi^\infty - \lambda + w + c
\; = \; (\tilde\psi^\infty - \lambda)^n + c \; = \; \psi^\infty + c.
\end{eqnarray*}
This finishes the proof.
\endproof
\\

\section{Examples and Counterexamples}
\label{exandcex}

In this section we illustrate the performance of the iterations proposed and discussed in this paper with various examples.

\subsection{Comparison of solution methods}\label{sec:compex}
The examples in Section \ref{sec:compex} are meant to illustrate different approaches for the formulation and solution of infinite horizon optimal control problems using dynamic programming. In particular, they emphasize the need for a terminal penalty function and highlight the benefits of using the $\hat{T}$ and $\check{T}$ operators for their solution.
\subsubsection{Need for terminal penalty function}
\begin{figure}[htb]
\centerline{
\includegraphics[width=8cm]{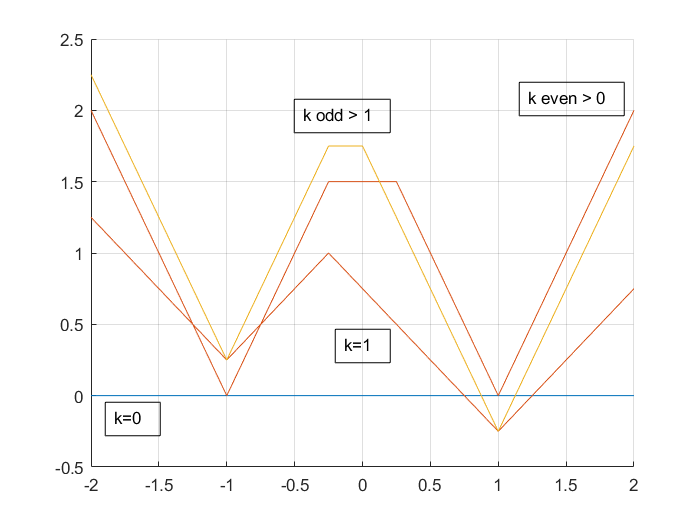} }
\label{nopenalty}
\caption{Sequence of cost-to-go functions $T^k \psi$, from $\psi=0$. }
\end{figure}

We consider the following scalar linear system:
\begin{equation}
\label{scalarlinear}
x^+ = - x + u
\end{equation}
along with state $x$ taking values in $\mathbb{X} = [-2,2]$, and input constraints $\mathbb{U} (x) = [-2+x,2+x]$.
The stage cost is piecewise linear and defined as:
\begin{equation}
\ell (x,u) = \min \left  \{ |x-1| - \frac{1}{4}, |x+1| + \frac{1}{4} \right \} + |u|.
\end{equation}
 Notice that the state-dependent part of the cost has two local minima, at $x$ equal $-1$ and $+1$. Moreover, for $u=0$ solutions are $2$-periodic and fulfill $x(t)= (-1)^t x(0)$.
It is possible to show that the optimal average cost is $0$, achieved by the solution $x(t) = (-1)^t$ corresponding to $u(t)=0$.
We show that using $\psi=0$ does not lead to a convergent sequence of cost-to-go functions. See Fig. \ref{nopenalty}.
 In particular, $T^k \psi$ converges to a period $2$ oscillation between two distinct piecewise linear functions after $2$ iterations.
Accordingly the optimal state-feedback (which is bang-bang) does not converge and will differ at least in some regions of state-space depending on whether an horizon of odd or even length is considered.

\begin{figure}[htb]
\centerline{
\includegraphics[width=8cm]{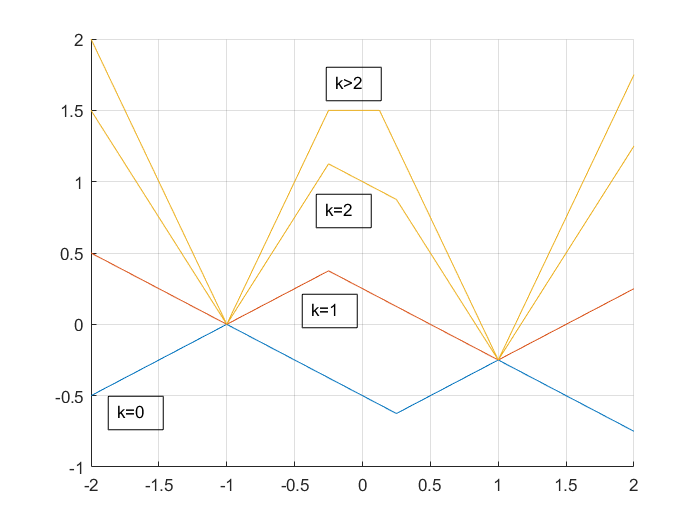} 
\includegraphics[width=8cm]{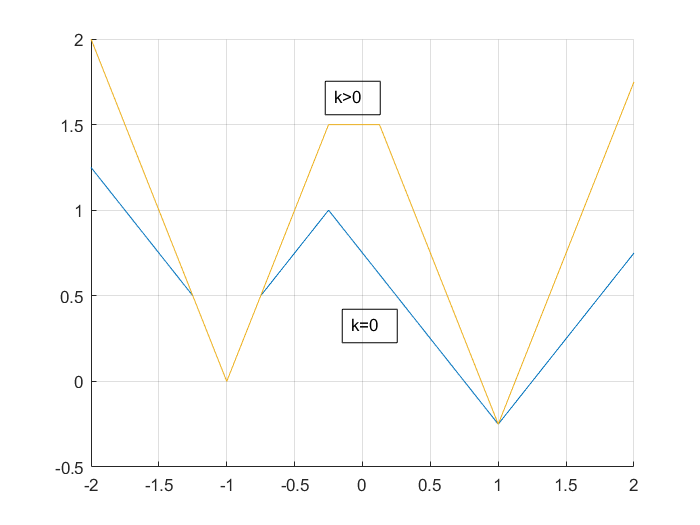} }
\caption{Sequence of cost-to-go functions $T^k \psi$, from $\psi=-\lambda_1$ (left) and $\psi=-\lambda_1$ (right) }
\label{withpenalty}
\end{figure}

In order to obtain meaningful infinite horizon costs and feedback policies we need to use a suitable penalty function for the final state. In particular by letting $\psi=- \lambda$ where $\lambda$ is a storage function.
For the considered example one can show that the function:
\[ \lambda_1(x) = \min \left \{ |x-1|+ \frac{1}{2}, |x+1| \right \} / 2 \]
is a storage function. Fig.\ \ref{withpenalty}(left) shows that the iteration initialized with $\psi=\lambda_1$ converges.
Notice that the cost monotonically converges in $3$ steps to its infinite horizon value. It is well known that storage functions need not be unique. 
For instance the following function is another storage function:
\[   \lambda_2 (x) = - \min \left \{ |x+1| +\frac{1}{4}, |x-1| - \frac{1}{4}, 2 |x+1|    \right \} \]
Our results show that any storage function can be used in order to define a suitable infinite horizon cost, provided this exists finite.
We show in Fig.\ \ref{withpenalty}(right) how choosing a different penalty function $\psi=-\lambda_2$ still leads, for this particular example, to the same infinite horizon cost, with convergence in just one time step.
\subsubsection{Solution with use of $\hat{T}$ operator}
\begin{figure}[htb]
\centerline{
\includegraphics[width=8cm]{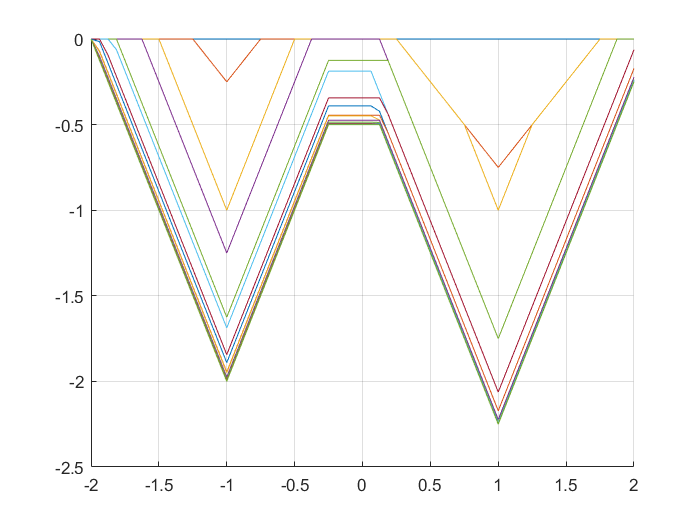}
\includegraphics[width=8cm]{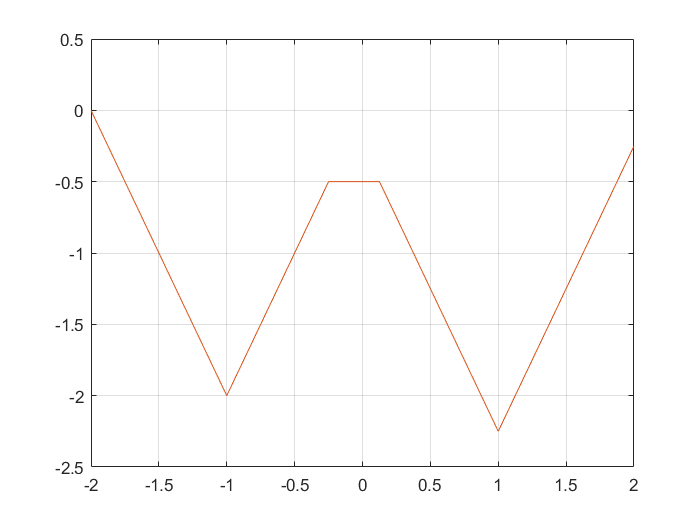} }
\caption{Sequence $\hat{T}^k \psi$ from initialisation $\psi=0$ (left) and limiting function (right)}
\label{that}
\end{figure}

\begin{figure}[htb]
\centerline{
\includegraphics[width=8cm]{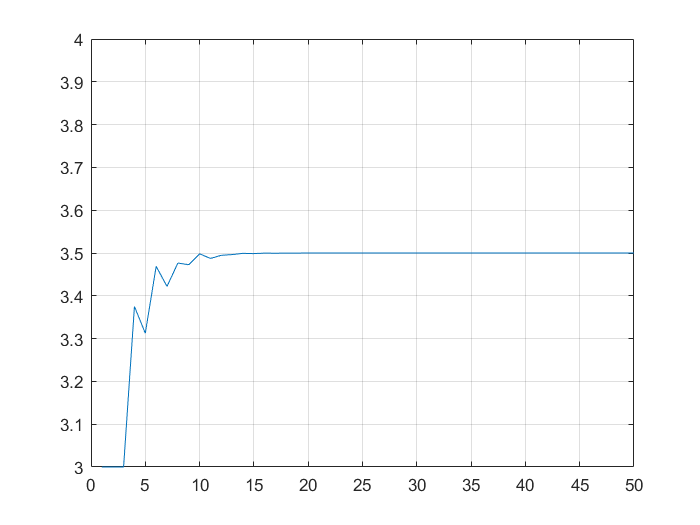} }
\caption{Sequence of shifts $c (\hat{T}^k \psi,T \hat{T}^k \psi)$.}
\label{translations}
\end{figure}

\begin{figure}[htb]
\centerline{
\includegraphics[width=8cm]{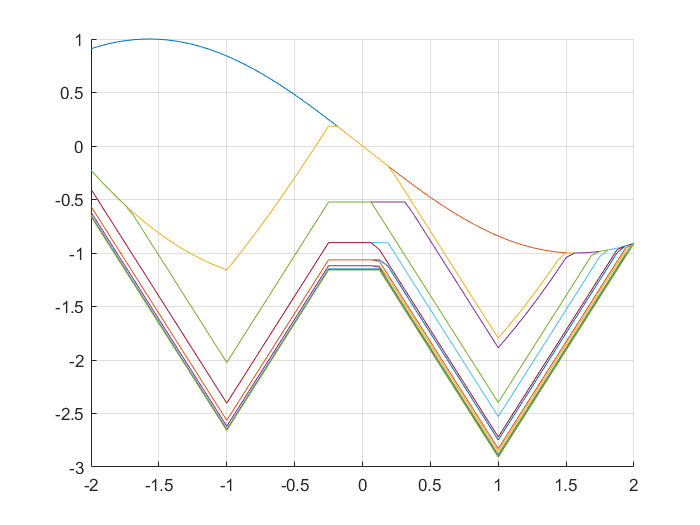} }
\caption{Sequence $\hat{T}^k \psi$ from initialisation $\psi=- \sin(x)$}
\label{sinite}
\end{figure}

We consider below the same system and constraints as in the previous example, namely
\begin{equation}
\label{scalarlinearagain}
x^+ = - x + u
\end{equation}
along with state $x$ taking values in $\mathbb{X} = [-2,2]$, and input constraints $\mathbb{U} (x) = [-2+x,2+x]$.
The stage cost is merely a shifted version of the previous piecewise linear cost:
\begin{equation}
\ell (x,u) = \min \left  \{ |x-1| - \frac{15}{4}, |x+1| - \frac{13}{4} \right \} + |u|.
\end{equation}
Rather then applying ad hoc considerations trying to figure out the optimal average performance (which in this case is $-7/2$) and correspondingly shifting $\ell$ in order to make the problem into its previous version with optimal $0$ average, we
directly apply the operator $\hat{T}$ to an arbitrary initialization $\psi(x)=0$.
We show in Fig.\ \ref{that}, the resulting non-increasing sequence of functions $\hat{T}^k \psi$, and the corresponding limit, which is a solution of the shifted Bellman Equation.
The value of shift applied $c (\hat{T}^k \psi,T \hat{T}^k \psi)$ is displayed in Fig. \ref{translations}. Notice that the shifts converge to $7/2$, which is indeed the positive translation needed in order to compensate for the optimal infinite horizon average performance of $-7/2$.
To highlight the power of the $\hat{T}$ iteration, which simultaneously adjusts to the right value of shift and asymptotic cost, we show in Fig. \ref{sinite} its evolution for a different initialisation $\psi(x)=- \sin(x)$.
\subsubsection{Solution with $\check{T}$ operator}
\begin{figure}[htb]
\centerline{
\includegraphics[width=8cm]{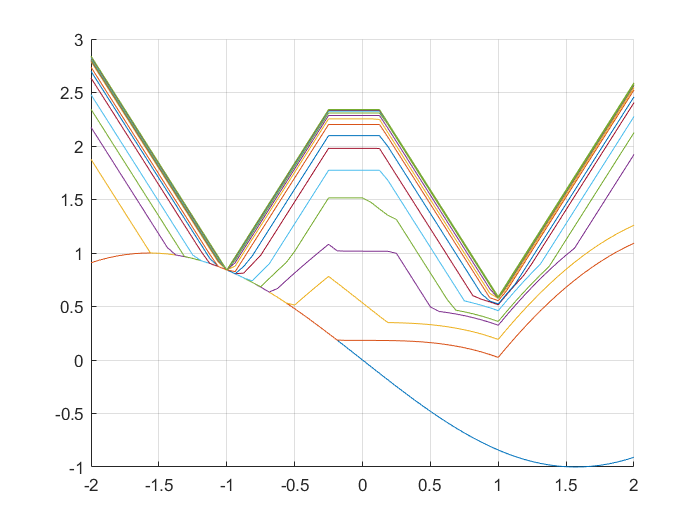} 
\includegraphics[width=8cm]{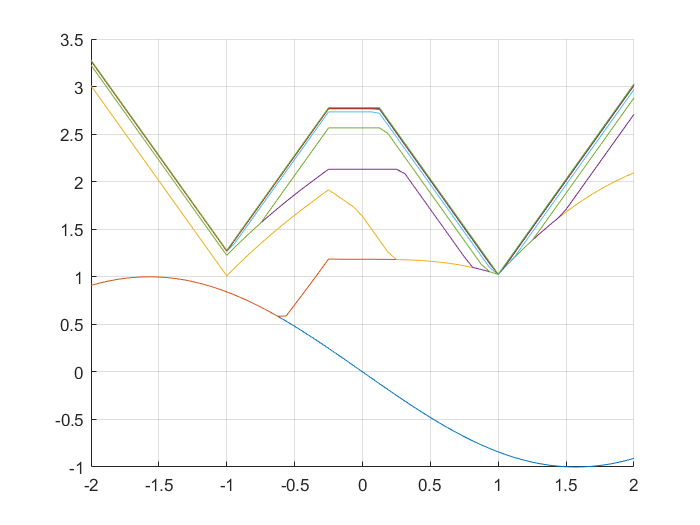} }
\caption{Sequence $\check{T}^k \psi$ from initialisation $\psi=-\sin(x)$ with shift term $c$ (left) and shift term $\tilde c$ with $\alpha=3/4$ (right)}
\label{tche}
\end{figure}

We provide next numerical evidence of convergence using the $\check{T}$ operator in Fig.\ \ref{tche}(left).
It is also interesting to remark that both $\hat{T}$ and $\check{T}$ operators show robustness with respect to the definition of the 
shift term $c( \psi, T \psi)$. Specifically, any strict convex combination ($\alpha \in (0,1)$ ):
\[ \tilde{c}( \psi_1, \psi_2) := \alpha \max_{x \in \mathbb{X}} \big [ \psi_1 (x) - \psi_2 (x) \big ]    + (1- \alpha ) \min_{x \in \mathbb{X}}  \big [ \psi_1 (x) - \psi_2 (x) \big ] \]
yields convergence, although at possibly different speed. 
To this end we show the iteration corresponding to $\alpha=3/4$ in Fig.\ \ref{tche}(right). 

\subsection{Non uniqueness of optimal solutions}
The following examples illustrate non-uniqueness phenomena arising when dealing with infinite horizon control problems. In particular, they emphasize non uniqueness of the fixed-points of the Bellman Equation and/or of the associated optimal feedback policies.
\subsubsection{Example with multiple solutions of the Bellman Equation}
\label{nondissipativeex}
Consider the scalar linear system:
\begin{equation}
x^+ = - x + u
\end{equation}
along with the state constraint: $\mathbb{X} = [-2,2]$ and input constraints $\mathbb{U} (x) = [-2+x,2 + x]$.
We consider a piecewise linear stage cost defined as:
\begin{equation}
\ell (x,u) = \varepsilon x + |u|
\end{equation}
for some constant $\varepsilon$ which will need to be sufficiently small.
Any function $\psi(x) = \alpha |x| + \varepsilon x/2$ is a solution of the (shifted) Bellman Equation, as long as $0 \leq \alpha< 1- \varepsilon$.
In fact:
\[ T \psi = \min_{u \in \mathbb{U} (x) } \varepsilon x + |u| + \alpha | -x + u | +  \varepsilon (u-x) /2 = \varepsilon x / 2 + \min_{u \in \mathbb{U}(x) } |u| + \alpha | -x + u | + \varepsilon u. \]
We notice that if $0 \leq \alpha<1- \varepsilon$ then the optimal value is achieved for $u=0$, since the slope of the absolute value of $|u|$ dominates the slope of the other terms.
In particular, substituting $u=0$ yields $T \psi (x) = \alpha |x| + \varepsilon x/2$.
Hence there are infinitely many (even continuous) solutions to the shifted Bellman Equation (\ref{shiftedBE}) (although the associated optimal feedback policies happen to be the same). We remark that because of Theorem \ref{thm:unique} this implies that the problem is not strictly dissipative. 
\subsubsection{Example with multiple optimal feedback policies}
\begin{figure}[htb]
\centerline{
\includegraphics[width=6cm]{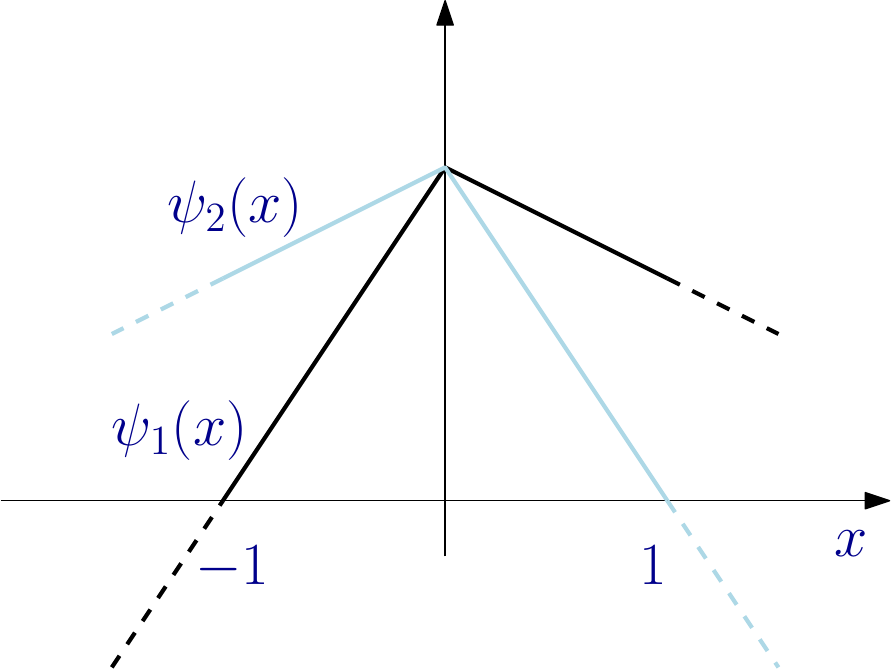} }
\caption{Multiple solutions of Bellman Equation}
\label{doubleweight}
\end{figure}

We consider the following scalar linear system:
\begin{equation}
x^+ = x + u
\end{equation}
along with the state constraint $\mathbb{X} = [-1,1]$ and input constraints $\mathbb{U} (x) =[-1-x,1-x]$.
We consider a piecewise linear stage cost defined as:
\begin{equation}
\ell(x,u) = 1 - |x| + |u|/2.
\end{equation}
Notice that, for each given $x \in \mathbb{X}$, $u=0$ minimizes the stage cost and makes $x$ into an equilibrium for the system. Hence, maximizing $|x|$ so as to minimize $\ell$, the optimal average performance is achieved for  the equilibrium solutions $x= \pm 1$ provided a zero input is applied.
Consider the following terminal penalty functions:
\begin{equation}
\begin{array}{rl}
\psi_1 (x) &= 1 - |x| + (1+x)/2 \\
\psi_2 (x) &= 1 - |x| + (1-x)/2
\end{array}
\end{equation}
As seen in Fig.\ \ref{doubleweight}, the functions $\psi_1$ and $\psi_2$ assign different terminal costs to the two optimal equilibria. In particular $\psi_1$ favours $-1$, with $0$ terminal cost, while $\psi_2$ favours
$+1$.

Both functions fulfill the Bellman Equation.
In fact:
\begin{eqnarray*}
\min_{u \in \mathbb{U} (x) } \ell(x,u) + \psi_1 (f(x,u)) & = & 
\min_{u \in [-1-x,1-x]} 1 - |x| + |u|/2 + [ 1 - |x+u| + (1+x+u)/2 \\
 & = & 1 - |x| + (1+x)/2,
\end{eqnarray*}
which is achieved for $u_1^*(x) = -1 -x$.
Similarly one can show that $u_2^*(x)=1-x$ achieves the optimum for $\psi_2$ and that $\psi_2$ is a solution of the Bellman Equation.
Notice that:
\[ \hat{\psi} (x) = \frac{3}{2} ( 1 -|x| ) = \min \{ \psi_1(x), \psi_2 (x) \} \]
is also a legitimate choice of terminal penalty function. In fact, this is the infimum element in $\Psi$, and is therefore the terminal penalty function that corresponds to the cheapest infinite horizon transient cost.
As shown in Proposition \ref{sameaverage}, feedback policies corresponding to different fixed-points of the shifted Bellman Equation, share the same infinite horizon average cost.
Notice, in addition, that for any constants $c_1$ and $c_2$, the
function:
\[ \psi(x) = \min \{ \psi_1 (x) +c_1, \psi_2 (x) + c_2 \}, \]
is a fixed-point of the shifted Bellman Equation.
In fact, in this case, it can be shown that every fixed point of the shifted Bellman Equation is of this form.
This result is likely to admit an extension to more general control set-ups.
\subsection{Regularity of fixed-points of Bellman Equation}
The following examples are meant to illustrate potential discontinuity and unboundedness issues of the fixed-point of the (shifted) Bellman Equation.
\subsubsection{Example with lower semi-continuous solution of the Bellman Equation}
\begin{figure}[htb]
    \centering
    \includegraphics[width=6cm]{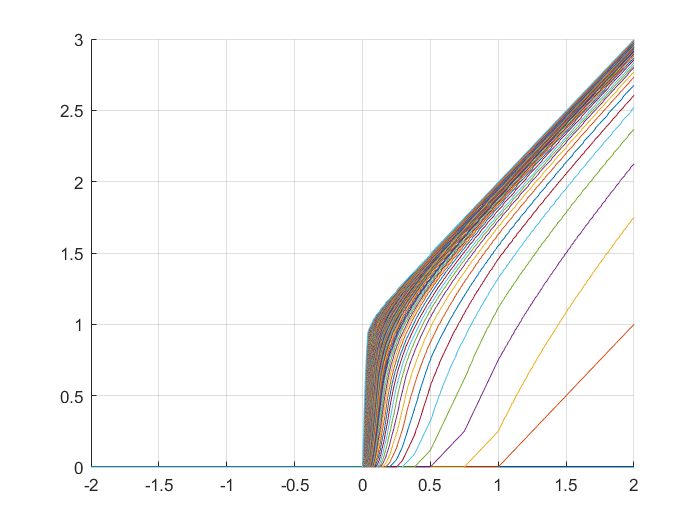}
    \includegraphics[width=6cm]{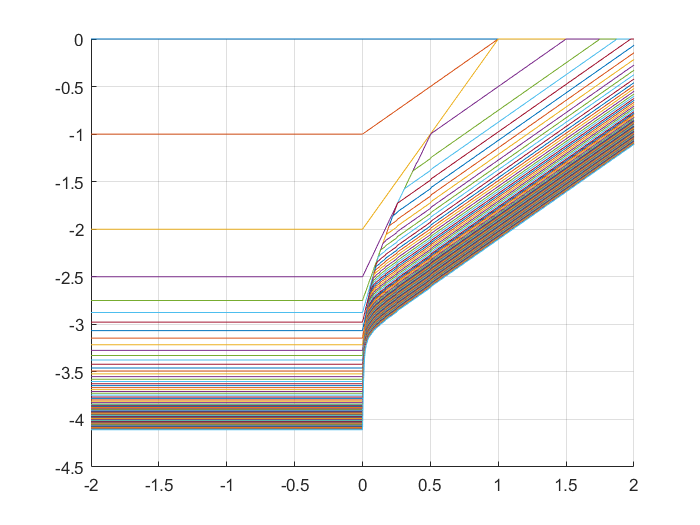}
    \caption{Iteration of $\check{T}^k \psi$ and $\hat{T}^k \psi$}
    \label{numLSC}
\end{figure}

Consider the following bilinear scalar system:
\begin{equation}
x^+ = x (1 + u )
\end{equation}
with state taking values in $\mathbb{X} =[-2,2]$ and input constraints:
\[  \mathbb{U} (x) = [ -2,0 ]. \]
Let the stage cost be piecewise linear defined according to:
\[ \ell(x,u) = \max \{ 0 , x \} + |u|. \]
Notice that for $u=0$ every point is an equilibrium. Hence, simply letting $u=0$ whenever the initial condition is $\leq 0$ achieves the minimum average cost. If the initial condition is positive, the 
best control action is $u=-1$. Indeed, an input $u \leq -1$ is needed in order to leave the set of positive states and enter the negative semi-axis, where the optimal average performance can be achieved. Hence, the best choice, given the penalty $|u|$ on inputs, is to have
$u=-1$. Moreover, waiting to apply such a control action does not pay off as the same cost will need to be incurred at some point in the future in order to switch to negative states.
The following function is a lower semi-continuous solution of the associated Bellman Equation:
\[  \psi (x) = \left \{ \begin{array}{cl} 0 & \textrm{if } x \leq 0 \\
                                 1 + x & \textrm{if } x > 0
\end{array}
\right .
\] 
which is achieved for the following control policy: 
\[   u^*(x)  = \left \{ \begin{array}{cl} 0 & \textrm{if } x \leq 0 \\
                                 -1 & \textrm{if } x > 0
\end{array}
\right . \]
We show in Fig. \ref{numLSC}, how the iterations of the operators $\hat{T}$ and $\check{T}$ behave when initialised from $\psi(x)=0$.

\noindent
It is worth pointing out that while both sequences seem to asymptotically approximate the correct `shape' of infinite-horizon cost, the theory confirms that $\hat{T}^k \psi (x)$ cannot be bounded, since its pointwise limit is known to be at least upper semi-continuous, which is not the case for the fixed point in the considered example.

\subsubsection{Example with unbounded infinite horizon cost}
Consider the following bilinear scalar system:
\begin{equation}
\label{biscalar}
x^+ = x u
\end{equation}
with state $x \in [0,1]:= \mathbb{X}$ and $u \in [1/2,1]$. Consider the stage cost:
\begin{equation}
\ell(x,u) = |u-1| + |x|.
\end{equation}
We claim that the optimal average cost is $0$. In fact,
the control sequence $u(t) = 1/2$ for $t = 0 \ldots K-1$ and $u(t) = 1$ for $t >K$ yields:
$x(K) = x(0)/2^K$, and $x(t)=x(K)$ for $t \geq K$. Notice that 
$\ell(x(t),u(t)) = |x(0)|/2^K <= 1/ 2^K$ for all $t \geq K$. Hence the  average cost can be made less or equal than $1/2^K$ for any positive integer $K$, and this, together with the inequality $\ell(x,u) \geq 0$, proves $0$ optimal average cost. 
We show next that the optimal cost is unbounded. \\
By induction, $x(k) = x(0) \prod_{t=0}^{k-1} u(t)$.
For the infinite horizon cost to be bounded we need to find an input such that $x(k) \rightarrow 0$ as $k \rightarrow + \infty$.
Hence, the input needs to fulfill $\prod_{t=0}^{k-1} u(t) \rightarrow 0$.
On the other hand:
\[ \prod_{t=0}^{k-1} u(t)= e^{ \sum_{t=0}^{k-1} \log ( u(t) ) } \]
and therefore, for the cost to be bounded we need:
\[ \sum_{t=0}^{k-1} \log ( u(t) ) \rightarrow - \infty \] 
as $k \rightarrow + \infty$.
However, on the interval $[1/2,1]$, concavity of the $\log$ function yields:
\[ \log (u ) \geq \log(2) ( u -1 ). \]
Using the inequality above shows:
\[ \sum_{t=0}^{k-1} \log ( u(t) ) \geq
\log(2) \sum_{t=0}^{k-1} (u(t) -1 ). \]
As a consequence, for the infinite horizon cost to be bounded we need:
\[ \sum_{t=0}^{k-1} (u(t) -1 ) \rightarrow - \infty, \]
as $k \rightarrow + \infty$. This, however, contradicts boundedness of the cost as $\ell (x,u) \geq 1 - u$. \\

It is worth pointing out that the optimal steady state for the considered example is $x_s=0$ and $u_s=1$.
This steady state is not reachable in finite time, though. Notice also that this is trivially a dissipative system with storage function $\lambda(x)=0$
due to the non-negativity of the cost.
 As a consequence no bounded fixed-point of the shifted Bellman Equation exists. 

\subsubsection{Example with continuous and discontinuous fixed points}
\label{notfixedpoint}
\begin{figure}[htb]
    \centering
    \includegraphics[width=8cm]{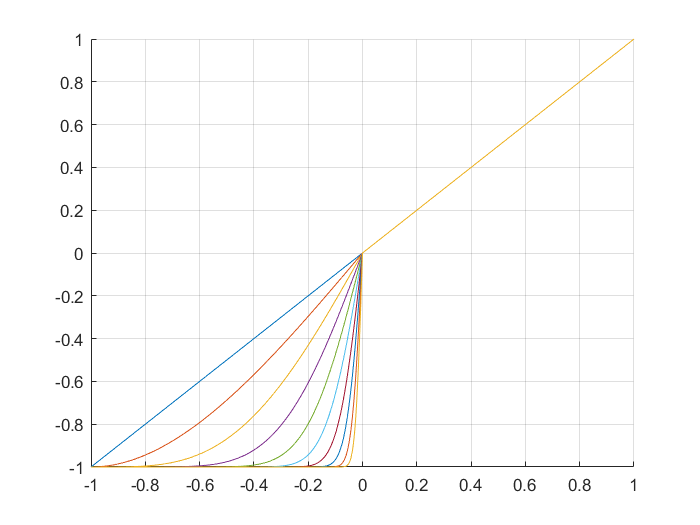}
    \caption{Numerical $\hat{T}^k x$ iteration}
    \label{figthat}
\end{figure}

Consider the autonomous nonlinear system:
\[  x^+ = \frac{3}{2} x - \frac{1}{2} x^3, \]
along with the cost functional $\ell(x,u)=0$.
Choose $\mathbb{X} = [-1,1]$ which is a forward invariant set for the dynamics, with $3$ equilibria in $-1,0$ and $1$ respectively.
The equilibrium in $0$ is antistable, while the equilibria in $\pm 1$ are asymptotically stable with basin of attraction $(0,1)$ and $(-1,0)$ respectively.
Clearly, $\psi(x) \equiv 0$ is a fixed point of the Bellman Equation. 
Any function of the form:
\[ \psi(x) = \left \{ \begin{array}{cl} c_{1} & x < 0 \\
c_2 & x = 0 \\ c_3 & x>0 \end{array} \right . \]
is also a fixed point. Consider next an arbitrary continuous increasing initialisation of $\psi$ of the $\hat{T}$ and $\check{T}$ maps.
It can be seen that $T \psi$ is also increasing, as $f(x)$ is such in the interval $[-1,1]$. As a consequence $\hat{T} \psi$ and $\check{T} \psi$ are also increasing.
Moreover, $T \psi (0) = \psi (0)$ and $T \psi ( \pm 1) = \psi (\pm 1)$.
Thus, $\hat{T} \psi (0) - \hat{T} \psi (-1) = \psi(0)- \psi (-1)$ and $\hat{T} \psi (1) - \hat{T} \psi(0) = \psi(1) - \psi(0)$.
By induction then, $\hat{T}^k \psi (x)$ is increasing with respect to $x$ for all $k$ and so is $\check{T}^k \psi(x)$. 
It can be shown that for $\psi(x)=x$ it holds $c(\hat{T}^k \psi,T \hat{T}^k \psi)=0$ for all $k$. In particular,
$\hat{T}^k \psi$ converges to:
\[ \hat{\psi} (x) = \left \{ \begin{array}{cl} -1 & x < 0 \\
 x & x \geq 0 \end{array} \right . \]
Numerical simulations indeed confirm this claim, see Fig.\ \ref{figthat}.
This shows that even if $\hat{T}$ (or $\check{T}$) admit continuous fixed points, the iteration of $\hat{T}^k \psi$ does not necessarily converge to a fixed point of the Bellman Equation.
Similarly, considering the iteration $\check{T}^k \psi$, for the same initial function $\psi(x)=x$, it holds $c(\check{T}^k \psi, T \check{T}^k \psi)=0$ for all $k$ and $\check{T}^k \psi$ converges to:
\[ \check{\psi} (x) = \left \{ \begin{array}{cl} 1 & x > 0 \\
 x & x \leq 0 \end{array} \right .\]
\subsubsection{Example with upper semi-continuous fixed point}
\begin{figure}[htb]
    \centering
    \includegraphics[width=8cm]{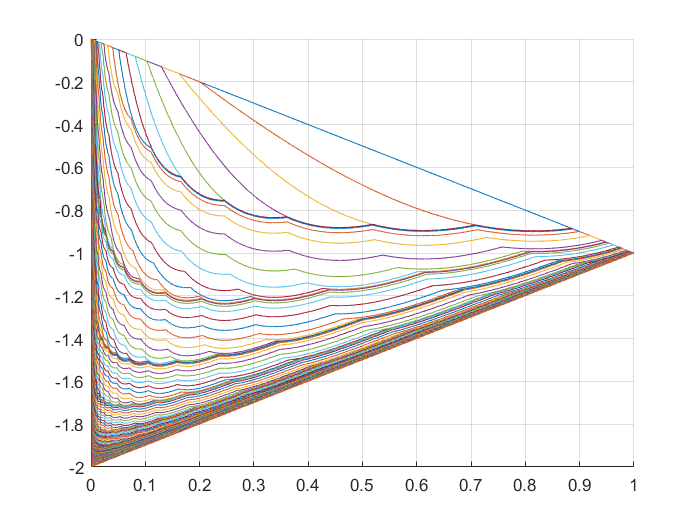}
    \caption{Numerical $\hat{T}^k (- x)$ iteration}
    \label{figusc}
\end{figure}

We slightly modify the previous example to include a scalar control input and induce an upper semi-continuous fixed point.
Consider the nonlinear system:
\[  x^+ = u \left (\frac{3}{2} x - \frac{1}{2} x^3 \right )=: f(x,u), \]
with state-space $\mathbb{X}=[0,1]$, scalar input $u$ constrained in $\mathbb{U}(x) = [0,1]$ along with the cost functional 
\[ \ell(x,u)= |u-1| - f(x,u) + x.\]
Notice that:
\[ \sum_{k=0}^{T-1} \ell(x(k),u(k)) = x(0) - x(T) + \sum_{k=0}^{T-1} |u(k)-1|. \]
Hence, the optimal average performance is $0$, achieved for $u(\cdot)=1$.
The function $\bar{\psi}$ defined below:
\[ \bar{\psi} (x) = \left \{ \begin{array}{rl} -2 + x & \textrm{for } x \in (0,1] \\ 0 & \textrm{for } x=0 \end{array} \right. \]
is a fixed point of the Bellman Equation.
To see this, notice, assuming $x \neq 0$:
\begin{eqnarray*}
T \bar{\psi} (x) & = & \min_{u \in [0,1]} |u-1| + x - f(x,u) + \bar{\psi} ( f(x,u) ) \\
& = & \min \left \{ 1 + x + \bar{\psi} (0),  \inf_{u \in (0,1]} |u-1| + x - f(x,u) + \bar{\psi} ( f(x,u) ) \right \} \\
& = & \min \left \{ 1 + x + \bar{\psi} (0),  \inf_{u \in (0,1]} |u-1| + x  - 2   \right \} \; = \; -2 + x.    
\end{eqnarray*}
For $x=0$, it is easy to verify $ T \bar{\psi} (0) = 0$. We show in Fig.\ \ref{figusc} the iteration converging to $\bar{\psi}$.
Notice that, despite $\bar{\psi}$ being upper semi-continuous, not admitting a minimum in $[0,1]$, and the discontinuity point $x=0$ being reachable from all states in $\mathbb{X}$ within a single step, still the minimum in the definition of the operator $T \bar{\psi}$ is achieved.
More in general we see that the iteration $\hat{T}^k \psi$ converges, for $x>0$, to $\psi(1)-1 + x$. \\
\emph{
We don't have any examples of optimal control problems where the only fixed points are upper semi-continuous (and not continuous), or where the minimum $T \bar{\psi}$ is not achieved. }
It is worth pointing out that $\bar{\psi}(x)=x$ is also a fixed point of the Bellman Equation.


\subsection{Complex optimal regime of operation}
We consider examples where the optimal average performance is not achieved at steady-state, but for more exotic type of behaviours. It is worth pointing out that dealing with a terminal penalty function allows to treat such examples without the need of an a priori known terminal absorbing state or terminal absorbing set. Moreover, the optimal regime of operation does not entail a constant (or zero) optimal stage cost in steady-state.
\subsubsection{
Example with chaotic optimal regime}
\begin{figure}[htb]
\centerline{
\includegraphics[width=6cm]{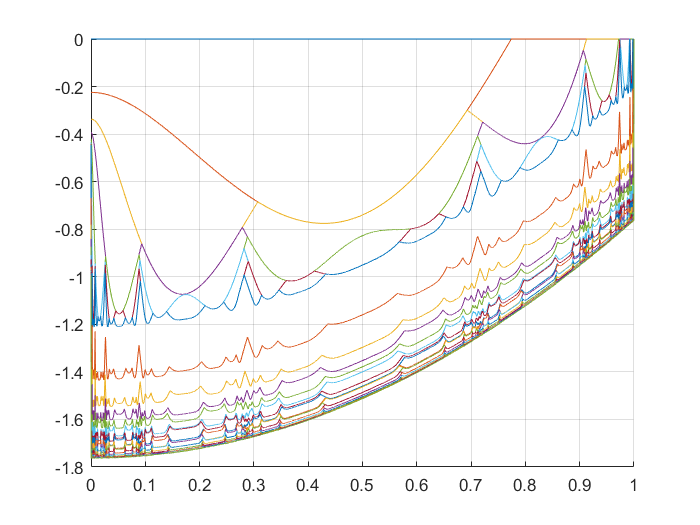}
\includegraphics[width=6cm]{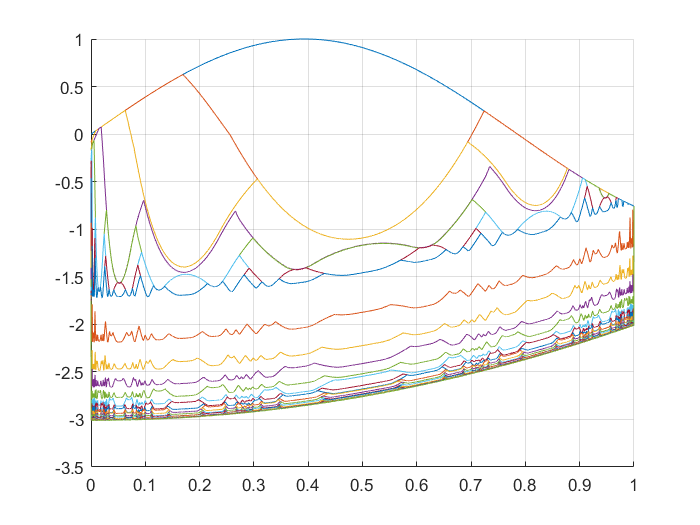} }
\caption{Iteration using the $\hat{T}$ operator from $\psi(x)=0$ and $\psi(x)= \sin(4x)$}
\label{chaos1min}
\end{figure}

\begin{figure}[htb]
\centerline{
\includegraphics[width=6cm]{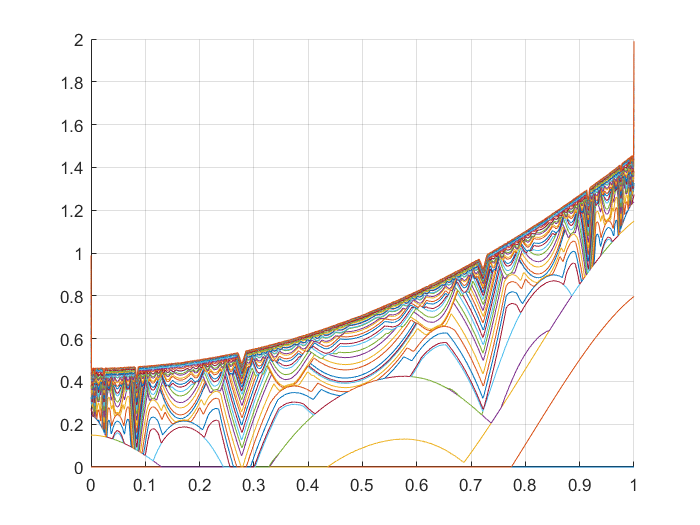}
\includegraphics[width=6cm]{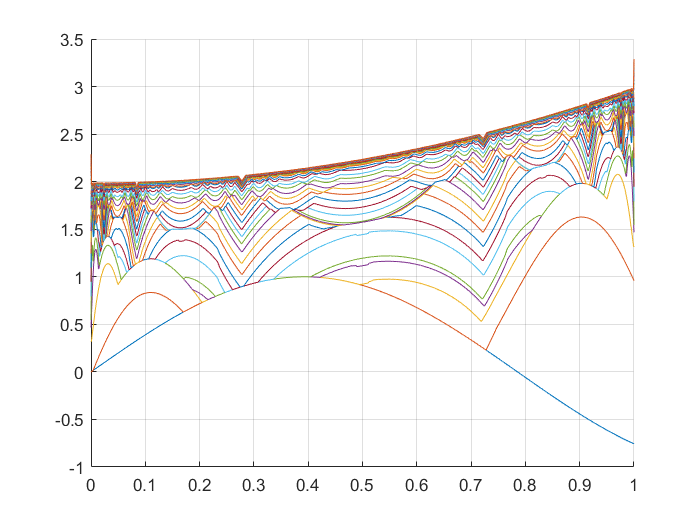} }
\caption{Iteration using the $\check{T}$ operator from $\psi(x)=0$ and $\psi(x)= \sin(4x)$}
\label{chaos2}
\end{figure}

Consider the scalar nonlinear system:
\begin{equation}
\label{logisticmap}
x^+ = u x (1-x)
\end{equation}
with scalar state $x \in \mathbb{X} := [0,1]$ and input $u \in \mathbb{U} : = [0,4]$.
We consider the stage-cost:
\[ \ell(x,u) = x^2 - [u x (1-x)]^2 + |u-18/5|. \]
Notice that $\ell (x(k),u(k)) = x(k)^2 - x(k+1)^2 + |u(k)-18/5|$.
Therefore, along arbitrary solutions we have:
\[  \sum_{k=0}^{T-1} \ell(x(k),u(k)) = x(0)^2 - x(T)^2 + \sum_{k=0}^{T-1} |u(k)-18/5|. \]
In particular, computing asymptotic time averages we see:
\[    \lim_{T \rightarrow + \infty}  \frac{ \sum_{k=0}^{T-1} \ell(x(k),u(k))}{T} = \lim_{T \rightarrow + \infty} \frac{ x(0)^2 - x(T)^2 + \sum_{k=0}^{T-1} |u(k)-18/5|}{T} \]
\[ \qquad =  \lim_{T \rightarrow + \infty} \frac{ \sum_{k=0}^{T-1} |u(k)-18/5|}{T}. \]
The optimal average performance is therefore $0$, and is achieved for instance, for any input $u(k)$ converging to $18/5$. Notice that for $u=18/5$ the considered dynamical system is known to have chaotic solutions.
Moreover $u(k) = 18/5$ is potentially an optimal infinite horizon control policy.
This policy corresponds to the fixed point $\bar{ \psi } (x) = x^2$ of the Bellman Equation.
Indeed,
\[ T \bar{ \psi} = \min_{u \in [0,4]  } x^2 - [u x (1-x)]^2 + |u-18/5| + \bar{ \psi } ( u x (1-x) ) \] 
\[ \qquad \qquad = \min_{u \in [0,4]} x^2 + |u-18/5| = x^2 = \bar{ \psi} (x). \]
Numerical solution using the $\hat{T}$ operator is shown in Fig.
\ref{chaos1min}, starting from two distinct initializations, $\psi(x) =0$ and $\psi(x) = \sin(4x)$.
The optimal average performance is correctly estimated to be $0$ and $\hat{T}^k \psi$ converges to a shifted version of $x^2$ in both cases.
The numerical solution using the $\check{T}$ operator is slightly different and is shown in Fig. \ref{chaos2}.
\begin{figure}[h!]
\centerline{
\includegraphics[width=8cm]{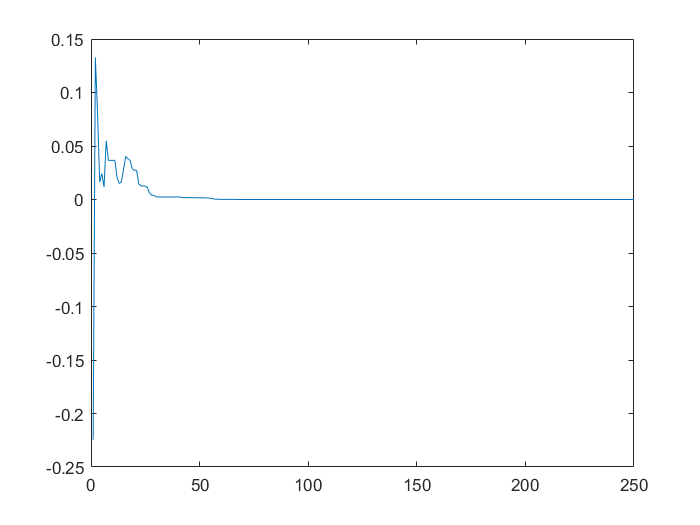}
}
\caption{Sequence $c ( \check{T}^k \psi, T \check{T}^k \psi)$, for $\psi=0$ }
\label{chaos3}
\end{figure}
While it is hard to write an explicit analytic solution of the limiting function, due to the presence of somewhat unexpected spikes, we believe that the numerical result hint at the presence of multiple solutions to corresponding the Bellman Equation.
These solutions match $x^2$ for most of the interval $[0,1]$ but appear to allow for piecewise linear spikes that might correspond to transient costs in regions which are not visited by the chaotic attractor.
It seems more plausible that these be true solutions rather than artifacts due to numerical approximations.
The optimal average performance is identified with very good precision in both cases. In particular, for the $\hat{T}$ iteration the error is lower than $10^{-16}$.
See Fig.\ \ref{chaos3} for the shift sequence achieved for the $\check{T}$ operator when $\psi(x) \equiv 0$.

\subsubsection{Two-dimensional example with periodic optimal regime}
\begin{figure}[htb]
\centerline{
\includegraphics[width=8cm]{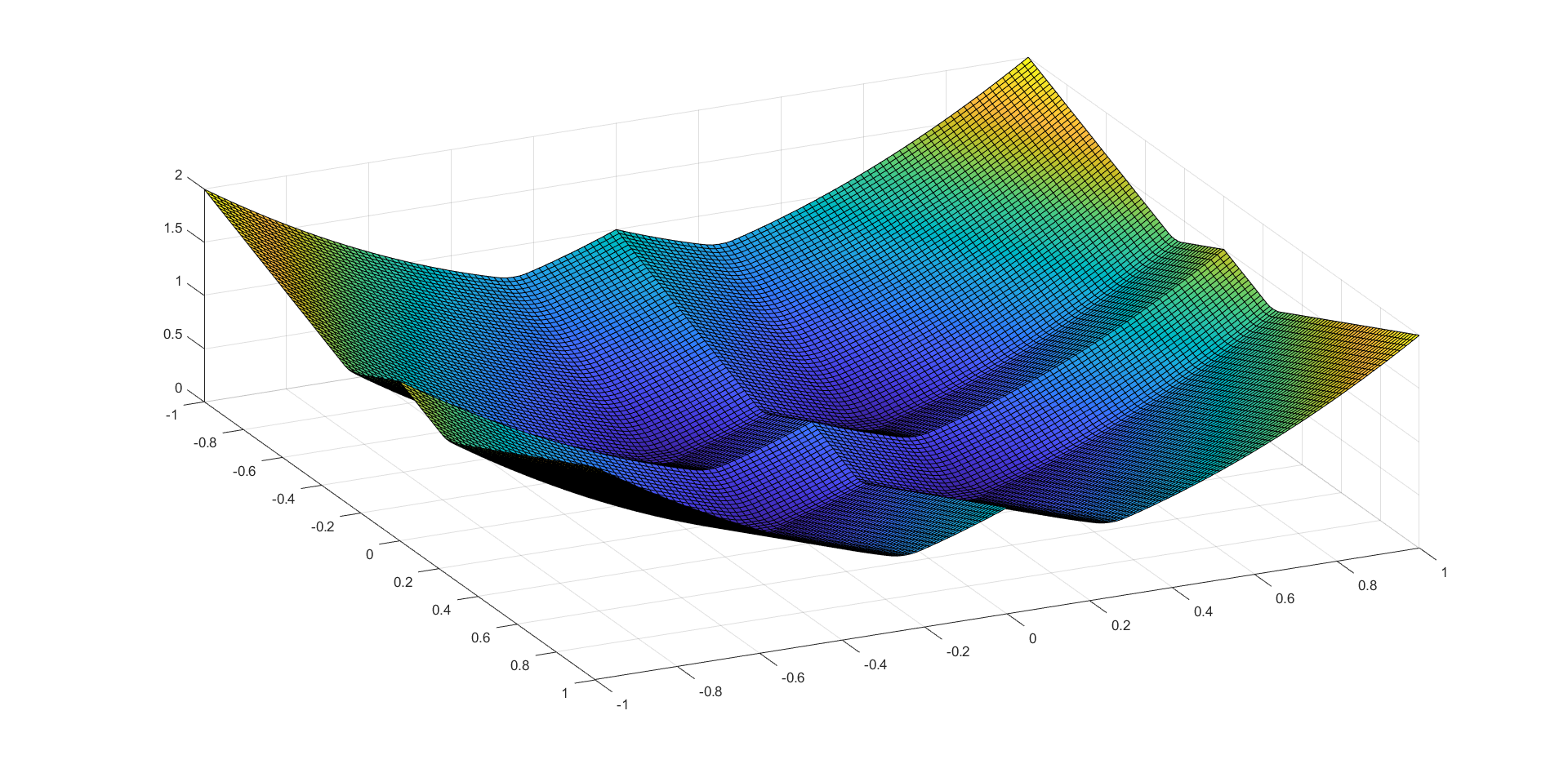} }
\caption{Fixed point of the 2d Bellman Equation}
\label{2dfirst}
\end{figure}

\begin{figure}[htb]
\centerline{
\includegraphics[width=8cm]{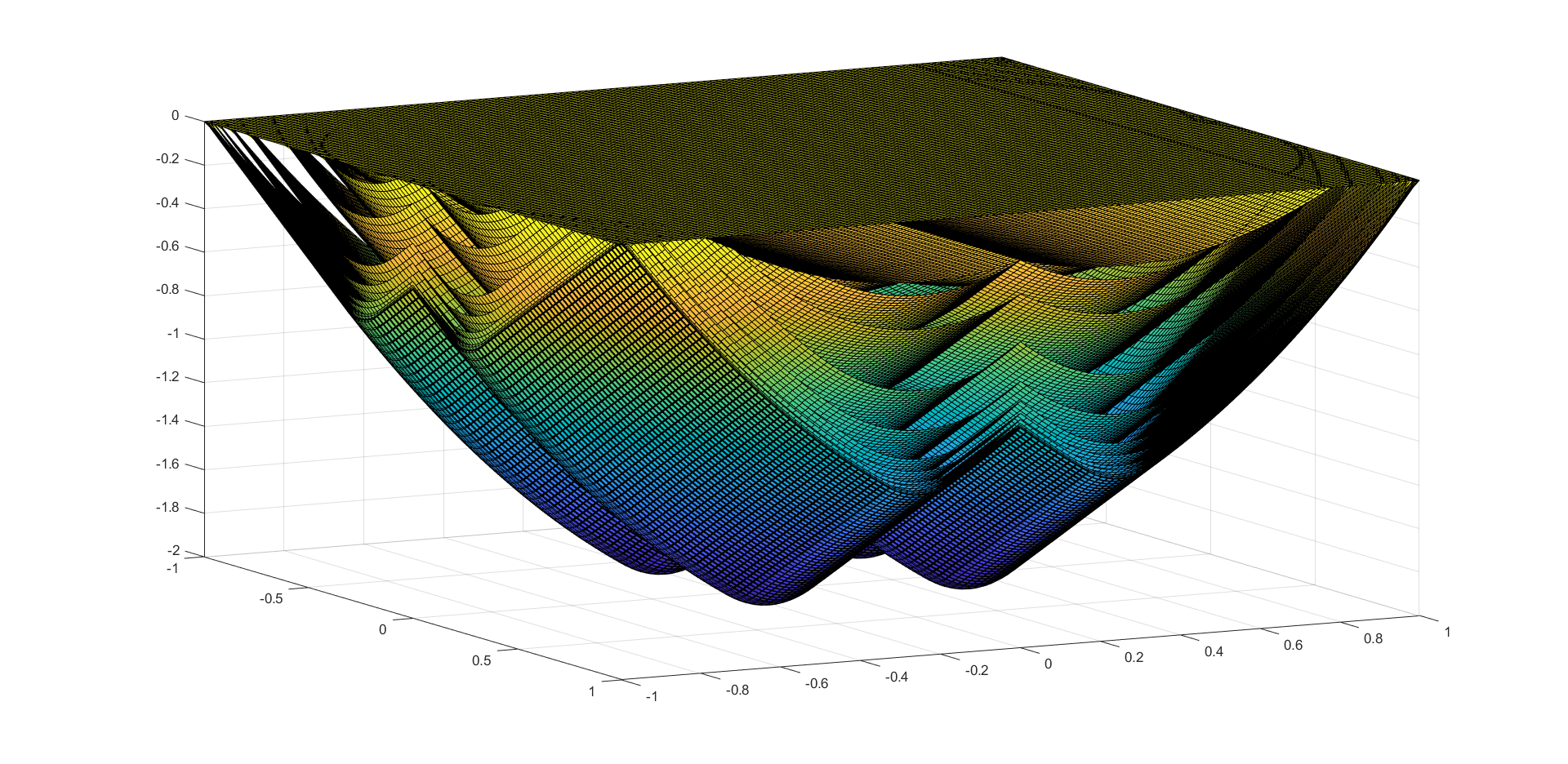} 
\includegraphics[width=8cm]{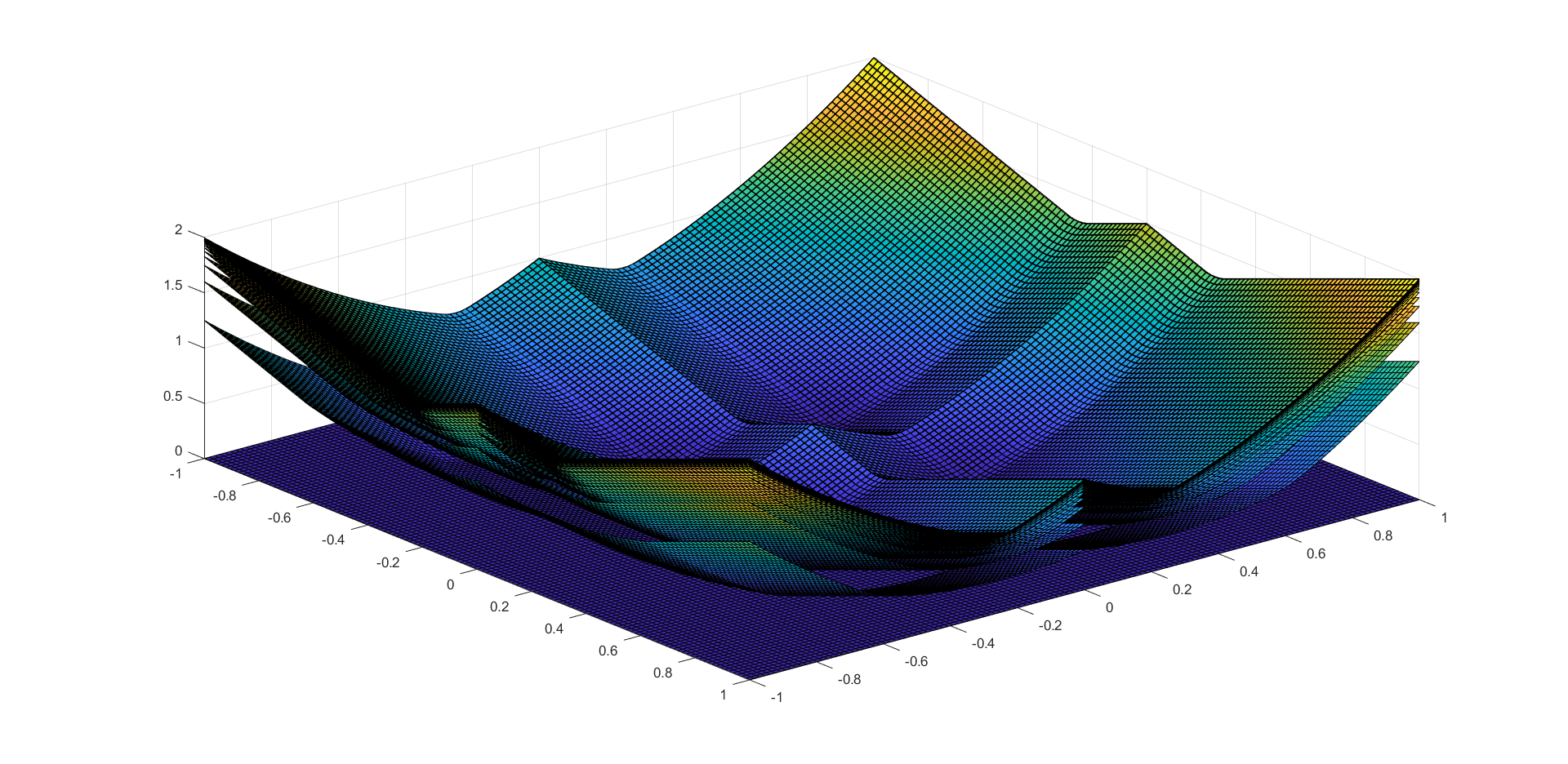} }
\caption{Iterations of the 2d $\hat{T}$ operator (left) and of the 2d $\check{T}$ operator (right)}
\label{2dsecond}
\end{figure}

We consider next the following two-dimensional linear system:
\begin{equation}
\label{2dsys}
x+ = \left [ \begin{array}{cc} 0 & 1 \\ -1 & 0 \end{array} \right ] \, x + \left [ \begin{array}{c} 1 \\ 0 \end{array} \right ] \, u,
\end{equation}
with state $x \in \mathbb{X}:= [-1,1]^2$, and input $u \in \mathbb{U} (x) := [-1-x_2,1-x_2]$.
Consider the stage-cost
\[ \ell(x,u)= |u| + x_1^2 - |x_1|/2. \]
Notice that this cost is not positive definite. In particular, the optimal average performance can be expected to be negative, as the zero solution is feasible with zero input, yielding $0$ average cost. However,  the stage cost can be made negative for some values of $x_1 \neq 0$.  
The zero-input responses of the system are (feasible) period $4$ oscillations. Moreover the system is controllable, which guarantees an optimal average performance independent of the initial condition (and regardless of the adopted stage cost $\ell(x,u)$).
We show in Fig. \ref{2dfirst} a fixed point of the shifted Bellman Equation.
The iterations resulting from the $\hat{T}$ operator and the $\check T$ operator starting in $\psi^0\equiv 0$ are shown in Fig.\ \ref{2dsecond}.

\subsection{Inefficiency of exponential discounting factors}
\label{motivatingexample}
\begin{figure}[htb]
    \centering
    \includegraphics[width=8cm]{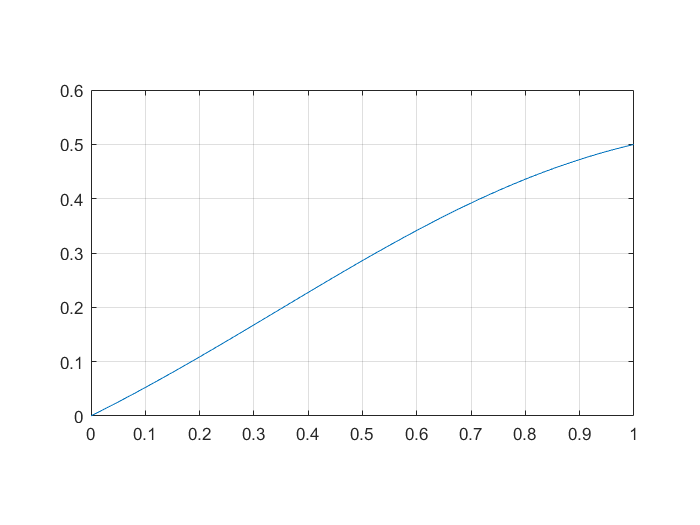}
    \caption{Equilibrium $x_{e}$ as a function of $\gamma$ in $[0,1]$}
    \label{xeqga}
\end{figure}

We end our example section with an example of a discounted optimal control problem, which shows that ensuring well-posedness of infinite horizon optimal control problems by means of discounting can have unwanted side effects, making the proposed approach via the shifted Bellman Equation an attractive alternative. To this end,
we consider a scalar infinite horizon linear quadratic optimal control problem with exponential discounting. In particular, the system's dynamics are given as:
\begin{equation}
    x^+ = (x+u)/2,
\end{equation} 
with $x$ and $u$ taking values in $\mathbb{R}$.
The stage cost is:
\[ \ell(x,u) = (x-1)^2 + u^2. \]
Since this choice will not give rise to bounded costs over an infinite horizon we use a discounting factor $\gamma \in (0,1)$:
\[ J_\gamma = \sum_{k=0}^{+\infty} \gamma^k \ell (x(k),u(k)). \]
The optimal infinite horizon cost fulfills the following Bellman Equation:
\[ J^*_{\gamma} (x) = \min_{u \in \mathbb{R}} \ell(x,u) + \gamma J^*_{\gamma} ( f(x,u) ). \]
It is possible to show that this equation admits a solution:
\[ J^*_{\gamma} (x) = \alpha x^2 + \beta x + \delta \]
where $\alpha$, $\beta$ and $\delta$ fulfill the conditions:
\[
\begin{array}{rcl}
\alpha (\gamma ) &=& \gamma - 2 + \sqrt{ \gamma^2 + 4 } \\
\beta (\gamma) &=& - \frac{2 \alpha(\gamma) \gamma + 8}{ \alpha(\gamma) \gamma + 4 - 2 \gamma} \\
\delta (\gamma) &=& \frac{4 \alpha(\gamma) \gamma + 16 - \beta^2(\gamma) \gamma^2   }{(4 \alpha(\gamma) \gamma + 16)(1 - \gamma)}
\end{array}
\]
The optimal feedback is affine in $x$ and expressed as:
\[  u^*(x) = - \frac{ \beta \gamma + \alpha \gamma x }{ \alpha \gamma + 4}. \]
This feedback globally asymptotically stabilizes a unique equilibrium $x_e(\gamma)$:
\[ x_e (\gamma)  = - \frac{ \beta(\gamma) \gamma }{2 \alpha(\gamma) \gamma +4 }. \]
Notice that the optimal average performance is achieved at equilibrium, for $x=1/2$ and $u=1/2$, which yields $V^{avg} = (1/2)^2 + (1/2)^2 = 1/2$.
On the other hand, the equilibrium $x_e (\gamma)$ only approaches the value $1/2$ as $\gamma \rightarrow 1$ (see Fig. \ref{xeqga}). This shows that the long run average performance achieved by introducing a discounting factor is in general suboptimal. Moreover, the discounting factor introduces a non existent trade-off between optimising transient cost and steady-state (average) costs which persist for $\gamma$ arbitrarily close to $1$. This trade-off can be avoided by the approach pursued in this paper.
On the other hand, any feedback $u=k(x)$ (for instance affine, $u= k_1 x + k_2$) which stabilizes the equilibrium
$1/2$, clearly achieves optimal average performance (and is therefore optimal with respect to the cost functional
$J^{\textrm{avg}}$), but, at the same time, it is not necessarily optimal from the point of view of transient costs.  We refer to \cite{GGHKW18,GruK21} for more examples of this kind and an in-depth study of the stability properties of discounted optimal equilibria.

\section{Conclusions and outlook}
\label{finalsection}
Two novel recursion operators are proposed for the simultaneous computation of value functions and minimal average asymptotic cost in discrete-time infinite horizon optimal control problems. The recursive formulas can be readily applied when average asymptotic cost is independent of initial conditions, a situation referred to as the \emph{ergodic case} in \cite{gaitsgory1}. 
The approach renders dynamic programming techniques invariant with respect to additive constants on the stage cost, as it is naturally the case in the finite horizon case, for infinite horizon control problems. 
The recursions converge, under fairly relaxed technical assumptions, to fixed-points of a shifted Bellman Equation, whose shift value is not a priori determined but is asymptotically computed alongside the value function. The approach removes the need for \emph{absorbing states} and zero cost conditions on the absorbing sets which have often hindered the applicability such techniques, or the need for discounting factors which introduce unnecessary trade offs between transient cost and asymptotic average performance.
While the approach is developed for the case of deterministic systems only, its extension to stochastic settings appears of potential interest. Finally, this may serve as a first step in understanding the more general question of a shift-invariant approach to infinite horizon optimal control problems in the \emph{non-ergodic} case, \cite{gaitsgory1, borkar}.

\newpage
 
\appendix

\section{Appendix: Technical results}
 \label{preliminaryresults}
 In order to analyse the convergence properties of the newly introduced operators $\check{T}$ and $\hat{T}$ it is useful to explore inequalities involving the $\max$ and $\min$ operators applied to a finite set of functions. The next two lemmas provide such tools.
\begin{lem}
\label{boundtop1}
Let $\psi_i \in \mathcal{C}( \mathbb{X} )$ for $i=1,\ldots,N$.
Then, the following holds:
\[  \max_{x \in \mathbb{X}}  \max_{i \in \{1, \ldots,N \} } \psi_i (x)  -  \min_{x \in \mathbb{X}}  \max_{i \in \{1, \ldots,N \} }  \psi_i (x)  \leq \max_{i \in \{1, \ldots, N \} } \left \{ \max_{x\in \mathbb{X}} \psi_i (x) - \min_{x \in \mathbb{X}} \psi_i (x) \right \}. \]
 \end{lem}
\emph{Proof.} Let $x^*$ in $\mathbb{X}$ be such that:
\[ \psi_{\bar{\imath}} (x^*) = \max_{x \in \mathbb{X}} \max_{ i \in \{1, \ldots, N \} } \psi_i (x), \]
 for some $\bar{\imath}$ in $\{1, \ldots, N\}$.
By monotonicity of the $\min$ operator, we see that:
\[ \min_{x \in \mathbb{X}} \max_{i \in \{1, \ldots, N \} } \psi_i (x) \geq \min_{x} \psi_{\bar{\imath}} (x). \]
Combining the latter inequality with the previous equality yields:
\[ \max_{x \in \mathbb{X}}   \max_{i \in \{1, \ldots,N \} }  \psi_i (x)  -  \min_{x \in \mathbb{X}} \max_{i \in \{1, \ldots, N \} } \psi_i (x) \]
\[ \qquad \leq \psi_{\bar{\imath} } (x^* ) - \min_{x \in \mathbb{X}} \psi_{\bar{\imath} } (x) = \max_{x \in \mathbb{X}} \psi_{\bar{\imath}} (x) -  \min_{x \in \mathbb{X}} \psi_{\bar{\imath} } (x) \]
\[ \qquad \qquad \leq \max_{i \in \{1, \ldots, N \} }
\left \{  \max_{x \in \mathbb{X}} \psi_{i} (x) -  \min_{x \in \mathbb{X}} \psi_{i } (x) \right \}. \]
\endproof

The following lemma provides a similar bound for the $\min$ operator.  
\begin{lem}\label{boundtopbis}
Let $\psi_i$, be continuous functions of $x \in \mathbb{X}$, for $i=1,\ldots,N$.
Then the following holds:
\[  \max_{x \in \mathbb{X}}  \min_{i \in \{1, \ldots,N \} }  \psi_i (x)  -  \min_{x \in \mathbb{X}}  \min_{i \in \{1, \ldots,N \} }  \psi_i (x)  \leq \max_{i \in \{1, \ldots, N \} } \left \{ \max_{x \in \mathbb{X}} \psi_i (x) - \min_{x \in \mathbb{X}} \psi_i (x) \right \}. \]
 \end{lem} 
\emph{Proof.} Let $x^*$ in $\mathbb{X}$ be such that:
\[ \psi_{\bar{\imath}} (x^*) = \min_{x \in \mathbb{X}} \min_{ i \in \{1, \ldots, N \} } \psi_i (x), \]
 for some $\bar{\imath}$ in $\{1, \ldots, N\}$.
By monotonicity of the $\max$ operator, we see that:
\[ \max_{x \in \mathbb{X}} \min_{i \in \{1, \ldots, N \} } \psi_i (x) \leq \max_{x \in \mathbb{X}} \psi_{\bar{\imath}} (x). \]
Combining the above inequalities imply:
\[ \max_{x \in \mathbb{X}}   \min_{i \in \{1, \ldots,N \} }  \psi_i (x)  -  \min_{x \in \mathbb{X}} \min_{i \in \{1, \ldots, N \} }\psi_i (x) \]
\[ \qquad \leq  \max_{x \in \mathbb{X}} \psi_{\bar{\imath} } (x) - \psi_{\bar{\imath} } (x^* )  = \max_{x \in \mathbb{X}} \psi_{\bar{\imath}} (x) -  \min_{x \in \mathbb{X}} \psi_{\bar{\imath} } (x) \] 
\[ \qquad \qquad
\leq \max_{i \in \{1, \ldots, N \} }
\left \{  \max_{x \in \mathbb{X}} \psi_{i} (x) -  \min_{x \in \mathbb{X}} \psi_{i } (x) \right \}. \]
\endproof

Existence of fixed points of the shifted Bellman Equation can be used to establish useful upper and lower bounds on the rate of growth of the $T^k$ operator applied to any initial condition $\psi \in \mathcal{C}(\mathbb{X})$. This is stated in the following lemma.
\begin{lem}
\label{boundtop2}
Assume that there exists a continuous solution $\bar{\psi}$ to the shifted Bellman Equation, viz. $T \bar{\psi} = \bar{ \psi} + c$, for some $c \in \mathbb{R}$.
Then, for any positive integer $k$, and any function $\psi \in \mathcal{C}(x)$, the following holds:
\[ T^k \psi  (x) \leq \bar{ \psi} (x) + k c + \max_{x \in \mathbb{X}} [ \psi(x) - \bar{ \psi } (x) ] \]
\[ T^k \psi (x) \geq \bar{ \psi} (x) + k c + \min_{x \in \mathcal{X}} [ \psi(x) - \bar{ \psi} (x) ]. \]
\end{lem}
\emph{Proof.} To see the first inequality, notice:
\[ \psi (x) \leq \bar{ \psi} (x) + \max_{x \in \mathbb{X}} [ \psi (x) - \bar{ \psi} (x) ]. \]
Hence, exploiting monotonicity of the $\min$ operator we get:
\[ T^k \psi (x) = \min_{u( \cdot ), x(\cdot): x(0) = x } 
\sum_{t = 0}^{k-1} \ell (x(k),u(k)) + \psi ( x(k) ) \]
\[ \quad \leq \min_{u( \cdot ), x(\cdot): x(0) = x } 
\sum_{t = 0}^{k-1} \ell (x(k),u(k)) + \bar{\psi} ( x(k) ) + \max_{x \in \mathbb{X}} [ \psi(x) - \bar{ \psi } (x) ]  \]
\[ \qquad = \bar{ \psi} (x) + k c + \max_{x \in \mathbb{X}} [ \psi(x) - \bar{ \psi } (x) ]. \]
The second inequality can be proved along similar lines. \endproof

A direct consequence of Lemma \ref{boundtop2} is that:
\[ \max_{x \in \mathbb{X}} T^k \psi (x) - \min_{x \in \mathbb{X}} T^k \psi (x) \leq  \left [ \max_{x \in \mathbb{X}}  \bar{ \psi} (x) - \min_{x \in \mathbb{X}} \bar{ \psi} (x) \right ] \] \[ \qquad \qquad + \left [    \max_{x \in \mathbb{X}} [ \psi(x) - \bar{ \psi } (x) ]- \min_{x \in \mathbb{X}} [ \psi(x) - \bar{ \psi } (x) ]        \right ]. \]
Moreover, we can state the following corollary:
\begin{cor}
Assume there exists a continuous solution $\bar{\psi}$ to the shifted Bellman Equation, viz. $T \bar{\psi} = \bar{ \psi} + c$, for some $c \in \mathbb{R}$. Then, for any $\psi \in \mathcal{C}( \mathbb{X} )$ the following holds:
\[ \lim_{k \rightarrow + \infty} \frac{ T^k \psi (x) }{k} = c. \]
\end{cor}
\emph{Proof.} The result follows dividing by $k$ both sides of the inequalities in Lemma \ref{boundtop2}, and taking the limit as $k \rightarrow + \infty$. \endproof

Notice that, by construction, if the
sequence $\hat{T}^k \psi$ is bounded it converges to an upper semi-continuous function.
Analogously, if $\check{T}^k \psi$ is bounded it converges to a lower semi-continuous function.
If a continuous fixed point of the shifted Bellman Equation exists, both iterations might be suitable for determining such function, however, if no continuous fixed point exists, then it is not a priori clear which operator might be most suitable for the analysis. In fact, fixed points of the shifted Bellman Equation might be both upper or lower semi-continuous (or neither), despite the operator $T$
being in principle defined only on lower semi-continuous functions.

The next lemma shows that iterates of the $\hat{T}$ operator have a bounded excursion between their maximum and minimum value, provided a continuous fixed-point of the Bellman Equation exists. 

\begin{lem}
\label{boundednessT}
Assume that there exists a continuous solution to the shifted Bellman Equation, viz. $T \bar{\psi} = \bar{ \psi} + c$, for some $c \in \mathbb{R}$.
Then, the solution $\hat{T}^k \psi$ fulfills the bound:
\[ \max_{x \in \mathbb{X}} \hat{T}^k \psi(x) - \min_{x \in \mathbb{X}} \hat{T}^k \psi (x) \leq  \left [ \max_{x \in \mathbb{X}} \bar{ \psi} (x) - \min_{x \in \mathbb{X}} \bar{ \psi} (x) \right ]  \] 
\[ \qquad \qquad \qquad + \left [    \max_{x \in \mathbb{X}} [ \psi(x) - \bar{ \psi } (x) ]- \min_{x \in \mathbb{X}} [ \psi(x) - \bar{ \psi } (x) ]        \right ]. \]
\end{lem}
\emph{Proof.} To see this, notice that, by the $\min$-commutativity property, a simple induction argument shows, $\hat{T}^k \psi( x) = \min_{h \in \{0, \ldots, k \} } T^h \psi(x) + c_h$, for suitable values of $c_h \in \mathbb{R}$ and $c_0=0$.
By Lemma \ref{boundtopbis} we have:
\[  \max_{x \in \mathbb{X}} \hat{T}^k \psi(x) - \min_{x \in \mathbb{X}} \hat{T}^k \psi (x) \leq \max_{i \in \{1, \ldots, k \} } \left \{ \max_{x \in \mathbb{X}} [ T^i \psi(x) + c_i ] - \min_{x \in \mathbb{X}} [T^i \psi (x) +c_i] \right \}. \]
Canceling out the constant terms and exploiting Lemma \ref{boundtop2} yields:
\[  \max_x \hat{T}^k \psi(x) - \min_x \hat{T}^k \psi (x) \leq \max_{i \in \{1, \ldots, k \} } \left \{ \max_{x \in \mathbb{X}} T^i \psi(x)  - \min_{x \in \mathbb{X}} T^i \psi (x)  \right \} \]
\[ \qquad \leq  \left [ \max_{x \in \mathbb{X}} \bar{ \psi} (x) - \min_{x \in \mathbb{X}} \bar{ \psi} (x) \right ] + \left [    \max_{x \in \mathbb{X}} [ \psi(x) - \bar{ \psi } (x) ]- \min_{x \in \mathbb{X}} [ \psi(x) - \bar{ \psi } (x) ]        \right ]. \]
This last inequality completes the proof of the lemma.\endproof

Our subsequent analysis will rely on a combination of monotonicity and Lyapunov-based arguments. To this end it is useful to show that $\hat{T}$ and $\check{T}$ operators yield non-increasing iterations according to suitable Lyapunov functionals.
Exploiting Lemma \ref{boundtopbis} yields the following:
\begin{lem}
\label{Lyapunovoperator2}
Assume that there exists a continuous solution to the shifted Bellman Equation, viz. $T \bar{\psi} = \bar{ \psi} + c$, for some $c \in \mathbb{R}$.
Define the Lyapunov functional:
\begin{equation}
V( \psi) := \max_{x \in \mathbb{X}} [\psi (x) - \bar{ \psi} (x)] - \min_{x \in \mathbb{X}} [\psi (x) - \bar{ \psi} (x)].
\end{equation}
Then, for any continuous $\psi$ the following holds:
\begin{equation}
V ( \hat{T} \psi ) \leq V (\psi).
\end{equation}
\end{lem}
\emph{Proof.} Let $\psi \in \mathcal{C}( \mathbb{X})$ be arbitrary.
The inequality can be derived as follows:
\[ V ( \hat{T} \psi ) = \max_{x \in \mathbb{X}} \Big [ \hat{T} \psi (x) - \bar{ \psi } (x) \Big ] -  \min_{x \in \mathbb{X}} \Big [ \hat{T} \psi (x) - \bar{ \psi } (x) \Big ] \]
\[ \qquad = \max_{x \in \mathbb{X} } \Big [ \min \{ \psi(x), T \psi (x) + c (\psi, T \psi) \} - \bar{ \psi} (x) \Big ] \]
\[ \qquad \qquad -   \min_{x \in \mathbb{X} } \Big [  \min \{ \psi(x), T \psi (x) + c (\psi, T \psi) \} - \bar{ \psi} (x) \Big ] \]
\[ \qquad =  \max_{x \in \mathbb{X} } \min \{ \psi(x) - \bar{ \psi} (x) , T \psi (x) - \bar{ \psi} (x) + c (\psi, T \psi) \}  \]
\[ \qquad \qquad -   \min_{x \in \mathbb{X} } \min \{ \psi(x) - \bar{ \psi} (x), T \psi (x) - \bar{ \psi} (x)+ c (\psi, T \psi) \}  \]
\[ \qquad \leq \max \left \{   \max_{x \in \mathbb{X}} \big [ \psi (x) - \bar{ \psi} (x) \big ] - \min_{x \in \mathbb{X}} \big [ \psi (x) - \bar{ \psi} (x) \big], \right . \]
\[ \qquad \qquad \qquad \left . \max_{x \in \mathbb{X} } \Big [ T \psi (x) - \bar{ \psi} (x) + c (\psi, T \psi) \Big ] - \min_{x \in \mathbb{X}}  \Big [ T \psi (x) - \bar{ \psi} (x) + c (\psi, T \psi) \Big ] \right \} \]
\[ \qquad = \max \left \{ V( \psi) , \max_{x \in \mathbb{X} } [ T \psi (x) - T \bar{ \psi } (x)] - \min_{x \in \mathbb{X}} [ T \psi(x) - T \bar{ \psi} (x) ] \right \} \leq V( \psi ), \]
where the first inequality follows by Lemma \ref{boundtopbis}, and the second follows because $d (T \psi,T \bar{ \psi} ) \leq d ( \psi, \bar{ \psi} )$. \endproof

\begin{lem}
\label{Lyapunovoperator}
Assume that there exists a continuous solution to the shifted Bellman Equation, viz. $T \bar{\psi} = \bar{ \psi} + c$, for some $c \in \mathbb{R}$.
Define the Lyapunov functional:
\begin{equation}
V( \psi) := \max_{x \in \mathbb{X}} [\psi (x) - \bar{ \psi} (x)] - \min_{x \in \mathbb{X}} [\psi (x) - \bar{ \psi} (x)].
\end{equation}
Then, for any continuous $\psi$ the following holds:
\begin{equation}
V ( \check{T} \psi ) \leq V (\psi).
\end{equation}
\end{lem}
\emph{Proof.} The inequality can be derived as follows:
\[ V ( \check{T} \psi ) = \max_{x \in \mathbb{X}} \Big [ \check{T} \psi (x) - \bar{ \psi } (x) \Big ] -  \min_{x \in \mathbb{X}} \Big [ \check{T} \psi (x) - \bar{ \psi } (x) \Big ] \]
\[ \qquad = \max_{x \in \mathbb{X} } \max \{ \psi(x), T \psi (x) + c (\psi, T \psi) \} - \bar{ \psi} (x) \]
\[ \qquad \qquad -   \min_{x \in \mathbb{X} } \max \{ \psi(x), T \psi (x) + c (\psi, T \psi) \} - \bar{ \psi} (x) \]
\[ \qquad =  \max_{x \in \mathbb{X} } \max \{ \psi(x) - \bar{ \psi} (x) , T \psi (x) - \bar{ \psi} (x) + c (\psi, T \psi) \}  \]
\[ \qquad \qquad -   \min_{x \in \mathbb{X} } \max \{ \psi(x) - \bar{ \psi} (x), T \psi (x) - \bar{ \psi} (x)+ c (\psi, T \psi) \}  \]
\[ \qquad \leq \max \left \{   \max_{x \in \mathbb{X}} \big [ \psi (x) - \bar{ \psi} (x) \big ] - \min_{x \in \mathbb{X}} \big [ \psi (x) - \bar{ \psi} (x) \big], \right . \]
\[ \qquad \qquad \qquad \left . \max_{x \in \mathbb{X} } \Big [ T \psi (x) - \bar{ \psi} (x) + c (\psi, T \psi) \Big ] - \min_{x \in \mathbb{X}}  \Big [ T \psi (x) - \bar{ \psi} (x) + c (\psi, T \psi) \Big ] \right \} \]
\[ \qquad = \max \left \{ V( \psi) , \max_{x \in \mathbb{X} } [ T \psi (x) - T \bar{ \psi } (x)] - \min_{x \in \mathbb{X}} [T \psi(x) - T \bar{ \psi} (x)] \right \} \leq V( \psi ), \]
where the first inequality follows by Lemma \ref{boundtop1}, and the second follows because $d (T \psi,T \bar{ \psi} ) \leq d ( \psi, \bar{ \psi} )$. \endproof

An alternative Lyapunov functional for the operator $\hat{T}$ can be stated as follows:
\begin{equation}
\label{previouslyap}
    W( \psi ) := d ( \psi, T \psi ).
\end{equation}
The following lemma proves that this is non increasing along iterations of $\hat{T}(\cdot)$.
\begin{lem}
\label{newlyap}
Consider the function $W(\psi)$ defined in (\ref{previouslyap}). For any real valued continuous function $\psi:\mathbb{X} \rightarrow \mathbb{R}$ the following holds:
\[ W ( \hat{T} \psi ) \leq W (\psi). \]
\end{lem}
\emph{Proof.} To prove the lemma consider the following inequalities:
\begin{equation}
    \label{exploited1}
 T \psi(x) + c( \psi,T \psi)  - d( \psi, T \psi ) \leq \hat{T} \psi (x) \leq T \psi (x) + c( \psi, T \psi), \quad \forall \, x \in \mathbb{X}. 
\end{equation}
In addition, by definition of $\hat{T} \psi$, we see that:
\[   \psi(x) - d (\psi,T \psi) \leq  \hat{T} \psi(x) \leq \psi(x). \]
By monotonicity and translation invariance, applying the $T$ operator to all sides of the former inequality yields:
\begin{equation}
    \label{exploited2}
    T \psi(x) - d(\psi,T \psi)      \leq                 T \hat{T} \psi (x) \leq T \psi(x) , \qquad \forall x \in \mathbb{X}. 
\end{equation}    
We are now ready to estimate $W (\hat{T} \psi)$ by combining inequalities
(\ref{exploited1}) and (\ref{exploited2}):
\[ W (\hat{T} \psi ) = \frac{1}{2} \left ( \max_{x \in \mathbb{X}} 
[ \hat{T} \psi (x) - T \hat{T} (\psi) (x) ] -
\min_{x \in \mathbb{X}} [ \hat{T} \psi(x) - T \hat{T} ( \psi)(x) ] \right )\]
\[\leq \frac{1}{2} \left ( \max_{x \in \mathbb{X}} 
[ T \psi (x) + c ( \psi,T\psi)  - T \psi(x) + d( \psi, T \psi)  ] \right . \qquad \qquad \qquad \]
\[ \qquad \qquad \left . -
\min_{x \in \mathbb{X}} [ T \psi(x) + c( \psi,T \psi) - d(\psi,T\psi) - T  \psi(x) ] \right ) = d (\psi,T\psi) = W(\psi).    \]
\endproof

\begin{rem}
\label{usedinproofbelow}
The same argument used to prove Lemma \ref{newlyap} can be used to prove the following decoupled inequalities:
\begin{equation}
\label{goingdown}
\max_{x \in \mathbb{X} } [ \hat{T} \psi (x) - T \hat{T} \psi (x) ]
\leq \max_{x \in \mathbb{X}} [ \psi (x) - T \psi (x) ]
\end{equation}
and:
\begin{equation}
\label{goingup}
\min_{x \in \mathbb{X} } [\hat{T} \psi (x) - T \hat{T} \psi (x)]
\geq \min_{x \in \mathbb{X}} [\psi (x) - T \psi (x)].
\end{equation}
\end{rem}

Our analysis indicates that regardless of whether the sequence of functions $\hat{T}^k \psi(x)$ converges, the real-valued sequence of shifts applied, $c( \hat{T}^k \psi, T \hat{T}^k \psi )$ is always bounded and convergent. 

\begin{lem}
\label{cproperties}
The sequence $c ( \hat{T}^k \psi, T \hat{T}^k \psi )$ is bounded and convergent, viz. there exists $\hat{c}_{\infty} \in \mathbb{R}$ such that:
\[ \lim_{k \rightarrow + \infty} c(\hat{T}^k \psi, T \hat{T}^k \psi ) = \hat{c}_{\infty}. \]
\end{lem}
\emph{Proof.} By induction, and Remark \ref{usedinproofbelow} we have that
the real-valued sequence:
$\max_{x \in \mathbb{X}} [\hat{T}^k \psi (x) - T \hat{T}^k \psi(x)] $ is monotonically non-increasing, (and bounded from below
by  $\min_{x \in \mathbb{X} } [\hat{T}^k \psi (x) - T \hat{T}^k \psi (x) ]$ ).
Similarly, $\min_{x \in \mathbb{X} } [\hat{T}^k \psi (x) - T \hat{T}^k \psi (x)]$ is monotonically non-decreasing (and bounded from above by $\max_{x \in \mathbb{X}} [\hat{T}^k \psi (x) - T \hat{T}^k \psi(x)] $).
Hence, both sequences admit a limit:
\[ M = \lim_{k \rightarrow + \infty} \max_{x \in \mathbb{X}} [\hat{T}^k \psi (x) - T \hat{T}^k \psi(x)] \]
\[ m = \lim_{k \rightarrow + \infty} \min_{x \in \mathbb{X}} [\hat{T}^k \psi (x) - T \hat{T}^k \psi(x)] \]
By definition of $c(\hat{T}^k \psi,T \hat{T}^k \psi )$ we see that:
\[ \lim_{k \rightarrow + \infty} c(\hat{T}^k \psi, T \hat{T}^k \psi ) = (M+m)/2 := \hat{c}_{\infty}, \]
which completes the proof of the lemma. \endproof

We turn next to establishing similar inequalities for the $\check{T}$ operator.
\begin{lem}
\label{newlyap2}
Consider the function $W$ defined in (\ref{previouslyap}). For any real valued continuous function $\psi:\mathbb{X} \rightarrow \mathbb{R}$ the following holds:
\[ W ( \check{T} \psi ) \leq W (\psi). \]
\end{lem}
\emph{Proof.} To see the inequality consider that we have:
\begin{equation}
    \label{exploited1bis}
 T \psi(x) + c( \psi,T \psi) \leq \check{T} \psi (x) \leq T \psi (x) + c( \psi, T \psi) + d ( \psi, T \psi), \quad \forall \, x \in \mathbb{X}. 
\end{equation}
In addition, by definition of $\check{T} \psi$, we see that:
\[   \psi(x) \leq  \check{T} \psi(x) \leq \psi(x) + d ( \psi, T \psi) \]
By monotonicity and translation invariance, applying the $T$ operator to all sides of the former inequality yields:
\begin{equation}
    \label{exploited2bis}
    T \psi(x)      \leq                 T \check{T} \psi (x) \leq T \psi(x) + d (\psi,T \psi ), \qquad \forall x \in \mathbb{X}. 
\end{equation}    
We are now ready to bound from above $W (\hat{T} \psi)$ by combining inequalities
(\ref{exploited1bis}) and (\ref{exploited2bis}):
\[ W (\check{T} \psi ) = \frac{1}{2} \left ( \max_{x \in \mathbb{X}} 
[ \check{T} \psi (x) - T \check{T} (\psi) (x) ] -
\min_{x \in \mathbb{X}} [ \check{T} \psi(x) - T \check{T} ( \psi)(x) ] \right )\]
\[\leq \frac{1}{2} \left ( \max_{x \in \mathbb{X}} 
[ T \psi (x) + c ( \psi,T\psi)  + d ( \psi,T \psi) - T \psi(x)   ] \right . \qquad \qquad \qquad \]
\[ \qquad \qquad \left . -
\min_{x \in \mathbb{X}} [ T \psi(x) + c( \psi,T \psi) - T  \psi(x) - d ( \psi, T \psi) ] \right ) = d (\psi,T\psi) = W(\psi).    \]
\endproof

\begin{rem}
\label{usedinproofbelow2}
The same argument used to prove Lemma \ref{newlyap2} can also be used to prove the following decoupled inequalities:
\begin{equation}
\label{goingdown2}
\max_{x \in \mathbb{X} } [\check{T} \psi (x) - T \check{T} \psi (x)]
\leq \max_{x \in \mathbb{X}} [\psi (x) - T \psi (x)]
\end{equation}
and:
\begin{equation}
\label{goingup2}
\min_{x \in \mathbb{X} } [\check{T} \psi (x) - T \check{T} \psi (x)]
\geq \min_{x \in \mathbb{X}} [\psi (x) - T \psi (x)].
\end{equation}
\end{rem}
A similar proof as in Lemma \ref{cproperties} allows to conclude the following result:
\begin{lem}
\label{cproperties2}
The sequence $c ( \check{T}^k \psi, T \check{T}^k \psi )$ is bounded and convergent, viz. there exists $\check{c}_{\infty} \in \mathbb{R}$ such that:
\[ \lim_{k \rightarrow + \infty} c(\check{T}^k \psi, T \check{T}^k \psi ) = \check{c}_{\infty}. \]
\end{lem}

It seems important to relate the value of $\hat{c}_{\infty}$ and $\check{c}_{\infty}$ with the optimal average infinite horizon cost, viz. $V^{avg}$.
The following result shows that $-\hat{c}_{\infty}$ is always an upper-bound to the optimal average cost.
\begin{lem}
\label{hatclemma1}
Assume that a fixed-point of the shifted Bellman Equation exists, viz. $T \bar{\psi} = \bar{\psi} + c$ for some $c \in \mathbb{R}$. Then, for any $\psi \in \mathcal{C}(\mathbb{X})$ it holds:
\begin{equation}
    c + \hat{c}_{\infty} \leq 0.
\end{equation} 
\end{lem}
\emph{Proof.}
We argue by contradiction. Assume that $c+ \hat{c}_{\infty} > 0$. Then, there exists $\varepsilon>0$ and $Q \in \mathbb{N}$ such that:
\begin{equation}
    \label{quotednext}
    c + c (\hat{T}^k \psi, T \hat{T}^k \psi ) \geq \varepsilon > 0 
    \end{equation}
for all $k \geq Q$. Moreover, there exists $N \in \mathbb{N}$ such that for any $x \in \mathbb{X}$
\begin{equation}
\label{quotednext2}
N \varepsilon \geq - \min_{x \in \mathbb{X}} [ \hat{T}^Q \psi(x) - \bar{\psi} (x)] - \bar{\psi} (x) + \psi(x). 
\end{equation}
We claim that, under such assumptions, $\hat{T}^k \psi(x)$ converges to a fixed-point within a finite number of iterations.
In fact, for any $m \geq N$ we see that:
\[ \hat{T}^{Q + m} \psi (x)  =  \hat{T}^m ( \hat{T}^Q \psi (x) ) \qquad \qquad \qquad \qquad \qquad  \]
\[ = \min_{\tau \in \{0, \ldots, m \}}  \left \{ T^{\tau} \hat{T}^Q \psi(x) + \min_{ S \subseteq \{0,1,\ldots,m-1 \}, |S|= \tau }  \left [ \sum_{s \in S} c ( \hat{T}^{s+Q} \psi, T \hat{T}^{s+Q} \psi ) \right ] \right \} \]
\[=  \min_{\tau \in \{0, \ldots, N\}}  \left \{ T^{\tau} \hat{T}^Q \psi(x) + \min_{ S \subseteq \{0,1,\ldots,N-1 \}, |S|= \tau}   \left [ \sum_{s \in S} c ( \hat{T}^{s+Q} \psi, T \hat{T}^{s+Q} \psi ) \right ] \right \},
\]
where the last equality holds because for $\tau \geq N$ application of Lemma \ref{boundtop2} and inequality (\ref{quotednext}) yields:
\[ T^{\tau} \hat{T}^Q \psi (x) + \min_{S \subseteq \{0,1,\ldots,m-1 \}, |S|= \tau}  \left [ \sum_{s \in S} c ( \hat{T}^{s+Q} \psi, T \hat{T}^{s+Q} \psi ) \right ] \]
\[ \qquad \geq \bar{ \psi} (x) + \tau c  + \min_{x \in \mathbb{X}} [ \hat{T}^Q \psi (x) - \bar{\psi} (x) ] + \tau ( \varepsilon - c ) \] 
\[ \qquad = \bar{ \psi} (x) + \tau \varepsilon  + \min_{x \in \mathbb{X}} [ \hat{T}^Q \psi (x) - \bar{\psi} (x) ] 
\geq \psi(x). \]
Hence $\hat{T}^{Q+N} \psi(x) = \lim_{k \rightarrow + \infty} \hat{T}^k \psi (x)$ where convergence is in a finite number of steps (uniform over $x \in \mathbb{X}$). Moreover, 
\[ \hat{T} ( \hat{T}^{Q+N} \psi(x) ) = \hat{T}^{Q + N +1 } \psi(x) = \hat{T}^{Q+N} \psi(x). \]
Therefore, $\hat{T}^{Q+N} \psi(x)$ is a (continuous) fixed point of the $\hat{T}$ operator, and by
virtue of Proposition \ref{fixedpointequalshifted} it is a solution of the shifted Bellman Equation for some $c = - c ( \hat{T}^{Q+N} \psi, T \hat{T}^{Q+N} \psi)$. This implies $c + \hat{c}_{\infty} =0$, which is a contradiction. 
\endproof

Whenever the sequence $\hat{T}^k \psi$ is pointwise convergent, one can show that also the converse inequality holds, and therefore $-\hat{c}_{\infty}$ equals the optimal average performance.
The next lemma is instrumental in deriving such result.

\begin{lem}
\label{monotonesequenceminimum2} Let $J_k (u): \mathbb{U} \rightarrow \mathbb{R}$ be a monotonically non-increasing sequence of continuous functions,  converging pointwise to $\hat{J}(u)$, and let $\mathbb{U}$ be a compact set. Then the following holds:
\[ \lim_{k \rightarrow + \infty} \min_{u \in \mathbb{U} } J_k (u) = \inf_{u \in \mathbb{U}} \hat{J} (u). \]
\end{lem}
\emph{Proof.} Remark that the function $\hat{J}$ is upper semi-continuous, but not necessarily lower semi-continuous. Hence its minimum might, a priori, not be well-defined.
By monotonicity of the minimum operator:
\[  \min_{u \in \mathbb{U}} J_{k} (u) \geq \min_{u \in \mathbb{U}} J_{k+1} (u), \]
for all $k \in \mathbb{N}$. Hence, the limit
$\lim_{k \rightarrow + \infty} \min_{u \in \mathbb{U} } J_k (u)$, exists.
Moreover, by monotonicity of the $\inf$ operator we see:
\[  \min_{u \in \mathbb{U}} J_{k} (u) =  \inf_{u \in \mathbb{U}} J_{k} (u) \geq \inf_{u \in \mathbb{U}} \hat{J} (u), \]
which holds for all $k \in \mathbb{N}$. Letting $k$ go to infinity in the previous inequality shows:
\[  \lim_{k \rightarrow + \infty} \min_{u \in \mathbb{U} } J_k (u) \geq \inf_{u \in \mathbb{U}} \hat{J} (u). \]
We need to show the converse inequality. To this end, denote by $u_n$ any sequence in $\mathbb{U}$ such that $\hat{J} (u_n) - 2^{-n} \leq \inf_{u \in \mathbb{U}} \hat{J}(u)$.
Clearly, for any $n$, there exists $k_n > n$ such that $J_{k_n} (u_n) \leq \hat{J} (u_n) + 2^{-n}$.
Overall we see:
\[ \inf_{u \in \mathbb{U}} \hat{J}(u) \geq \hat{J} (u_n) - 2^{-n} \geq J_{k_n} (u_n) - 2 \cdot 2^{-n}
\geq \min_{u \in \mathbb{U} } J_{k_n} (u) - 2 \cdot 2^{-n}. \]
Letting $n$ go to infinity in the previous inequality yields:
\[ \inf_{u \in \mathbb{U}} \hat{J}(u) \geq \lim_{k \rightarrow + \infty} \min_{u \in \mathbb{U}} J_k (u) .\]
This completes the proof of the lemma.\endproof

It is sometimes useful to consider the extension of operator $T$ to functions $\psi$ bounded from below (and non-necessarily continuous).
To this end, if $\psi: \mathbb{X} \rightarrow \mathbb{R}$ is bounded from below, we denote by $T \psi$ the following:
\[ T \psi (x) = \inf_{u \in \mathbb{U}(x) } \ell(x,u) + \psi ( f(x,u) ). \]

\begin{lem}
\label{subfixedpoint} 
Assume that the function $\hat{T}^k \psi (x)$ converges pointwise to $\hat{ \psi} (x)$, bounded from below.
Then the following holds:
\[ \hat{ \psi} (x) \leq T \hat{ \psi} + \hat{c}_{\infty}. \] 
\end{lem}
\emph{Proof.} To prove the lemma notice that:
\[ \hat{ \psi} (x) = \lim_{k \rightarrow + \infty} \hat{T}^{k+1} \psi (x) \]
\[ \qquad = \lim_{k \rightarrow + \infty} \min \{ \hat{T}^k \psi (x), T \hat{T}^k \psi(x) + c ( \hat{T}^k \psi, T\hat{T}^k \psi) \} \]
\[ \qquad \leq 
\lim_{k \rightarrow + \infty} T \hat{T}^k \psi(x) + c ( \hat{T}^k \psi, T\hat{T}^k \psi)  = T \hat{ \psi} (x) + \hat{c}_{\infty}. \]
where the last equality follows by applying Lemma \ref{monotonesequenceminimum2} to the sequence of $x$-parameterized functions
$J_k(x,u) := \ell(x,u) + \hat{T}^k \psi (f(x,u))$.\endproof

We are now ready to prove the converse inequality.

\begin{lem} 
\label{hatclemma2}
Assume the sequence $\hat{T}^k \psi (x)$ to be pointwise convergent to some bounded function $\hat{ \psi} (x)$. If a fixed point $\bar{\psi}$ of the shifted Bellman Equation exists, viz. $T \bar{\psi} = \bar{\psi} + c$ for some $c \in \mathbb{R}$ and some $\bar{\psi}: \mathbb{X} \rightarrow \mathbb{R}$,
the following holds:
\[ 0 \leq c + \hat{c}_{\infty}. \] 
\end{lem}
\emph{Proof.} By Lemma \ref{subfixedpoint}, we see that $\hat{ \psi} (x) \leq T \hat{ \psi} + \hat{c}_{\infty}.$
Monotonicity of $T$ together with shift-invariance yields, by induction for $k \in \mathbb{N}$:
\[ \hat{\psi}(x) \leq T^k \hat{\psi} (x) + k \hat{c}_{\infty}. \]
In particular then, for any continuous $\psi \geq \hat{\psi}$:
\[ \hat{\psi}(x) \leq T^k \psi(x) + k \hat{c}_{\infty}. \]
Dividing both sides of the previous inequality by $k$ and letting $k$ tend to infinity yields:
\[ 0 \leq c + \hat{c}_{\infty}. \]
\endproof

A similar analysis can be carried out with respect to the iteration $\check{T}^k \psi$ and the corresponding limiting value of the shift $\check{c}_{\infty}$. As a matter of fact,
not all results extend along the same lines, due to the lack of formula (\ref{minformula}).
We first state the analogue of Lemma \ref{monotonesequenceminimum2}.

\begin{lem}
\label{monotonesequenceminimum} Let $J_k (u): \mathbb{U} \rightarrow \mathbb{R}$ be a monotonically non-decreasing sequence of (lower semi-)continuous functions,  converging pointwise to $\check{J}(u)$, and let $\mathbb{U}$ be a compact set. Then the following holds:
\[ \lim_{k \rightarrow + \infty} \min_{u \in \mathbb{U} } J_k (u) = \min_{u \in \mathbb{U}} \check{J} (u). \]
\end{lem}
\emph{Proof.} Note that the function $\check{J}$ is lower semicontinuous, hence its minimum is well defined.
By monotonicity of the minimum operator:
\[  \min_{u \in \mathbb{U}} J_{k} (u) \leq \min_{u \in \mathbb{U}} J_{k+1} (u), \]
for all $k \in \mathbb{N}$. Hence, the limit
$\lim_{k \rightarrow + \infty} \min_{u \in \mathbb{U} } J_k (u)$, exists.
Moreover, again by monotonicity of the $\min$ operator we see:
\[  \min_{u \in \mathbb{U}} J_{k} (u) \leq \min_{u \in \mathbb{U}} \check{J} (u), \]
which holds for all $k \in \mathbb{N}$. Letting $k$ go to infinity in the previous inequality shows:
\[  \lim_{k \rightarrow + \infty} \min_{u \in \mathbb{U} } J_k (u) \leq \min_{u \in \mathbb{U}} \check{J} (u). \]
We need to show the converse inequality. To this end, denote by $u_n$ any element of $\mathbb{U}$ such that $J_n (u_n) = \min_{u \in \mathbb{U}} J_n(u)$.
For any $k \in \mathbb{N}$ and any $n \geq k$ we see that $J_n (u_n) \geq J_k (u_n)$. In particular, then:
\[ \limsup_{n \rightarrow + \infty} J_n (u_n) \geq \limsup_{n \rightarrow + \infty} J_k (u_n) \geq J_k (u^*) \]
for some limit point $u^* \in \mathbb{U}$ of the sequence $u_n$.
Hence:
\[  \lim_{k \rightarrow + \infty} \min_{u \in \mathbb{U} } J_k (u) =  \limsup_{n \rightarrow + \infty} J_n (u_n) \geq J_k (u^*) \]
for all $k \in \mathbb{N}$, and letting $k$ go to infinity in the right hand side of the previous inequality yields:
\[  \lim_{k \rightarrow + \infty} \min_{u \in \mathbb{U} } J_k (u) \geq \check{J} (u^*) \geq \min_{u \in \mathbb{U}} \check{J}(u). \]
This completes the proof of the lemma. \endproof

\begin{cor}
\label{mincor}
Assume that $\check{T}^k \psi(x)$ converges point-wise to a lower semi-continuous limit $\check{\psi} (x)$, for all $x \in \mathbb{X}$. Applying Lemma \ref{monotonesequenceminimum} to the $x$-parameterised sequence of cost functions:
\[ J_k (x,u) := \ell(x,u) + \check{T}^k \psi ( f(x,u) ) \]
admitting the limit:
\[ \check{J} (x,u) := \ell(x,u) + \check{\psi} (f(x,u)), \]
with $u \in \mathbb{U}(x)$, yields the following point-wise convergence result:
\[ \lim_{k \rightarrow + \infty} T \check{T}^k \psi (x) = T \check{\psi} (x) \]
\end{cor}

\begin{lem}
\label{supfixedpoint}
Assume that $\check{T}^k \psi(x)$ converges pointwise to a lower semi-continuous limit $\check{\psi} (x)$, for all $x \in \mathbb{X}$. Then, $\check{\psi} (x)$ fulfills:
\[ \check{ \psi} (x) \geq T \check{\psi} (x) + \check{c}_{\infty}. \]
If in addition the limit $\check{\psi}(x)$ is continuous, then it is a fixed point of a shifted Bellman Equation. 
\end{lem}
\emph{Proof.} For all $k \in \mathbb{N}$ we see:
\[ \check{\psi} (x) \geq \check{T}^{k+1} \psi (x) \geq
T \check{T}^k \psi (x) + c ( \check{T}^k \psi, T \check{T}^k \psi ). \]
Hence, by Corollary \ref{mincor}, letting $k \rightarrow + \infty$ in the right-hand side of the latter inequality yields:
\[ \check{ \psi} (x) \geq T \check{\psi} (x) + \check{c}_{\infty}. \]
In addition, if $\check{\psi} \in \mathcal{C}( \mathbb{X})$ then, by Dini's theorem, convergence is uniform and $\check{\psi}$ is a fixed point of $\check{T}$ by continuity of the $\check{T}$ operator in the topology of uniform convergence. \endproof

We are now ready to state the analogue of Lemma \ref{hatclemma2}.
\begin{lem} 
\label{checkclemma2}
Assume that the sequence $\check{T}^k \psi (x)$ be pointwise convergent to some bounded function $\check{ \psi} (x)$.
If a fixed point $\bar{\psi}$ of the shifted Bellman Equation exists, viz. $T \bar{\psi} = \bar{\psi} + c$ for some $c \in \mathbb{R}$, the following holds:
\[ 0 \geq c + \check{c}_{\infty}. \] 
\end{lem}
\emph{Proof.} By Lemma \ref{supfixedpoint}, we see that $\check{ \psi} (x) \geq T \check{ \psi} + \check{c}_{\infty}.$
Monotonicity of $T$ together with shift-invariance yields, by induction for $k \in \mathbb{N}$:
\[ \check{\psi}(x) \geq T^k \check{\psi} (x) + k \check{c}_{\infty}. \]
In particular then, for any continuous $\psi \leq \hat{\psi}$:
\[ \check{\psi}(x) \geq T^k \psi(x) + k \check{c}_{\infty}. \]
Dividing both sides of the previous inequality by $k$ and letting $k$ tend to infinity yields:
\[ 0 \geq c + \check{c}_{\infty}. \]
\endproof

A stronger claim can be achieved when the $\hat{T}^k \psi$ and $\check{T}^k \psi$ sequences admit a continous limit.
\begin{lem}
\label{fixcontinuous}
Let $\psi(x)$ be a continuous function, and assume that 
$\hat{T}^{k} \psi(x)$ ( or $\check{T}^k \psi (x) $ ) converges point-wise to a continuous limit $\hat{\psi} (x)$ ( $\check{\psi} (x)$, respectively ). Then,
$\hat{\psi}$ ($\check{\psi}$, respectively) is a fixed point of the shifted Bellman Equation. (a similar argument holds for $\check{T}$ ).
\end{lem}
\emph{Proof.}
By construction $\hat{T}^k \psi$ is monotone non-increasing with respect to $k$. Hence, by Dini's Theorem, convergence to $\hat{\psi}$
is uniform.
The result follows by \emph{continuity} of the $T ( . ) $
and $c( . , . )$ operators with respect to the topology of uniform convergence. \endproof

The convergence properties of $\hat{T}^k \psi$ and $\check{T}^k \psi$ sequences will be established through a combination of Lasalle-style and monotonicity-based arguments. The following lemmas are crucial to understand the implication of certain Lyapunov functionals being constant along iterations of the $\hat{T}$ and $\check{T}$ maps. 

\begin{lem}
\label{laslikelemma}
Let $\psi$ be a continuous function such that:
\[   \min_{x \in \mathbb{X}} \hat{T} \psi (x) -
T \hat{T} \psi (x) = \min_{x \in \mathbb{X}} \psi(x) - T \psi (x). \]
Then the following holds:
\begin{itemize}
    \item the sets achieving the minimum are nested:
    \[ \arg \min_{x \in \mathbb{X}} 
    \hat{T} \psi (x) -
T \hat{T} \psi (x) \subseteq \arg \min_{x \in \mathbb{X}} \psi(x) - T \psi (x)
    \] 
    \item the operator $\hat{T}$ does not alter the value of $\psi$ in the $\arg \min$ set:
    \[ \psi (x_m) = \hat{T} \psi (x_m) \; \qquad \forall \, x_m \in \arg \min_{x \in \mathbb{X}} 
    \hat{T} \psi (x) -
T \hat{T} \psi (x) \]
\item $T \hat{T} \psi(x_m) = T\psi (x_m)$ for all $ x_m \in \arg \min_{x \in \mathbb{X}} 
    \hat{T} \psi (x) -
T \hat{T} \psi (x)$.
\end{itemize}
\end{lem}
\emph{Proof.} To prove the lemma notice that inequality
(\ref{goingup}) holds, and can be derived from the following inequalities:
$T \hat{T} \psi (x) \leq T \psi(x)$ and
$\hat{T} \psi (x) \geq T \psi (x) + c(\psi,T\psi) - d ( \psi,T \psi)$. If $(\ref{goingup})$ is an equality, both previous inequalities need to be fulfilled non strictly for any
$x_m \in \arg \min_{x \in \mathbb{X}} \hat{T} \psi (x) - T \hat{T} \psi (x) $. Hence, it holds:
\[ \hat{T} \psi (x_m) = T \psi (x_m) + c ( \psi,T \psi) - d ( \psi, T \psi), \]
and
\[ T \hat{T} \psi (x_m) = T \psi (x_m).\]
Since $d ( \psi,T \psi) \geq 0$, the first equality proves that
$\hat{T} \psi (x_m) = \psi (x_m)$.
Moreover, by assumption:
\[ \min_{x \in \mathbb{X}} \psi (x) - T \psi(x)
 = \hat{T}\psi (x_m) - T \hat{T} \psi(x_m) = \psi(x_m) - T \psi(x_m). \]
This shows that $x_m \in \arg \min_{x \in \mathbb{X}} \psi(x) - T \psi (x).$ Since $x_m$ was arbitrary to start with, inclusion of the $\arg \min$ sets follows, which concludes the proof of the lemma.  \endproof

\begin{cor}
\label{convergeupper}
Assume that a continuous fixed point of $\bar{\psi}$ of the shifted Bellman Equation exists. If for some continuous $\psi$ and all $k \in \mathbb{N}$ it holds
\begin{equation}
    \label{hypcorref}
 \min_{x \in \mathbb{X}} \hat{T}^k \psi (x) -
T \hat{T}^k \psi (x) = \min_{x \in \mathbb{X}} \psi(x) - T \psi (x)  
\end{equation}
then, $\lim_{k \rightarrow + \infty} \hat{T}^k \psi$ exists and is an upper-semicontinuous function.  
\end{cor}
\emph{Proof}.
By virtue of Lemma \ref{laslikelemma}, if equation (\ref{hypcorref}) holds there exists $x_m \in \mathbb{X}$ such that $\hat{T}^k \psi (x_m) = \psi (x_m)$ for all $k \in \mathbb{N}$. In particular,
\[ \max_{x \in \mathbb{X}} \hat{T}^k \psi (x)
 \geq \hat{T}^k \psi (x_m) = \psi (x_m) \geq \min_{x \in \mathbb{X}} \psi(x), \]
 which in combination with the inequality proved in Lemma \ref{boundednessT} and existence of a fixed point of the shifted Bellman Equation
implies boundedness and pointwise convergence of the $\hat{T}$ iteration.
Moreover, as $\hat{T}^k \psi$ is non-increasing the limiting function is upper-semicontinuous. \\

A symmetric argument can be used to establish the following lemma.

\begin{lem}
\label{laslikelem2}
Let $\psi$ be a continuous function such that:
\begin{equation} 
\label{hypcor2}
\max_{x \in \mathbb{X}} \check{T} \psi (x) -
T \check{T} \psi (x) = \max_{x \in \mathbb{X}} \psi(x) - T \psi (x). \end{equation}
Then the following holds:
\begin{itemize}
    \item the sets achieving the maximum are nested:
    \[ \arg \max_{x \in \mathbb{X}} 
    \check{T} \psi (x) -
T \check{T} \psi (x) \subseteq \arg \max_{x \in \mathbb{X}} \psi(x) - T \psi (x)
    \] 
    \item the operator $\check{T}$ does not alter the value of $\psi$ in the $\arg \max$ set:
    \[ \psi (x_m) = \check{T} \psi (x_m) \; \qquad \forall \, x_m \in \arg \max_{x \in \mathbb{X}} 
    \check{T} \psi (x) -
T \hat{T} \psi (x) \]
\item $T \check{T} \psi(x_m) = T\psi (x_m)$ for all $ x_m \in \arg \max_{x \in \mathbb{X}} 
    \check{T} \psi (x) -
T \check{T} \psi (x)$.
\end{itemize}
\end{lem}
A version of Corollary \ref{convergeupper} can be proved for the $\check{T}$ operator.
\begin{cor}
\label{convergelower}
Assume that a continuous fixed point of $\bar{\psi}$ of the shifted Bellman Equation exists. If for some continuous $\psi$ and all $k \in \mathbb{N}$ it holds
\begin{equation}
    \label{hypcor}
 \max_{x \in \mathbb{X}} \check{T}^k \psi (x) -
T \check{T}^k \psi (x) = \max_{x \in \mathbb{X}} \psi(x) - T \psi (x)  
\end{equation}
then, $\lim_{k \rightarrow + \infty} \check{T}^k \psi$ exists and is a lower-semicontinuous function.  
\end{cor}
\emph{Proof}. By virtue of Lemma \ref{laslikelem2}, if equation (\ref{hypcor2}) holds there exists $x_m \in \mathbb{X}$ such that $\check{T}^k \psi (x_m) = \psi (x_m)$ for all $k \in \mathbb{N}$. In particular,
\begin{equation}
\label{boundingmin}
\min_{x \in \mathbb{X}} \check{T}^k \psi (x)
 \leq \check{T}^k \psi (x_m) = \psi (x_m) \leq \max_{x \in \mathbb{X}} \psi(x). 
 \end{equation}
We show next that the sequence $\check{T}^k \psi$ is bounded from above:
\[ \max_{x \in \mathbb{X}} \check{T}^k \psi (x) = \max_{x \in \mathbb{X}} [ \check{T}^k \psi (x) - \bar{ \psi} (x) + \bar{\psi} (x) ]\qquad \qquad \qquad \]
\[ \leq \max_{x \in \mathbb{X}} [ \check{T}^k \psi (x) - \bar{ \psi} (x) ] + \max_{x \in \mathbb{X}} \bar{\psi} (x) \]
\[ = d ( \check{T}^k \psi, \bar{\psi}  ) + \min_{x \in \mathbb{X}} [ \check{T}^k \psi (x) - \bar{ \psi} (x) ] + \max_{x \in \mathbb{X}} \bar{\psi} (x) \]
\[ \leq d( \psi, \bar{\psi} ) + \min_{x \in \mathbb{X}} \check{T}^k \psi (x)  - \min_{x \in \mathbb{X}} \bar{ \psi} (x)  + \max_{x \in \mathbb{X}} \bar{\psi} (x) \]
 \[ \leq  d( \psi, \bar{\psi} ) + \max_{x \in \mathbb{X}} \psi(x)  - \min_{x \in \mathbb{X}} \bar{ \psi} (x)  + \max_{x \in \mathbb{X}} \bar{\psi} (x), \]
 where the last inequality follows by (\ref{boundingmin}) and the former one by Lemma \ref{Lyapunovoperator}.
 Hence, pointwise convergence of the $\check{T}^k \psi$ sequence to a lower semi-continuous function follows by boundedness and monotonicity (viz. by $\check{T}^k \psi$ being non-decreasing in $k$). \endproof

\section{Appendix: Additional proofs}
\subsection{Proof of Proposition \ref{sameaverage}}

\label{sameaverageproof}
Let $k \in \mathbb{N}$ be arbitrary.
By induction it is possible to see that:
\[ T^k \psi_1 + k c_1  = \psi_1 .\] 
The claim is trivial for $k=1$. Assume this holds for $k$, we will show it is true for $k+1$:
\[ T^{k+1} \psi_1 = T ( T^k \psi_1 ) = T ( \psi_1 - k c_1 ) = T \psi_1 - k c_1 = \psi_1 - c_1 - k c_1 = \psi_1 - (k+1) c_1. \]
A similar argument applies to $\psi_2$. In particular then:
\[ \lim_{k \rightarrow + \infty} \frac{ T^k \psi_i (x_0) }{k} = - c_i, \qquad i = 1,2. \]
 Moreover, we know that:
\[  V_k^{\psi_1} (x_0) = T^k \psi_1 (x_0) , \] 
for all $k \in \mathbb{N}$ and all $x_0$.
Assume that $|\psi_1(x) - \psi_2 (x)| \leq M$ for all $x \in \mathbb{X}$, which is always fulfilled for sufficiently large $M$ due to boundedness of $\psi_1$ and $\psi_2$.
Then $V_k^{\psi_2} (x_0) \leq V_k^{\psi_1} (x_0) + M$ as it follows by remarking that the optimal solution relative to the terminal penalty $\psi_1$ can be used as a feasible solution to estimate the optimal cost of the problem with terminal 
cost $\psi_2$.
A symmetric argument also yields $V_{k}^{\psi_1} (x_0) \leq V_k^{\psi_2} (x_0) +M$.
This shows: $| V_k^{\psi_1} (x_0) - V_k^{\psi_2} (x_0)| \leq M$ for all $k$.
We may then divide by $k$ and let $k$ go to infinity to realize:
\[ | c_1 - c_2| = \lim_{k \rightarrow + \infty}  \frac{ | V_k^{\psi_1} (x_0) - V_k^{\psi_2} (x_0)| }{k} \leq \lim_{k \rightarrow + \infty} \frac{M}{k} = 0, \]
which completes our proof. \endproof

\subsection{Proof of formula (\ref{minformula})}
\label{formulaproof}
The formula is trivially fulfilled for $k=0$, remarking that by definition $\sum_{s \in \emptyset} (\cdot) = 0$.  
In fact:
\[ \min_{\tau \in \{0 \}} \left \{ T^\tau \psi + \min_{S = \emptyset} \sum_{s \in S} c( \hat{T}^s \psi,T \hat{T}^{s} \psi )  \right \} =
\psi = \hat{T}^0\psi. \]
Arguing by induction, and assuming the formula true for an arbitrary value of $k$, we can derive it for $k+1$ according to the following steps:
\[
\hat{T}^{k+1} \psi = \hat{T}^k ( \hat{T} \psi ) = \qquad \qquad \qquad \qquad \qquad \qquad \qquad \qquad \]
\[
= \min_{\tau \in \{0, \ldots, k\} } \left \{ T^{\tau} ( \hat{T} \psi ) +
\min_{S \subseteq \{ 0, \ldots, k-1 \}: |S|= \tau} \sum_{s \in S} c ( \hat{T}^s \hat{T} \psi, T \hat{T}^s \hat{T} \psi ) \right \} \qquad \qquad \]
\[ = \min_{\tau \in \{0, \ldots, k\} } \left \{ T^{\tau} ( \min \{ \psi, T \psi + c ( \psi, T \psi ) \} ) +
\min_{S \subseteq \{ 0, \ldots, k-1 \}: |S|= \tau} \sum_{s \in S}  c ( \hat{T}^{s+1} \psi, T \hat{T}^{s+1}  \psi ) \right \}
\]
\[ = \min_{\tau \in \{0, \ldots, k\} } \left \{ \min \{ T^{\tau} \psi, T^{\tau+1} \psi + c ( \psi, T \psi ) \}  +
\min_{S \subseteq \{ 1, \ldots, k \}: |S|= \tau} \sum_{s \in S} c ( \hat{T}^{s} \psi, T \hat{T}^{s}  \psi ) \right \}
\]
\[ = \min_{\tau \in \{0, \ldots, k\} } \left \{ \min \left \{ T^{\tau} \psi  + 
\min_{S \subseteq \{ 1, \ldots, k \}: |S|= \tau} \sum_{s \in S} c ( \hat{T}^{s} \psi, T \hat{T}^{s}  \psi ), \right . \right . \qquad \qquad  \]
\[ \qquad \qquad \qquad \left . \left . T^{\tau+1} \psi + c ( \psi, T \psi )   +
\min_{S \subseteq \{ 1, \ldots, k \}: |S|= \tau} \sum_{s \in S} c ( \hat{T}^{s} \psi, T \hat{T}^{s}  \psi ) \right \} \right \}
\]
\[ = \min \left \{ \min_{\tau \in \{0, \ldots, k\} } \left \{ T^{\tau} \psi  + 
\min_{S \subseteq \{ 1, \ldots, k \}: |S|= \tau} \sum_{s \in S} c ( \hat{T}^{s} \psi, T \hat{T}^{s}  \psi ) \right \}, \right . \qquad \qquad  \]
\[ \qquad \left .\min_{\tau \in \{ 0, \ldots, k \}} \left \{  T^{\tau+1} \psi + c ( \psi, T \psi )   +
\min_{S \subseteq \{ 1, \ldots, k \}: |S|= \tau} \sum_{s \in S} c ( \hat{T}^{s} \psi, T \hat{T}^{s}  \psi ) \right \} \right \}
\]
\[ = \min \left \{ \min_{\tau \in \{0, \ldots, k\} } \left \{ T^{\tau} \psi  + 
\min_{S \subseteq \{ 1, \ldots, k \}: |S|= \tau} \sum_{s \in S} c ( \hat{T}^{s} \psi, T \hat{T}^{s}  \psi ) \right \}, \right . \qquad \qquad  \]
\[ \qquad \left .\min_{\tau \in \{ 1, \ldots, k+1 \}} \left \{  T^{\tau} \psi + c ( \psi, T \psi )   +
\min_{S \subseteq \{ 1, \ldots, k \}: |S|= \tau-1} \sum_{s \in S} c ( \hat{T}^{s} \psi, T \hat{T}^{s}  \psi ) \right \} \right \}
\]
\[ = \min \left \{ \min_{\tau \in \{0, \ldots, k\} } \left \{ T^{\tau} \psi  + 
\min_{0 \notin S \subseteq \{0, 1, \ldots, k \}: |S|= \tau} \sum_{s \in S} c ( \hat{T}^{s} \psi, T \hat{T}^{s}  \psi ) \right \}, \right . \qquad \qquad  \]
\[ \qquad \quad \quad \left .\min_{\tau \in \{ 1, \ldots, k+1 \}} \left \{  T^{\tau} \psi + 
\min_{0 \in S \subseteq \{0, 1, \ldots, k \}: |S|= \tau} \sum_{s \in S} c ( \hat{T}^{s} \psi, T \hat{T}^{s}  \psi ) \right \} \right \}
\]
\[ = \min \left \{ \psi, \min_{\tau \in \{1, \ldots, k\} } \left \{ T^{\tau} \psi  + 
\min_{0 \notin S \subseteq \{0, 1, \ldots, k \}: |S|= \tau} \sum_{s \in S} c ( \hat{T}^{s} \psi, T \hat{T}^{s}  \psi ) \right \}, \right . \qquad \qquad  \]
\[ \qquad \quad \quad \left .\min_{\tau \in \{ 1, \ldots, k \}} \left \{  T^{\tau} \psi + 
\min_{0 \in S \subseteq \{0, 1, \ldots, k \}: |S|= \tau} \sum_{s \in S} c ( \hat{T}^{s} \psi, T \hat{T}^{s}  \psi ) \right \}, \right .\]
\[  \left . \qquad \qquad \qquad  
T^{k+1} \psi + \sum_{s =0}^{k+1}c ( \hat{T}^{s} \psi, T \hat{T}^{s}  \psi )   \right \}
\]
\[ = \min \left \{ \psi, \min_{\tau \in \{1, \ldots, k\} } \left \{ T^{\tau} \psi  + 
\min_{S \subseteq \{0, 1, \ldots, k \}: |S|= \tau} \sum_{s \in S} c ( \hat{T}^{s} \psi, T \hat{T}^{s}  \psi ) \right \}, \right . \qquad \qquad  \]
\[ \left .  \qquad \qquad \qquad  
T^{k+1} \psi + \sum_{s =0}^{k+1}c ( \hat{T}^{s} \psi, T \hat{T}^{s}  \psi  )  \right \}
\]
\[ = \min_{\tau \in \{0,1, \ldots, k+1\} } \left \{ T^{\tau} \psi  + 
\min_{S \subseteq \{0, 1, \ldots, k \}: |S|= \tau} \sum_{s \in S} c ( \hat{T}^{s} \psi, T \hat{T}^{s}  \psi ) \right \}.  \qquad \qquad \]
\endproof

\bibliographystyle{unsrt}
\bibliography{biblioinfinite}

\begin{thebibliography}{10}

\bibitem{bellman}
Richard Bellman.
\newblock The theory of dynamic programming.
\newblock {\em Bulletin of the American Mathematical Society}, 60:503–516,
  1954.

\bibitem{bertsekas1}
Dimitri~P. Bertsekas.
\newblock {\em Dynamic Programming and Optimal Control}, volume~I.
\newblock Athena Scientific, 4th edition, 2017.

\bibitem{bertsekas3}
Dimitri~P. Bertsekas.
\newblock {\em Abstract Dynamic Programming}.
\newblock Athena Scientific, 2nd edition, 2018.

\bibitem{infinite}
Haurie Dean A.~Carlson, Alain~B and Arie Leizarowitz.
\newblock {\em Infinite Horizon Optimal Control: Deterministic and Stochastic
  Systems}.
\newblock Springer, 1991.

\bibitem{stokey}
N.~Stokey and R.E. Lucas.
\newblock {\em Recursive Methods in Economic Dynamics}.
\newblock Harvard University Press, Cambridge, MA, 1989.

\bibitem{economic2}
Lars Ljungqvist and Thomas Sargent.
\newblock {\em Recursive Macroeconomic Theory}.
\newblock MIT Press, 3rd edition, 2012.

\bibitem{bertsekas2}
Dimitri~P. Bertsekas.
\newblock {\em Dynamic Programming and Optimal Control: Approximate Dynamic
  Programming}, volume~II.
\newblock Athena Scientific, 4th edition, 2012.

\bibitem{BarS18}
Andrew Barto and Richard~S. Sutton.
\newblock {\em Reinforcement Learning: An Introduction}.
\newblock MIT Press, 2nd edition, 2018.

\bibitem{ADBB17}
Kai Arulkumaran, Marc~Peter Deisenroth, Miles Brundage, and Anil~Anthony
  Bharath.
\newblock Deep reinforcement learning: A brief survey.
\newblock {\em IEEE Signal Processing Magazine}, 34(6):26--38, 2017.

\bibitem{willems1}
J.C. Willems.
\newblock Dissipative dynamical systems part i: General theory.
\newblock {\em Arch. Rational Mech. Anal.}, 45:321--–351, 1972.

\bibitem{willems2}
J.C. Willems.
\newblock Dissipative dynamical systems part ii: Linear systems with quadratic
  supply rates.
\newblock {\em Arch. Rational Mech. Anal.}, 45:352–--393, 1972.

\bibitem{moylan}
P.J. Moylan and B.D.O. Anderson.
\newblock Nonlinear regulator theory and an inverse optimal control problem.
\newblock {\em IEEE Transactions on Automatic Control}, 18:460–465, 1973.

\bibitem{dissipativity3}
Lars Gr\"une and Matthias~A. M\"uller.
\newblock On the relation between strict dissipativity and turnpike properties.
\newblock {\em Systems \& Control Letters}, 90:45--53, 2016.

\bibitem{dissipativity4}
Lars Gr\"une.
\newblock Dissipativity and optimal control: Examining the turnpike phenomenon.
\newblock {\em IEEE Control Systems Magazine}, 42(2):74--87, 2022.

\bibitem{dissipativity1}
Rishi~Amrit David~Angeli and James~B. Rawlings.
\newblock On average performance and stability of economic model predictive
  control.
\newblock {\em IEEE Transactions on Automatic Control}, 57(7):1615--1626, 2012.

\bibitem{dissipativity2}
David~Angeli Matthias A.~Müller and Frank Allg\"ower.
\newblock On necessity and robustness of dissipativity in economic model
  predictive control.
\newblock {\em IEEE Transactions on Automatic Control}, 60(6):1671--1676, 2015.

\bibitem{finlay}
Luke Finlay, Vladimir Gaitsgory, and Ivan Lebedev.
\newblock Duality in linear programming problems related to deterministic long
  run average problems of optimal control.
\newblock {\em SIAM Journal on Control and Optimization}, 47(4):1667--1700,
  2008.

\bibitem{gaitsgory0}
Vladimir Gaitsgory, Alex Parkinson, and Ilya Shvartsman.
\newblock Linear programming formulations of deterministic infinite horizon
  optimal control problems in discrete time.
\newblock {\em Discrete \& Continuous Dynamical Systems - B},
  22(10):3821--3838, 2017.

\bibitem{MULLERgrune}
Matthias~A. Müller and Lars Grüne.
\newblock Economic model predictive control without terminal constraints for
  optimal periodic behavior.
\newblock {\em Automatica}, 70:128--139, 2016.

\bibitem{muller}
M.~A. M\"{u}ller.
\newblock Dissipativity in economic model predictive control: beyond
  steady-state optimality.
\newblock In {\em Recent advances in model predictive control---theory,
  algorithms, and applications}, volume 485 of {\em Lect. Notes Control Inf.
  Sci.}, pages 27--43. Springer, Cham, [2021] \copyright 2021.

\bibitem{gaitsgory1}
Vivek~S. Borkar, Vladimir Gaitsgory, and Ilya Shvartsman.
\newblock Lp formulations of discrete time long-run average optimal control
  problems: The nonergodic case.
\newblock {\em SIAM Journal on Control and Optimization}, 57(3):1783--1817,
  2019.

\bibitem{AlBM07}
Olivier Alvarez, Martino Bardi, and Claudio Marchi.
\newblock Multiscale problems and homogenization for second-order
  {H}amilton-{J}acobi equations.
\newblock {\em J. Differential Equations}, 243(2):349--387, 2007.

\bibitem{GGHKW18}
Vladimir Gaitsgory, Lars Gr\"une, Matthias H\"oger, Christopher~M. Kellett, and
  Steven~R. Weller.
\newblock Stabilization of strictly dissipative discrete time systems with
  discounted optimal control.
\newblock {\em Automatica}, 93:311--320, 2018.

\bibitem{GruK21}
Lars Gr\"{u}ne and Lisa Kr\"{u}gel.
\newblock Local turnpike analysis using local dissipativity for discrete time
  discounted optimal control.
\newblock {\em Appl. Math. Optim.}, 84(suppl. 2):S1585--S1606, 2021.

\bibitem{borkar}
Vivek~S. Borkar and Vladimir Gaitsgory.
\newblock Linear programming formulation of long run average optimal control
  problem.
\newblock {\em Journal of Optimization Theory and Applications}, 181:101--125,
  2019.

\end{thebibliography}

\end{document}